\documentclass[pdflatex,sn-mathphys-num]{sn-jnl}

\usepackage{amsmath,amssymb,amsfonts}
\usepackage{mathtools}
\usepackage{booktabs}
\usepackage[expansion=false]{microtype}
\graphicspath{{figures/}{paperA_v3/figures/}{./}}

\theoremstyle{thmstyleone}
\newtheorem{theorem}{Theorem}[section]
\newtheorem{proposition}[theorem]{Proposition}
\newtheorem{lemma}[theorem]{Lemma}
\newtheorem{corollary}[theorem]{Corollary}
\theoremstyle{thmstylethree}
\newtheorem{definition}[theorem]{Definition}
\theoremstyle{thmstyletwo}
\newtheorem{remk}[theorem]{Remark}
\newtheorem*{openproblem}{Open Problem}

\newcommand{\pstep}[1]{\smallskip\noindent\textit{Step #1.}\ }
\newcommand{\ppart}[1]{\smallskip\noindent\textit{#1}\ }

\newcommand{\RR}{\mathbb{R}}
\newcommand{\EE}{\mathbb{E}}
\newcommand{\PP}{\mathbb{P}}
\newcommand{\Sph}[1]{S^{#1}}
\DeclareMathOperator{\rank}{rank}
\DeclareMathOperator{\im}{im}
\DeclareMathOperator{\spn}{span}
\DeclareMathOperator{\cone}{cone}
\DeclareMathOperator{\tr}{tr}
\DeclareMathOperator{\diag}{diag}
\DeclareMathOperator{\sgn}{sgn}
\DeclareMathOperator{\vol}{vol}
\DeclareMathOperator{\Cov}{Cov}
\DeclareMathOperator{\Haf}{Haf}
\DeclareMathOperator{\supp}{supp}
\DeclareMathOperator{\ind}{ind}
\DeclareMathOperator{\Lk}{Lk}
\DeclareMathOperator{\dist}{dist}
\DeclareMathOperator{\rec}{rec}
\newcommand{\transp}{^{\!\top}}
\newcommand{\pinv}{^{\dagger}}
\newcommand{\Norm}{\mathcal{N}}
\newcommand{\Haus}[1]{\mathcal{H}^{#1}}
\newcommand{\ip}[2]{\langle #1,\, #2\rangle}
\newcommand{\norm}[2][]{{#1\lVert #2 #1\rVert}}
\newcommand{\abs}[2][]{{#1\lvert #2 #1\rvert}}
\newcommand{\dd}{\mathop{}\!d}

\newcommand{\dimin}{D}                  % domain dimension
\newcommand{\ngen}{H}                   % number of walls
\newcommand{\wall}[1]{\Pi_{#1}}         % wall h
\newcommand{\arr}{\mathcal{A}}          % the arrangement
\newcommand{\cell}[1]{\mathcal{C}_{#1}} % open chamber

\newcommand{\Match}[1]{\mathfrak{M}_{#1}}
\newcommand{\Spec}{\mathbb{S}}
\newcommand{\Pprof}{\mathcal{P}}
\newcommand{\Rtube}{\mathsf{Q}}
\newcommand{\face}{\mathfrak{c}}
\newcommand{\psector}{\mathsf{S}}       % planar sector of a central fan

\AtBeginDocument{\emergencystretch=3em\relax\tolerance=1000\relax}
\allowdisplaybreaks

\newcommand{\PaperAppendices}{\begin{appendices}\counterwithin*{equation}{section}}
\newcommand{\PaperAppendicesEnd}{\end{appendices}}

\begin{document}

\title[Banchoff Curvature and Realisation Moduli of Arrangement Immersions]{Banchoff Curvature and Realisation Moduli of Piecewise-Affine Arrangement Immersions}

\author*[1]{\fnm{William} \sur{Nelson}}\email{wtnelson@google.com}
\affil*[1]{\orgname{Google DeepMind}, \orgaddress{\city{New York}, \country{USA}}}

% \unboldmath: sn-jnl's \Authorfont sets \boldmath, which leaks into the
% abstract and would set its inline maths in bold.
\abstract{\unboldmath We compute the curvature concentrated on the flats of a
piecewise-affine immersion of $\RR^{\dimin}$ into $\RR^n$ creased along an
affine hyperplane arrangement, and parametrise the immersions sharing its
chamber metrics. Suppose every transverse crease cone is embedded (the
positive-path-rank condition). Then the Gauss--Bonnet--Chern density of any
convex coordinate-wise smoothing that is affine outside a compact transition
zone reduces, on each transverse slice, to a Hafnian. This Hafnian integrates
to Banchoff's critical-point curvature on each even-codimension flat. The
curvature coincides with the Cheeger--M\"uller--Schrader density and equals an
alternating sum of Gaussian orthant probabilities determined by Schur
complements of the incident chamber metrics. In even dimension, under mild
non-degeneracy, total vertex curvature equals the recession curvature at
infinity and is unchanged when the walls are translated, even across
multi-wall collisions. Under generic codimension-two intersections,
realisations with prescribed chamber metrics form a compact rank-constrained
spectrahedron of generic local dimension $\max(0,\ngen-r(r+1)/2)$, where
$r = \max(0,\dimin+\ngen-n)$. At minimal target dimension the given immersion
is the unique realisation preserving orientation across each single wall, and
a positive semidefinite stress matrix certifies universal rigidity.}



\maketitle

\section{Introduction}\label{sec:intro}

\subsection{Setting and geometric context}\label{subsec:setting}

A smooth submanifold of Euclidean space carries pointwise Riemannian
curvature. A continuous piecewise-affine map creased along an affine hyperplane
arrangement, by contrast, is flat on every open chamber and concentrates all of
its curvature on the flats where the walls meet. Maps of this kind arise in
polyhedral differential geometry and rigid origami~\cite{banchoff1967,
demaineorourke2007}, and as multivariate continuous degree-$1$
splines~\cite{billera1988, ziegler1995}. We ask how the
Gauss--Bonnet--Chern curvature of a smooth regularisation concentrates onto the
flats of the arrangement, and to what extent the induced chamber metrics
determine the ambient immersion.

To be concrete, let $\arr = \{\wall{1},\dots,\wall{\ngen}\}$ be a finite
arrangement of pairwise distinct affine hyperplanes
$\wall{h} = \{x \in \RR^{\dimin} : w_h\transp x + b_h = 0\}$ in $\RR^{\dimin}$
($\dimin \ge 2$, cf.\ Remark~\ref{rem:PR-dichotomy}), with non-zero normals
$w_h$ and offset vector
$b = (b_1,\dots,b_{\ngen})\transp$. We study continuous piecewise-affine maps
$F : \RR^{\dimin} \to \RR^n$ ($n \ge \dimin$) of the form
\begin{equation}\label{eq:intro-map}
  F(x) = Lx + c_0 + \sum_{h=1}^{\ngen} a_h\,(w_h\transp x + b_h)_+,
\end{equation}
where $(z)_+ = \max(z,0)$, $L \in \RR^{n\times\dimin}$, $c_0 \in \RR^n$, and
$A = [a_1\ \cdots\ a_{\ngen}] \in \RR^{n\times\ngen}$. On each chamber of
$\arr$, indexed by its sign pattern $\sigma \in \{0,1\}^{\ngen}$, the map $F$ is
affine with differential $J_\sigma = L + \sum_h \sigma_h a_h w_h\transp$ and
constant pullback metric $g_\sigma = J_\sigma\transp J_\sigma$. Across each wall
$\wall{h}$, the differential $J_\sigma$ jumps by the rank-one matrix
$a_h w_h\transp$, so $g_\sigma$ jumps by a symmetric matrix of rank at most two.
When every chamber Jacobian $J_\sigma$ has rank $\dimin$ and $F$ satisfies the
positive-path-rank condition of Section~\ref{sec:hafnian}, its image is a
piecewise-flat immersed submanifold of $\RR^n$. (Positive path rank embeds
every transverse crease cone and excludes wall fold-overs.) Call $\arr$
codimension-$2$ generic if no two walls are parallel and no three walls share
a codimension-$2$ flat. In that case the ansatz~\eqref{eq:intro-map} exhausts
all continuous piecewise-affine maps subordinate to $\arr$
(Corollary~\ref{cor:spline-exhaustion}, after~\cite{billera1988,
billerarose1991}). In approximation theory and machine
learning,~\eqref{eq:intro-map} is the hinging-hyperplane or single-layer
rectified-linear (ReLU) representation with a linear skip
connection~\cite{breiman1993, wangsun2005}. Its smoothings replace
$(z)_+$ by softplus-type profiles. The profile class $\Pprof$ used below
(Definition~\ref{def:arr-immersion}) instead confines the transition to a
compact interval. This keeps the transition tubes of different walls apart at
a finite scale in Section~\ref{sec:global}. The local flatwise identity of
Theorem~\ref{thm:banchoff-orthant}, on the other hand, holds for softplus as
well (Remark~\ref{rem:profile-scope}). The companion
paper~\cite{nelsonpaperb2026} analyses the stochastic geometry of such maps for
random arrangements. Here we treat their deterministic polyhedral curvature and
isometric realisation moduli.

\subsection{Main results}\label{subsec:main-results}

Replacing $(z)_+$ by a convex smoothing $\rho_\beta(z) = \beta^{-1}\rho(\beta z)$
turns $F$ into a smooth map $F_\beta$. We work with its pullback metric
$g_\beta$ on $\RR^{\dimin}$, which avoids any ambiguity from ambient
self-intersections. At first sight the smoothed metric $g_\beta$ across a wall
tube depends on the full Gram matrix of $[L\ A]$, including cross-terms
invisible inside any open chamber. Yet its Gauss--Bonnet--Chern integral
across any even-codimension flat $\Phi_I = \bigcap_{h\in I}\wall{h}$ depends
solely on the incident chamber metrics (Remark~\ref{rem:bridge}). We prove
three main results.

\emph{Curvature at a flat.} Suppose positive path rank holds at an
even-codimension flat $\Phi_I$ ($|I| = 2k$). A unimodular nonlinear shear then
splits $g_\beta$ near $\Phi_I$ into a flat tangential block and a
$2k$-dimensional translation submanifold. This transverse slice is additively
separable, so its second fundamental form is coordinate-diagonal and the Chern
Pfaffian reduces to a Hafnian density $h_{2k}(s)$ on $[0,1]^{2k}$, independent
of $\rho$ and $\beta$ (Proposition~\ref{prop:hafnian-euler}). Integrating
$h_{2k}$ recovers Banchoff's critical-point curvature $K(\Phi_I)$ as the
alternating Gaussian orthant sum (Theorem~\ref{thm:banchoff-orthant})
\[
  K(\Phi_I)
  = \sum_{\tau\in\{\pm\}^{2k}} (-1)^{\#\tau_+}\,\PP(Y^\tau > 0),
  \qquad Y^\tau \sim \Norm(0, M_\tau),\quad
  M_\tau = (W_I g_\tau^{-1} W_I\transp)^{-1},
\]
which depends only on the Schur complements $M_\tau$ of the incident chamber
metrics and coincides with the Cheeger--M\"uller--Schrader curvature density
along $\Phi_I$ (Proposition~\ref{prop:schur-cms}).

\emph{Global vertex curvature.} When the arrangement is simple and positive
path rank holds at every face, the smoothed Gauss--Bonnet--Chern measures
converge weak-$*$ to the polyhedral curvature measure at local rate
$O(\beta^{-1})$ on compact sets (Theorem~\ref{thm:global-weak-star}). In top
degree $\dimin = 2k$ the convergence is global, and each vertex tube carries
exact mass $K(v)$ already at finite $\beta$. A cellular Morse count at infinity
then shows that $\sum_v K(v) = K_\infty$, the apex curvature of the recession
map $x \mapsto Lx + A(Wx)_+$, for every offset $b$ that contracts no bounded
edge to a point (Theorem~\ref{thm:offset-invariance}). So translating the walls
leaves the total vertex curvature unchanged, even across collisions of several
walls. For convex graphs this reproduces the subdifferential tiling of the
Newton zonotope (Corollary~\ref{cor:conic-zonotope}).

\emph{Realisation moduli.} When $\arr$ is codimension-$2$ generic, the Gram
matrices $G = [L\ A]\transp[L\ A]$ inducing prescribed chamber metrics
$\{g_\sigma\}$ form an affine slice $G_0 + \spn\{E_h\}$ spanned by rank-two
wall matrices $E_h$ (Theorem~\ref{thm:gram-spectrahedron}). Isometric
realisations in $\RR^n$ modulo $O(n)$ therefore form a compact rank-constrained
spectrahedron $\mathcal{M}(G_0,n)$ of generic local dimension
$\max(0, \ngen - \binom{r+1}{2})$, where $r = \max(0,\dimin+\ngen-n)$
(Theorems~\ref{thm:generic-local-dim} and~\ref{thm:compactness-stress}). At
minimal target dimension $n = \dimin$, the given immersion is the unique
realisation that preserves orientation across each single wall
(Theorem~\ref{thm:compactness-stress}(2)). In higher codimension, by
contrast, an arrangement with a missing chamber can carry a second realisation
that is immersive on every face it actually has
(Proposition~\ref{prop:vertex-flex-rigid}(4)). For $n < \dimin+\ngen$, a
positive semidefinite stress matrix certifies universal rigidity across all
continuous piecewise-affine maps subordinate to $\arr$ in every ambient
dimension $n' \ge n$ (Corollary~\ref{cor:stress-rigidity}).

Sections~\ref{sec:hafnian}--\ref{sec:moduli} prove these results in turn. The
main tools are the transverse shear splitting, the Hafnian evaluation of the
Kac--Rice density and the rank-two wall pencil developed below, used together
with classical degree-$1$ spline exhaustion
(Corollary~\ref{cor:spline-exhaustion}), polyhedral Morse theory on the
recession fan (Lemma~\ref{lem:pl-morse-infinity}) and the Schur-complement
determinantal chart (Lemma~\ref{lem:psd-schur-chart}).
Section~\ref{sec:conclusion} discusses open problems. Auxiliary bounds, global
collar estimates and algebraic genericity proofs are deferred to
Appendices~\ref{app:banchoff-aux}--\ref{app:moduli-proofs}.

\subsection{Related work}\label{subsec:related-work}

Polyhedral curvature is classically defined extrinsically, by averaging
critical-point indices of linear height functions~\cite{banchoff1967,
banchoff1983}, and intrinsically, via cone-angle Lipschitz--Killing
measures~\cite{cheeger1984}. For smooth immersions, the Gauss--Bonnet--Chern
form~\cite{chern1944} is the Gaussian average over normal directions of the
Hessian determinant of the height function~\cite{weyl1939, chernlashof1957},
so the Kac--Rice formula~\cite{adlertaylor2007} expresses its integral as an
expected signed count of critical points. These are the classical ingredients
of Section~\ref{sec:banchoff}. What is specific to arrangement immersions is
the mechanism that localises the integral onto the flats. Once the tangential
wall jumps are sheared away, the transverse slice is a translation
submanifold. Its second fundamental form is then coordinate-diagonal and its
Chern Pfaffian reduces to a Hafnian. Gaussian regression together with
Isserlis' theorem~\cite{isserlis1918} identifies this Hafnian pointwise with
the Kac--Rice density (Proposition~\ref{prop:hafnian-euler} and
Theorem~\ref{thm:banchoff-orthant}).

Coordinate-wise smoothing is not merely a convenient choice here. An
arrangement immersion is not smooth outside any compact set, since its wall
creases and lower-dimensional flats run to infinity. A smoothing must
therefore regularise every crease transversely along its whole length, and
along a positive-dimensional flat ($2k < \dimin$) this couples the tangential
and transverse directions. Additive separability of the profile lets the
unimodular shear of Lemma~\ref{lem:shear-splitting} decouple the flat
tangential factor from the $2k$-dimensional transverse slice. It also
produces a pointwise density $h_{2k}(s)$ on $[0,1]^{2k}$, positive for convex
graphs (Corollary~\ref{cor:conic-zonotope}). The integrated curvature itself is
Banchoff's and does not depend on the profile.

Polyhedral surfaces admit topology-preserving smoothings that concentrate
Gaussian curvature near the vertices~\cite{brehmkuhnel1982}. On compact
ambient subsets, normal-cycle theory bounds curvature discrepancies under
polyhedral or smooth approximation~\cite{fu1994, cohensteinermorvan2006}.
Working instead on $(\RR^{\dimin}, g_\beta)$ avoids
ambient self-intersection issues and yields the local $O(\beta^{-1})$ rate
directly from positive path rank. Globally, vertex indices on a closed
polyhedral manifold sum to its Euler characteristic~\cite{banchoff1967,
banchoffkuhnel1997, grunertkuhnelrote2023}, but an affine arrangement is
unbounded. Its face and chamber counts are
classical~\cite{zaslavsky1975, orlikterao1992}. More recently, transversality
of bent hyperplane arrangements~\cite{grigsbylindsey2022} and the
piecewise-linear Morse theory of rectified-linear
maps~\cite{grigsbylindseymasden2022} have been developed; see
also~\cite{zhangnaitzatlim2018, haninrolnick2019} on tropical duality and
linear regions. Using the Euler
characteristic valuation on polyhedra~\cite{mcmullen1975, schneiderweil2008},
we show that translating the walls (even across multi-wall circuit collisions)
preserves total vertex curvature. In two dimensions this gives a polyhedral
counterpart of the total-curvature theorems of Cohn-Vossen and
Huber~\cite{cohnvossen1935, huber1957}. For convex graphs it recovers the
subdifferential tiling of the Newton zonotope, a regular zonotopal tiling in
the sense of Bohne and Dress~\cite{richtergebertziegler1994},
\cite[\S7.5]{ziegler1995}.

On the algebraic side, continuous degree-$1$ splines on an arrangement are
governed by codimension-$2$ syzygies~\cite{billera1988}, a perspective recently
extended to decomposition polyhedra of piecewise-linear
functions~\cite{brandenburggrillohertrich2026}. Passing to chamber metrics turns
the rank-one wall jumps into a pencil of rank-two symmetric matrices. This
gives a piecewise-linear counterpart of the Gram spectrahedra of convex
algebraic geometry~\cite{chuaplaumannsinnvinzant2017}, \cite{blekherman2013}.
Its generic local dimension follows the rank bounds of Barvinok and
Pataki~\cite{barvinok1995, pataki1998}. Our universal-rigidity criterion
parallels the stress-matrix theory of Connelly, Gortler, and
Thurston~\cite{connelly1982, gortlerthurston2014}, in the Gale-matrix form of
Alfakih~\cite{alfakih2007, alfakih2018}. Unlike isometric piecewise-linear
embeddings that are free to re-triangulate the
domain~\cite{buragozalgaller1996}, our realisations keep every crease wall
fixed. General rigid-origami crease patterns satisfy loopwise affine
compatibility around vertices~\cite{belcastrohull2002, demaineorourke2007}.
In our setting each wall is a single flat jump across all the chambers it
separates. This couples all vertices at once, in contrast to a single
Kokotsakis vertex~\cite{kokotsakis1933, izmestiev2017}.
\section{Arrangement Immersions, Shear Splitting, and the Hafnian Density}\label{sec:hafnian}

\subsection{Arrangement immersions, local models, and standing conditions}\label{subsec:arr-local-models}

Throughout, $\Spec^d \supset \Spec_+^d \supset \Spec_{++}^d$ denote the space
of real $d\times d$ symmetric matrices and its positive semidefinite and
positive definite cones.

\begin{definition}
\label{def:arr-immersion}
Fix integers $\dimin \ge 1$, $\ngen \ge 1$, and $n \ge \dimin$, and set
$m := n - \dimin \ge 0$. Let $L \in \RR^{n\times\dimin}$, $c_0 \in \RR^n$,
$A = [a_1\ \cdots\ a_{\ngen}] \in \RR^{n\times\ngen}$,
$W = [w_1\ \cdots\ w_{\ngen}]\transp \in \RR^{\ngen\times\dimin}$ with
$w_h \ne 0$ ($1 \le h \le \ngen$), and $b \in \RR^{\ngen}$, and write
$z(x) := Wx + b$ with components $z_h(x) := w_h\transp x + b_h$.
\begin{enumerate}
  \item The \emph{piecewise-affine map} of $(L, c_0, A, W, b)$ is
    \begin{equation}\label{eq:arr-map}
      F(x) := Lx + c_0 + A(Wx+b)_+
      = Lx + c_0 + \sum_{h=1}^{\ngen} a_h\,(z_h(x))_+.
    \end{equation}
    Without loss of generality, we assume throughout that the walls
    $\wall{h} := \{x \in \RR^{\dimin} : z_h(x) = 0\}$ ($1 \le h \le \ngen$)
    are pairwise distinct affine hyperplanes forming the arrangement
    $\arr = \{\wall{1},\dots,\wall{\ngen}\}$.
    Co-oriented duplicate walls merge after adjusting $a_h$. Anti-oriented
    duplicates $(-w_h,-b_h)$ merge after adjusting $(L,c_0,a_h)$ via
    $(-z)_+ = (z)_+ - z$, or are absorbed into $(L,c_0)$ if the jumps
    cancel. Each non-empty open
    cell $\cell{\sigma} := \{x \in \RR^{\dimin} : (2\sigma_h - 1)z_h(x) > 0\
    \forall h\}$ ($\sigma \in \{0,1\}^{\ngen}$) is a \emph{chamber} of $\arr$,
    on which $F$ is affine, with differential and pullback metric given by
    \begin{equation}\label{eq:chamber-jacobian}
      J_\sigma := L + \sum_{h=1}^{\ngen} \sigma_h a_h w_h\transp
      \in \RR^{n\times\dimin},
      \qquad
      g_\sigma := J_\sigma\transp J_\sigma \in \Spec^{\dimin}.
    \end{equation}
    We call $F$ an \emph{arrangement immersion} if $\rank J_\sigma = \dimin$ on
    every non-empty chamber $\cell{\sigma}$.

  \item Let $\Pprof$ denote the class of \emph{smoothing profiles}
    $\rho \in C^\infty(\RR)$ satisfying
    \[
      \rho(z) = 0\ (z \le -1),\qquad
      \rho''(z) > 0\ (z \in (-1,1)),\qquad
      \rho'(z) = 1\ (z \ge 1),
    \]
    so $\supp\rho'' = [-1,1]$, $\rho'$ maps $(-1,1)$ diffeomorphically onto
    $(0,1)$, $\int_{-1}^1\rho''(z)\,\dd z = 1$, and $\rho(z) = z - c_\rho$ for
    $z \ge 1$ with $c_\rho := \int_{-1}^1 z\rho''(z)\,\dd z \in (-1,1)$. Note
    that there is no assumption about the symmetry of $\rho$. For $\beta > 0$,
    set
    $\rho_\beta(z) := \beta^{-1}\rho(\beta z)$. Since $0 \le \rho(z)
    \le z + 1$ on $[-1,1]$ and $|c_\rho| < 1$, we have
    $\sup_{z\in\RR}|\rho_\beta(z) - z_+| \le \beta^{-1}$. Define the
    \emph{smoothed map} $F_\beta$ and pullback metric
    $g_\beta(x) := DF_\beta(x)\transp DF_\beta(x)$ by
    \[
      F_\beta(x) := Lx + c_0 + \sum_{h=1}^{\ngen} a_h\,\rho_\beta(z_h(x)),\qquad
      DF_\beta(x) = L + \sum_{h=1}^{\ngen} \rho'(\beta z_h(x))\,a_h w_h\transp.
    \]
\end{enumerate}
\end{definition}

\begin{definition}
\label{def:local-model}
Fix an ordered $q$-tuple $I = (h_1,\dots,h_q)$ ($1 \le q \le \dimin$) and a
point $x_0 \in \Phi_I := \bigcap_{a=1}^q \wall{h_a}$ at which \emph{local
transversality} holds: $W_I := [w_{h_1}\ \cdots\ w_{h_q}]\transp \in
\RR^{q\times\dimin}$ has rank $q$ and $z_h(x_0) \ne 0$ for all $h \notin I$.
\begin{enumerate}
  \item Set $T\Phi_I := \ker W_I$, $\sigma^0_h := \mathbf{1}\{z_h(x_0) > 0\}$
    ($h \notin I$), and define the \emph{base Jacobian} and \emph{wall dual
    basis}
    \begin{gather*}
      J_0 := L + \sum_{h\notin I} \sigma^0_h a_h w_h\transp
      \in \RR^{n\times\dimin},\\
      V_I := W_I\transp(W_I W_I\transp)^{-1} = W_I\pinv = [v_1\ \cdots\ v_q]
      \in \RR^{\dimin\times q},
    \end{gather*}
    so $W_I V_I = I_q$ and $V_I\transp T\Phi_I = 0$. Fix an orthonormal basis
    $Q_\perp \in \RR^{\dimin\times(\dimin-q)}$ of $T\Phi_I$, and define the
    adapted affine chart $\iota_I : \RR^q\times\RR^{\dimin-q} \to \RR^{\dimin}$
    by $\iota_I(u,x_\perp) := x_0 + V_I u + Q_\perp x_\perp$, so $u_a =
    z_{h_a}(\iota_I(u,x_\perp))$. The \emph{transition strip} is
    $\Rtube_\beta := [-\beta^{-1},\beta^{-1}]^q$, written $\Rtube_\beta(I)$
    when the flat is emphasised. A set $U \subset \RR^{\dimin}$ is
    \emph{$\beta$-admissible at $I$} if $(2\sigma^0_h - 1)\beta z_h(x) \ge 1$
    for all $x \in U$ and $h \notin I$. Since $z_h(x_0) \ne 0$ for every
    $h \notin I$, an open ball $U_0 \ni x_0$ is $\beta$-admissible for all
    large $\beta$, with $U_0 := \RR^{\dimin}$ when $I = \{1,\dots,\ngen\}$.
    For any $\beta$-admissible $U \subset \RR^{\dimin}$ and $(u,x_\perp) \in
    \widehat{U} := \iota_I^{-1}(U)$,
    \begin{equation}\label{eq:local-model}
      F_\beta(\iota_I(u,x_\perp))
      = c_\beta + J_0 Q_\perp x_\perp
        + \sum_{a=1}^q \bigl(u_a J_0 v_a + \rho_\beta(u_a)\,a_{h_a}\bigr),
    \end{equation}
    where $c_\beta := L x_0 + c_0 + \sum_{h\notin I}\sigma^0_h(z_h(x_0) -
    \beta^{-1}c_\rho)a_h$.

  \item When $\rank(J_0 Q_\perp) = \dimin - q$, set $\mathcal{T}_I :=
    \im(J_0 Q_\perp) \subset \RR^n$ ($\dim\mathcal{T}_I = \dimin - q$),
    $P_{\mathcal{T}_I} := (J_0 Q_\perp)((J_0 Q_\perp)\transp(J_0 Q_\perp))^{-1}
    (J_0 Q_\perp)\transp$, and $P_I := I_n - P_{\mathcal{T}_I}$ onto
    $\mathcal{T}_I^\perp$ ($\dim\mathcal{T}_I^\perp = q + m$). Put
    $r_a := P_I J_0 v_a$, $p_a := P_I a_{h_a}$, $X_a(s_a) := r_a + s_a p_a$, and
    $X(s) := [X_1(s_1)\ \cdots\ X_q(s_q)] \in \RR^{n\times q}$ ($s \in [0,1]^q$),
    with transverse Gram matrix and normal projector
    \begin{equation}\label{eq:schur-projector}
      M(s) := X(s)\transp X(s) \in \Spec^q,\qquad
      \Pi^\perp(s) := P_I - X(s) M(s)^{-1} X(s)\transp.
    \end{equation}
\end{enumerate}
\end{definition}

Besides local transversality, Sections~\ref{sec:hafnian}--\ref{sec:global} use
two regularity conditions:
\begin{itemize}
  \item \emph{Positive path rank at $I$:}
    $\rank(J_0 + \sum_{a=1}^q s_a a_{h_a} w_{h_a}\transp) = \dimin$ on
    $[0,1]^q$, equivalently $\rank(J_0 Q_\perp) = \dimin - q$ and $M(s) \succ 0$
    on $[0,1]^q$. Then $\Pi^\perp(s)$ projects orthogonally onto
    $\mathcal{T}_I^\perp \cap (\im X(s))^\perp$ ($\rank\Pi^\perp(s) = m$), and
    we set $\lambda_0 := \min_{s\in[0,1]^q}\lambda_{\min}(M(s)) > 0$ and
    $\overline{\lambda}_0 := \max_{s\in[0,1]^q}\lambda_{\max}(M(s))$.
  \item \emph{Global genericity:}
    $\arr(b)$ is \emph{simple}, meaning that every non-empty $J \subset
    \{1,\dots,\ngen\}$ with $\Phi_J := \bigcap_{h\in J}\wall{h} \ne \emptyset$
    satisfies $|J| \le \dimin$ and $\rank W_J = |J|$. In addition, \emph{positive
    path rank holds at every face}: every non-empty chamber has
    $\rank J_\sigma = \dimin$, and for every non-empty face $\mathsf{F} \in
    \mathcal{F}(\arr(b))$ with active set $I(\mathsf{F}) := \{h : \mathsf{F}
    \subset \wall{h}\}$, the interpolated Jacobian
    $J_{\mathsf{F}}(s) := L + \sum_{h\notin I(\mathsf{F})}\mathbf{1}\{z_h > 0
    \text{ on }\mathsf{F}\}a_h w_h\transp + \sum_{h\in I(\mathsf{F})}s_h a_h
    w_h\transp$ has rank $\dimin$ for all $s \in [0,1]^{I(\mathsf{F})}$.
\end{itemize}
Here $\mathcal{F}(\arr(b))$ denotes the set of non-empty relatively open faces
of $\arr(b)$, $\mathsf{F}(y)$ the face containing $y \in \RR^{\dimin}$, and
$V(b)$ the vertex set of $\arr(b)$, with $I(v) := I(\{v\})$.

\begin{remk}\label{rem:scope-PR-GA}
Positive path rank at $I$ is equivalent to global injectivity of the transverse
piecewise-linear cone $\Psi_0(u) := \sum_{a=1}^q (u_a r_a + (u_a)_+ p_a) :
\RR^q \to \mathcal{T}_I^\perp$ (Lemma~\ref{lem:pr-injectivity}(1)), so it rules
out wall fold-overs (Remark~\ref{rem:PR-dichotomy}). It holds automatically for
graph lifts $x \mapsto (x, A'(Wx+b)_+)$, and, whenever $\rank L = \dimin$, for
all sufficiently small $\norm{A}$ (depending on $L$ and $W$). The
set $\{X(s) : s \in [0,1]^q\} = \operatorname{conv}\{X_\tau : \tau \in
\{-,+\}^q\}$ is Clarke's generalised Jacobian of $\Psi_0$ at
$0$~\cite{clarke1983}. Positive path rank is therefore the familiar
maximal-rank condition on the generalised Jacobian from nonsmooth analysis; in
the square case $m = 0$ it is coherent orientation~\cite{scholtes2012}. Here
the condition is also sufficient for global injectivity, because additive
separability realises every $X(\bar s)$ as an exact secant operator of
$\Psi_0$. Simplicity in global genericity excludes central arrangements with
$\ngen > \dimin$ (Lemma~\ref{lem:joint-genericity}(4)). Neither
Sections~\ref{sec:hafnian}--\ref{sec:banchoff} nor
Theorem~\ref{thm:offset-invariance} assumes simplicity.
\end{remk}

\subsection{Unimodular shear splitting and the Hafnian Euler density}\label{subsec:shear-hafnian}

\begin{figure}[t]
\centering
\includegraphics[width=\textwidth]{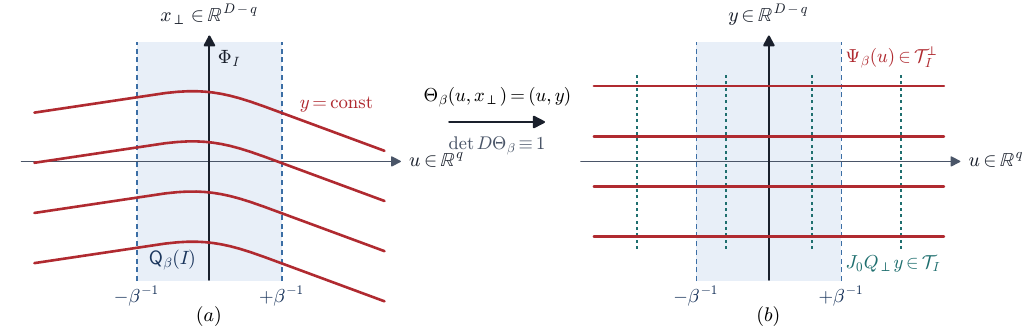}
\caption{Unimodular shear $\Theta_\beta$ near $\Phi_I$ (Lemma~\ref{lem:shear-splitting}). \emph{(a)}~When $P_{\mathcal{T}_I}a_{h_a} \ne 0$, the tangential part of $F_\beta$ bends across $\Rtube_\beta(I)$. \emph{(b)}~Setting $y = x_\perp + \sum_a \varphi_{a,\beta}(u_a)$ straightens the fibres with $\det D\Theta_\beta = 1$ and splits the image into $\mathcal{T}_I \oplus \mathcal{T}_I^\perp$.}
\label{fig:shear}
\end{figure}

When $P_{\mathcal{T}_I}a_{h_a} \ne 0$, the derivative $\partial_{u_a}F_\beta =
J_0 v_a + \rho'(\beta u_a)a_{h_a}$ has a $\mathcal{T}_I$-component that varies
with $u_a$. No linear change of coordinates can then orthogonalise $g_\beta$.
A separable shear with unit Jacobian does the job instead: on any
$\beta$-admissible region $U$ it splits $g_\beta$ into a Riemannian product
(Figure~\ref{fig:shear}).

\begin{lemma}
\label{lem:shear-splitting}
In Definitions~\ref{def:arr-immersion}--\ref{def:local-model}, assume local
transversality and $\rank(J_0 Q_\perp) = \dimin-q$, and set $(J_0 Q_\perp)\pinv
:= ((J_0 Q_\perp)\transp(J_0 Q_\perp))^{-1}(J_0 Q_\perp)\transp$. Define
$\varphi_{a,\beta} : \RR \to \RR^{\dimin-q}$, $\Psi_\beta : \RR^q \to
\mathcal{T}_I^\perp$, and $\Theta_\beta : \RR^q\times\RR^{\dimin-q} \to
\RR^q\times\RR^{\dimin-q}$ by $\Theta_\beta(u,x_\perp) := (u, y)$ with
$y := x_\perp + \sum_{a=1}^q\varphi_{a,\beta}(u_a)$ and
\[
  \varphi_{a,\beta}(u_a)
  := (J_0 Q_\perp)\pinv\bigl(u_a J_0 v_a + \rho_\beta(u_a)\,a_{h_a}\bigr),
  \qquad
  \Psi_\beta(u)
  := \sum_{a=1}^q \bigl(u_a r_a + \rho_\beta(u_a)\,p_a\bigr).
\]
Then for any $\beta$-admissible region $U \subset \RR^{\dimin}$ with coordinate
preimage $\widehat{U} := \iota_I^{-1}(U) \subset \RR^q\times\RR^{\dimin-q}$:
\begin{enumerate}
  \item $\Theta_\beta$ is a smooth diffeomorphism that preserves every
    transverse cylinder $\Omega\times\RR^{\dimin-q}$ ($\Omega \subset \RR^q$),
    with $D\Theta_\beta = \bigl(\begin{smallmatrix} I_q & 0 \\ * & I_{\dimin-q}
    \end{smallmatrix}\bigr)$, $\det D\Theta_\beta \equiv 1$, and inverse
    $\Theta_\beta^{-1}(u,y) = \bigl(u,\, y - \sum_{a=1}^q\varphi_{a,\beta}(u_a)\bigr)$.
  \item $\widetilde{F}_\beta := F_\beta \circ \iota_I \circ \Theta_\beta^{-1}$ on
    $\Theta_\beta(\widehat{U})$ satisfies
    $\widetilde{F}_\beta(\Theta_\beta(\widehat{U})) = F_\beta(U)$ and splits
    orthogonally into $\mathcal{T}_I \oplus \mathcal{T}_I^\perp$:
    \begin{equation}\label{eq:metric-split}
      \widetilde{F}_\beta(u, y)
      = c_\beta + J_0 Q_\perp y + \Psi_\beta(u),
      \qquad
      (\iota_I \circ \Theta_\beta^{-1})^* g_\beta(u,y)
      = g_{\Psi,\beta}(u) \oplus g_\perp,
    \end{equation}
    where $g_{\Psi,\beta}(u) := D\Psi_\beta(u)\transp D\Psi_\beta(u) =
    M(\rho'(\beta u))$ and $g_\perp := (J_0 Q_\perp)\transp(J_0 Q_\perp)
    \succ 0$.
  \item Under positive path rank at $I$, the normal space to
    $\widetilde{F}_\beta(\Theta_\beta(\widehat{U})) = F_\beta(U)$ at $(u,y) \in
    \Theta_\beta(\widehat{U})$ equals that of $\Psi_\beta(u)$ in
    $\mathcal{T}_I^\perp$, with orthogonal projector $\Pi^\perp(\rho'(\beta
    u))$, and $\mathrm{II}^{\widetilde{F}_\beta}(\partial_{y_\alpha}, \cdot) =
    0$, $\mathrm{II}^{\widetilde{F}_\beta}(\partial_{u_a}, \partial_{u_b}) =
    \mathrm{II}^{\Psi_\beta}(\partial_{u_a}, \partial_{u_b})$.
\end{enumerate}
\end{lemma}

\begin{proof}
Item~1 holds because $\varphi_{a,\beta}$ depends only on $u_a$. For Item~2,
\eqref{eq:local-model} and $I_n = (J_0 Q_\perp)(J_0 Q_\perp)\pinv + P_I$ rewrite
$u_a J_0 v_a + \rho_\beta(u_a)a_{h_a} = J_0 Q_\perp\varphi_{a,\beta}(u_a) +
u_a r_a + \rho_\beta(u_a)p_a$; summing over $a$ gives~\eqref{eq:metric-split},
since $\im\partial_y\widetilde{F}_\beta = \mathcal{T}_I \perp \mathcal{T}_I^\perp
\ni \partial_{u_a}\widetilde{F}_\beta = X_a(\rho'(\beta u_a))$. Item~3 follows
from $(\mathcal{T}_I \oplus \im X(\rho'(\beta u)))^\perp =
\im\Pi^\perp(\rho'(\beta u))$ and affinity of $\widetilde{F}_\beta$ in $y$.
\end{proof}

Since $g_\perp$ is constant, the curvature of $g_\beta$ reduces to that of the
translation submanifold $\Psi_\beta : \RR^q \to \mathcal{T}_I^\perp$. This
submanifold is additively separable, so its second fundamental form is
diagonal in the coordinates $u$. Our curvature sign convention is
$R(U,V)W = \nabla_U\nabla_V W -
\nabla_V\nabla_U W - \nabla_{[U,V]}W$ and $R_{abcd} =
\ip{R(\partial_{u_a},\partial_{u_b})\partial_{u_d}}{\partial_{u_c}}$, so that
$R_{1212} > 0$ on the round sphere.

\begin{lemma}
\label{lem:translation-gauss}
Assume positive path rank at $I$, let $\rho \in \Pprof$, and write
$s(u) := (\rho'(\beta u_1),\dots,\rho'(\beta u_q)) \in [0,1]^q$. Then
$\partial_{u_a}\partial_{u_b}\Psi_\beta \equiv 0$ for $a \ne b$, the second
fundamental form of $\Psi_\beta$ in $\mathcal{T}_I^\perp$ is
\begin{equation}\label{eq:translation-ii}
  \mathrm{II}(\partial_{u_a}, \partial_{u_b})
  = \delta_{ab}\,\rho''_\beta(u_a)\,\Pi^\perp(s(u))\,p_a,
\end{equation}
and with $N_{ab}(s) := \ip{\Pi^\perp(s)p_a}{\Pi^\perp(s)p_b} =
p_a\transp\Pi^\perp(s)p_b$, the Riemann tensor of $g_{\Psi,\beta}(u) = M(s(u))$
is
\begin{equation}\label{eq:binomial-gauss}
  R_{abcd}(u)
  = \rho''_\beta(u_a)\,\rho''_\beta(u_b)\,N_{ab}(s(u))\,
    \bigl(\delta_{ac}\delta_{bd} - \delta_{ad}\delta_{bc}\bigr).
\end{equation}
\end{lemma}

\begin{proof}
Additive separability gives $\partial_{u_a}\partial_{u_b}\Psi_\beta(u) =
\delta_{ab}\rho''_\beta(u_a)p_a$; projecting via $\Pi^\perp(s(u))$
yields~\eqref{eq:translation-ii}, and the Gauss equation $R_{abcd} =
\ip{\mathrm{II}_{ac}}{\mathrm{II}_{bd}} -
\ip{\mathrm{II}_{ad}}{\mathrm{II}_{bc}}$ yields~\eqref{eq:binomial-gauss}.
\end{proof}

For $N = (N_{ab}) \in \Spec^{2k}$, its \emph{Hafnian}~\cite{barvinok2016} sums
products over the set $\Match{2k}$ of all $(2k-1)!! = (2k)!/(2^k k!)$ perfect
matchings of $\{1,\dots,2k\}$:
\[
  \Haf(N)
  := \frac{1}{2^k k!}\sum_{\pi \in S_{2k}} \prod_{j=1}^k N_{\pi(2j-1)\pi(2j)}
  = \sum_{\{\{a_1,b_1\},\dots,\{a_k,b_k\}\} \in \Match{2k}}
    \prod_{j=1}^k N_{a_j b_j}.
\]
As an oriented $2k$-form, $e_{2k}(g)$ pulls back with sign $\sgn\det D\psi$
under a diffeomorphism $u = \psi(\tilde u)$. Replacing
$\dd u_1 \wedge \cdots \wedge \dd u_{2k}$ by $\dd u$ instead gives an
unoriented Riemannian $1$-density, invariant under arbitrary diffeomorphisms.
Below, integrals of $e_{2k}(g)$ over domains in $\RR^{2k}$ or $\RR^{\dimin}$
always mean integrals of this unoriented density (cf.\ $\mu_\beta^{(2k)}$
in~\eqref{eq:smoothed-gbc-measure}).

\begin{proposition}
\label{prop:hafnian-euler}
Let $q = 2k \ge 2$ be even.
\begin{enumerate}
  \item On an oriented $2k$-manifold $(U, g)$ with positively oriented
    coordinates $(u_1,\dots,u_{2k})$, the Gauss--Bonnet--Chern Euler
    $2k$-form~\cite{chern1944}, \cite[Vol.~V, Ch.~13]{spivak1999} is
    \begin{equation}\label{eq:coord-pfaffian}
      e_{2k}(g)
      = \frac{\mathsf{P}_{2k}(R)}{(8\pi)^k\,k!\,\sqrt{\det g}}\,
        \dd u_1 \wedge \cdots \wedge \dd u_{2k},
    \end{equation}
    normalised so that $\int_{\mathcal{X}} e_{2k}(g) = \chi(\mathcal{X})$ on
    every closed oriented $2k$-manifold $\mathcal{X}$, where
    \begin{equation}\label{eq:curvature-twoform}
      \mathsf{P}_{2k}(R)
      := \sum_{\pi,\pi' \in S_{2k}} \sgn(\pi)\sgn(\pi')
         \prod_{j=1}^k R_{\pi(2j-1)\pi(2j)\pi'(2j-1)\pi'(2j)}.
    \end{equation}
  \item If $R_{abcd} = \kappa_a \kappa_b N_{ab}(\delta_{ac}\delta_{bd} -
    \delta_{ad}\delta_{bc})$ for $\kappa_1,\dots,\kappa_{2k} \in \RR$ and
    $N \in \Spec^{2k}$, then
    \begin{equation}\label{eq:pfaffian-hafnian}
      \mathsf{P}_{2k}(R)
      = 2^{2k}\,k!\,\Haf(N)\,\prod_{a=1}^{2k} \kappa_a.
    \end{equation}
  \item Under positive path rank at $I$ ($|I| = 2k$) and $\rho \in \Pprof$, the
    Euler $2k$-form of $g_{\Psi,\beta}(u) = M(s(u))$ is
    \begin{equation}\label{eq:hafnian-euler}
    \begin{aligned}
      e_{2k}(g_{\Psi,\beta})(u)
      &= h_{2k}(\rho'(\beta u))
         \Bigl(\prod_{a=1}^{2k} \rho''_\beta(u_a)\Bigr)
         \dd u_1 \wedge \cdots \wedge \dd u_{2k},\\
      h_{2k}(s)
      &:= \frac{\Haf(N(s))}{(2\pi)^k\sqrt{\det M(s)}}.
    \end{aligned}
    \end{equation}
\end{enumerate}
\end{proposition}

\begin{proof}
Formula~\eqref{eq:coord-pfaffian} is standard~\cite{chern1944, spivak1999}.
Under $u = \psi(\tilde u)$, $\mathsf{P}_{2k}(R)$ scales by $(\det D\psi)^2$,
$(\det g)^{-1/2}$ by $|\det D\psi|^{-1}$, and $\dd u_1\wedge\cdots\wedge\dd
u_{2k}$ by $\det D\psi$. The net factor is $\sgn\det D\psi = +1$ on oriented
charts, and the unoriented density with $\dd u$ is invariant.

For~\eqref{eq:pfaffian-hafnian}, fix $\pi \in S_{2k}$. The product
in~\eqref{eq:curvature-twoform} is non-zero if and only if
$\{\pi'(2j-1),\pi'(2j)\} = \{\pi(2j-1),\pi(2j)\}$ for each $j$. Transposing
$\ell$ of these pairs multiplies both $\sgn(\pi')$ and the Kronecker product by
$(-1)^\ell$. Summing over the $2^k$ such $\pi'$ and over $\pi \in S_{2k}$ gives
$\mathsf{P}_{2k}(R) = 2^k(\prod_a\kappa_a)\sum_{\pi\in S_{2k}}\prod_{j=1}^k
N_{\pi(2j-1)\pi(2j)} = 2^{2k}k!\Haf(N)\prod_a\kappa_a$.

For~\eqref{eq:hafnian-euler}, set $\kappa_a = \rho''_\beta(u_a)$
in~\eqref{eq:binomial-gauss} and use $(8\pi)^k k! = 2^{2k}(2\pi)^k k!$.
Under $s_a = \rho'(\beta u_a)$, the resulting form $e_{2k}(g_{\Psi,\beta})$
pushes forward to $h_{2k}(s)\,\dd s$ on $[0,1]^{2k}$.
\end{proof}

\section{Identification with Banchoff and Lipschitz--Killing Curvature}\label{sec:banchoff}

We show that for every even codimension $2k \ge 2$, $\int_{[0,1]^{2k}}
h_{2k}(s)\,\dd s$ equals Banchoff's critical-point curvature
$K(\Phi_I)$~\cite{banchoff1967} of the transverse cone $\Psi_0(\RR^{2k})$ at
its apex, and coincides with the Cheeger--M\"uller--Schrader curvature density
$(2\pi)^{-k}\mathcal{R}^{2k}(\Phi_I)$ along $\Phi_I$.

\subsection{Kac--Rice critical-point formula and PL lower links}\label{subsec:kac-rice-links}

Fix an ordered tuple $I = (h_1,\dots,h_{2k})$ ($2k \ge 2$) satisfying local
transversality and positive path rank at $I$, together with $\rho \in \Pprof$
and $\beta > 0$. For $\xi \in \mathcal{T}_I^\perp$, define the outward ray
slopes $\ell_a^-(\xi) := -\ip{\xi}{r_a}$ and $\ell_a^+(\xi) :=
\ip{\xi}{r_a + p_a}$ of $\Psi_0$ and the signed critical-point indicator
\begin{equation}\label{eq:ray-slopes}
  \varepsilon_a(\xi)
  := \mathbf{1}\{\ell_a^-(\xi)>0,\,\ell_a^+(\xi)>0\}
     - \mathbf{1}\{\ell_a^-(\xi)<0,\,\ell_a^+(\xi)<0\}.
\end{equation}

\begin{lemma}
\label{lem:separable-crit}
For $\xi \in \mathcal{T}_I^\perp$, define
$f_{\xi,\beta}(u) := \ip{\xi}{\Psi_\beta(u)} = \sum_{a=1}^{2k}\bigl(-u_a\ell_a^-(\xi)
+ \rho_\beta(u_a)(\ell_a^+(\xi)+\ell_a^-(\xi))\bigr)$ ($u \in \RR^{2k}$) and
$E_{\mathrm{lin}} := \bigcup_{a=1}^{2k}\{\xi \in \mathcal{T}_I^\perp :
\ell_a^-(\xi)\ell_a^+(\xi) = 0\}$. Then $E_{\mathrm{lin}}$ is a finite union of
proper hyperplanes in $\mathcal{T}_I^\perp$ (of Lebesgue and Gaussian measure
zero), and for every $\xi \in \mathcal{T}_I^\perp \setminus E_{\mathrm{lin}}$,
$f_{\xi,\beta}$ is Morse on $\RR^{2k}$: if $\ell_a^-(\xi)\ell_a^+(\xi) > 0$ for
all $a$, it has a unique critical point $u^*(\xi,\beta) \in
(-\beta^{-1},\beta^{-1})^{2k}$, given by $\rho'(\beta u_a^*) =
\ell_a^-(\xi)/(\ell_a^-(\xi)+\ell_a^+(\xi)) \in (0,1)$, and otherwise none. In
both cases,
\begin{equation}\label{eq:separable-morse-sum}
  \sum_{u \in \RR^{2k} :\,\nabla f_{\xi,\beta}(u)=0}
    \sgn\det\operatorname{Hess} f_{\xi,\beta}(u)
  = \prod_{a=1}^{2k} \varepsilon_a(\xi).
\end{equation}
\end{lemma}

\begin{proof}
Positive path rank gives $r_a = X_a(0) \ne 0$ and $r_a + p_a = X_a(1) \ne 0$,
so $E_{\mathrm{lin}}$ is a finite union of proper hyperplanes. Writing
$f_{\xi,\beta}(u) = \sum_{a=1}^{2k} f_a(u_a)$, we have
\[
  f'_a(u_a) = -\ell_a^-(\xi) + (\ell_a^-(\xi)+\ell_a^+(\xi))\rho'(\beta u_a),
  \qquad
  f''_a(u_a) = (\ell_a^-(\xi)+\ell_a^+(\xi))\rho''_\beta(u_a).
\]
Fix $\xi \notin E_{\mathrm{lin}}$, so $\ell_a^\pm(\xi) \ne 0$. If
$\ell_a^-(\xi)+\ell_a^+(\xi) = 0$, then $f'_a \equiv -\ell_a^-(\xi) \ne 0$ has
no roots and $\varepsilon_a(\xi) = 0$. Otherwise $f'_a(u_a) = 0 \iff
\rho'(\beta u_a) = \ell_a^-(\xi)/(\ell_a^-(\xi)+\ell_a^+(\xi))$, and this value
lies in $(0,1)$ if and only if $\ell_a^-(\xi)\ell_a^+(\xi) > 0$.

The function $\rho'(\beta\,\cdot\,)$
increases strictly from $0$ to $1$ on $[-\beta^{-1},\beta^{-1}]$, with
$\rho''_\beta > 0$ on $(-\beta^{-1},\beta^{-1})$. So in that case there is a
unique root $u_a^* \in (-\beta^{-1},\beta^{-1})$, with $\sgn f''_a(u_a^*) =
\sgn(\ell_a^-(\xi)+\ell_a^+(\xi)) = \varepsilon_a(\xi)$, and none otherwise.
Multiplying across $a$ gives~\eqref{eq:separable-morse-sum}.
\end{proof}

\begin{lemma}
\label{lem:kac-rice-hafnian}
Let $\xi = P_I\zeta \sim \Norm(0, P_I)$ ($\zeta \sim \Norm(0,I_n)$), and set
$U_\xi(u) := \nabla f_{\xi,\beta}(u) = X(s(u))\transp\xi \in \RR^{2k}$ and
$\mathsf{D}_\xi(u) := \operatorname{Hess} f_{\xi,\beta}(u) =
\diag(\rho''_\beta(u_a)\ip{p_a}{\xi})_{a=1}^{2k}$. At each $u \in \RR^{2k}$,
$U_\xi(u) \sim \Norm(0, M(s(u)))$ has density $p_{U_\xi(u)}(0) =
(2\pi)^{-k}/\sqrt{\det M(s(u))}$. Conditional on $U_\xi(u) = 0$, the vector
$(\ip{p_a}{\xi})_{a=1}^{2k}$ has law $\Norm(0, N(s(u)))$. By Gaussian
regression~\cite{adlertaylor2007} and Isserlis' theorem~\cite{isserlis1918},
\begin{equation}\label{eq:kac-rice-pointwise}
  \EE\bigl[\det \mathsf{D}_\xi(u) \bigm| U_\xi(u) = 0\bigr]\,p_{U_\xi(u)}(0)
  = h_{2k}(s(u))\prod_{a=1}^{2k} \rho''_\beta(u_a),
\end{equation}
whose integral over $\RR^{2k}$ satisfies
\begin{equation}\label{eq:kac-rice-integral}
  \int_{\RR^{2k}} \EE[\det \mathsf{D}_\xi(u) \mid U_\xi(u) = 0]\,p_{U_\xi(u)}(0)\,\dd u
  = \EE_\xi\Bigl[\sum_{U_\xi(u)=0} \sgn\det \mathsf{D}_\xi(u)\Bigr]
  = \EE_\xi\Bigl[\prod_{a=1}^{2k} \varepsilon_a(\xi)\Bigr].
\end{equation}
\end{lemma}

\begin{proof}
Since $\im X(s) \subset \mathcal{T}_I^\perp$ and $\Pi^\perp(s)X(s) = 0$
by~\eqref{eq:schur-projector}, $\Cov(U_\xi(u)) = M(s) \succ 0$ and
$\Cov(\Pi^\perp(s)\xi, U_\xi(u)) = 0$, so $\Pi^\perp(s)\xi \sim
\Norm(0,\Pi^\perp(s))$ is independent of $U_\xi(u)$. Conditional on
$U_\xi(u) = 0$ we have $\xi = \Pi^\perp(s)\xi$, so
$Y := (\ip{p_a}{\xi})_{a=1}^{2k} \sim \Norm(0, N(s(u)))$ and Isserlis'
theorem~\cite{isserlis1918} gives $\EE[\prod_{a=1}^{2k}Y_a \mid U_\xi(u)=0] =
\Haf(N(s(u)))$, proving~\eqref{eq:kac-rice-pointwise}.

Since $\supp\rho''_\beta \subset [-\beta^{-1},\beta^{-1}]$, both
\eqref{eq:kac-rice-pointwise} and $\det\mathsf{D}_\xi(u)$ vanish outside
$\Rtube_\beta = [-\beta^{-1},\beta^{-1}]^{2k}$, so it suffices to
prove~\eqref{eq:kac-rice-integral} on $\Rtube_\beta$. Since $\mathsf{D}_\xi$
degenerates on $\partial\Rtube_\beta$, we cannot appeal to boundary Hessian
non-degeneracy; instead we verify the Kac--Rice
identity~\cite[Thm.~11.2.1]{adlertaylor2007} directly.
Let $\eta_{1,\delta}$ be the centred Gaussian density of variance $\delta^2$ on
$\RR$ and $\eta_\delta(z) := \prod_{a=1}^{2k}\eta_{1,\delta}(z_a)$. We evaluate
$\EE_\xi\int_{\Rtube_\beta}\eta_\delta(U_\xi(u))\det\mathsf{D}_\xi(u)\,\dd u$
in both orders of integration and let $\delta \downarrow 0$.

First integrate in $\xi$. From $P_I\xi = X(s)M(s)^{-1}U_\xi(u) +
\Pi^\perp(s)\xi$, conditional on $U_\xi(u) = z$ we have $\ip{p_a}{\xi} =
p_a\transp X(s)M(s)^{-1}z + \mathsf{W}_a$ with $\mathsf{W} \sim \Norm(0, N(s))$.
Since $M(s) \succeq \lambda_0 I_{2k}$ and $\norm{N(s)}_{\mathrm{op}} \le \tr N(s)
\le \sum_{a=1}^{2k}\norm{p_a}_2^2$, the function $z \mapsto
\EE[\det\mathsf{D}_\xi(u)\mid U_\xi(u)=z]$ is a polynomial of degree $\le 2k$
whose coefficients are bounded uniformly in $u \in \Rtube_\beta$. Moreover
\[
  \EE[\abs{\det\mathsf{D}_\xi(u)}\mid U_\xi(u)=z] \le C(1 + \norm{z}_2^{2k}),
  \qquad
  p_{U_\xi(u)}(z) \le C' e^{-c\norm{z}_2^2},
\]
with constants independent of
$u$, $z$, and $\delta$. Hence $\EE_\xi[\eta_\delta(U_\xi(u))
\abs{\det\mathsf{D}_\xi(u)}]$ is bounded uniformly on $\Rtube_\beta$ for
$\delta \in (0,1]$. As $z \mapsto \EE[\det\mathsf{D}_\xi(u)\mid
U_\xi(u)=z]\,p_{U_\xi(u)}(z)$ is continuous, it follows that
\[
  \EE_\xi\bigl[\eta_\delta(U_\xi(u))\det\mathsf{D}_\xi(u)\bigr]
  \longrightarrow
  \EE[\det\mathsf{D}_\xi(u)\mid U_\xi(u)=0]\,p_{U_\xi(u)}(0)
  \qquad (\delta \downarrow 0)
\]
boundedly on $\Rtube_\beta$.

Now integrate first in $u$ at fixed $\xi$. Write $f_{\xi,\beta}(u) =
\sum_a f_a(u_a)$ as in the proof of Lemma~\ref{lem:separable-crit}. The
coordinate $U_{\xi,a}(u) = f'_a(u_a)$ depends only on $u_a$, and
$\partial_{u_a}U_{\xi,a} = f''_a(u_a) = \mathsf{D}_{\xi,aa}(u)$. Hence both
$\eta_\delta(U_\xi(u)) = \prod_a\eta_{1,\delta}(f'_a(u_a))$ and
$\det\mathsf{D}_\xi(u) = \prod_a f''_a(u_a)$ factor over the coordinates.
Setting $G_\delta(z) := \int_{-\infty}^z\eta_{1,\delta}(t)\,\dd t$ and using
$f'_a(-\beta^{-1}) = -\ell_a^-(\xi)$, $f'_a(\beta^{-1}) = \ell_a^+(\xi)$, the
one-dimensional fundamental theorem of calculus gives, for each $a$,
\begin{align*}
  \int_{-\beta^{-1}}^{\beta^{-1}} \eta_{1,\delta}\bigl(f'_a(u_a)\bigr)\,f''_a(u_a)\,\dd u_a
  &= \Bigl[G_\delta\bigl(f'_a(u_a)\bigr)\Bigr]_{u_a=-\beta^{-1}}^{u_a=\beta^{-1}}\\
  &= G_\delta\bigl(\ell_a^+(\xi)\bigr) - G_\delta\bigl(-\ell_a^-(\xi)\bigr)
  = \int_{-\ell_a^-(\xi)}^{\ell_a^+(\xi)} \eta_{1,\delta}(z_a)\,\dd z_a,
\end{align*}
a signed change of variables: when $\ell_a^-(\xi) + \ell_a^+(\xi) < 0$, the
map $f'_a$ is decreasing and the oriented interval is traversed backwards.
Taking the product over $a$,
\[
  \int_{\Rtube_\beta} \eta_\delta\bigl(U_\xi(u)\bigr)\det \mathsf{D}_\xi(u)\,\dd u
  = \prod_{a=1}^{2k}\int_{-\ell_a^-(\xi)}^{\ell_a^+(\xi)} \eta_{1,\delta}(z_a)\,\dd z_a
  \in [-1,1],
\]
and as $\delta \downarrow 0$ each factor converges, for $\xi \notin
E_{\mathrm{lin}}$, to $\mathbf{1}\{\ell_a^-(\xi)>0,\,\ell_a^+(\xi)>0\} -
\mathbf{1}\{\ell_a^-(\xi)<0,\,\ell_a^+(\xi)<0\} = \varepsilon_a(\xi)$. The
product therefore converges to $\prod_a\varepsilon_a(\xi) =
\sum_{U_\xi(u)=0}\sgn\det\mathsf{D}_\xi(u)$ by Lemma~\ref{lem:separable-crit}.

Since the product is bounded by $1$, dominated convergence in $\xi$ shows that
$\EE_\xi\int_{\Rtube_\beta}\eta_\delta(U_\xi(u))\det\mathsf{D}_\xi(u)\,\dd u \to
\EE_\xi[\prod_a\varepsilon_a(\xi)]$. On the other hand, Fubini's theorem and
the bounded convergence on $\Rtube_\beta$ established above show that the same
quantity converges to $\int_{\Rtube_\beta}\EE[\det\mathsf{D}_\xi(u)\mid
U_\xi(u)=0]\,p_{U_\xi(u)}(0)\,\dd u$. This proves~\eqref{eq:kac-rice-integral}.
\end{proof}

Let $\mathcal{F}_q$ be the orthant fan on $\RR^q$, with non-apex cells
$C_\varsigma^\circ := \{u \in \RR^q : \sgn(u_a) = \varsigma_a\ \forall a\}$
($\varsigma \in \{-,0,+\}^q \setminus \{0\}$, $\dim C_\varsigma = |\varsigma| :=
\#\{a : \varsigma_a \ne 0\}$). For $f : \RR^q \to \RR$ linear on each orthant
with ray values $f(e_a) = \ell_a^+$ and $f(-e_a) = \ell_a^-$ (matching
\eqref{eq:ray-slopes} for $f = \xi\cdot\Psi_0$), Banchoff's apex
index~\cite{banchoff1967} is
\begin{equation}\label{eq:banchoff-index-def}
  \ind_0(f)
  := \sum_{\varsigma \in \{-,0,+\}^q} (-1)^{|\varsigma|}\,
     \mathbf{1}\bigl\{0\text{ is the strict maximum of }f\text{ on }\overline{C_\varsigma}\bigr\},
\end{equation}
where the $\varsigma = 0$ term equals $+1$. Because $S_1^{q-1} := \{u \in \RR^q
: \norm{u}_1 = 1\}$ is the join $\{\pm e_1\} * \cdots * \{\pm
e_q\}$~\cite{banchoff1967}, the strict-descent complex factors across
coordinates.

\begin{lemma}
\label{lem:lower-link-index}
Let $q \ge 1$, let $f : \RR^q \to \RR$ be linear on each orthant of
$\mathcal{F}_q$ with $f(\pm e_a) = \ell_a^\pm$, and let $\Delta_-(f)$ be the
full subcomplex of $S_1^{q-1}$ on $V_- := \{v \in \{\pm e_1,\dots,\pm e_q\} :
f(v) < 0\}$. Then:
\begin{enumerate}
  \item $0 = \arg\max_{\overline{C_\varsigma}} f$ strictly if and only if every
    active coordinate ray of $C_\varsigma$ lies in $\Delta_-(f)$, so
    \begin{equation}\label{eq:lower-link-euler}
      \ind_0(f)
      = 1 - \chi\bigl(\Delta_-(f)\bigr)
      = \prod_{a=1}^q \bigl(1 - \mathbf{1}\{\ell_a^+ < 0\} - \mathbf{1}\{\ell_a^- < 0\}\bigr).
    \end{equation}
  \item $\Lk_-(f) := \{u \in S_1^{q-1} : f(u) < 0\}$ deformation-retracts onto
    $|\Delta_-(f)|$, so $\ind_0(f) = 1 - \chi(\Lk_-(f))$
    (cf.~\cite{banchoff1967, grunertkuhnelrote2023}).
  \item When $\ell_a^\pm \ne 0$ for all $a$, the $a$-th factor
    in~\eqref{eq:lower-link-euler} equals $\varepsilon_a =
    \mathbf{1}\{\ell_a^->0,\ell_a^+>0\} - \mathbf{1}\{\ell_a^-<0,\ell_a^+<0\}$.
\end{enumerate}
\end{lemma}

\begin{proof}
Linearity of $f$ on $\overline{C_\varsigma} = \cone\{\varsigma_a e_a :
\varsigma_a \ne 0\}$ gives $0 = \arg\max_{\overline{C_\varsigma}} f \iff
f(\varsigma_a e_a) < 0$ for all $\varsigma_a \ne 0$, and this factors
\eqref{eq:banchoff-index-def} into~\eqref{eq:lower-link-euler}. The
simplices of $\Delta_-(f)$ are $\overline{C_\varsigma} \cap S_1^{q-1}$
($\varsigma \ne 0$, $\dim = |\varsigma|-1$), so $\chi(\Delta_-(f)) = 1 -
\ind_0(f)$.

For Item~2, write $u = \sum_{v\in\{\pm e_a\}} t_v(u)v$ on $S_1^{q-1}$, where
$t_v(u) := \max(\ip{u}{v},0)$. At most one of $t_{\pm e_a}(u)$ is positive.
Take $u \in \Lk_-(f)$, and note that $\nu(u) := \sum_{v\in V_-} t_v(u) > 0$.
Then $u$ and $\nu(u)^{-1}\sum_{v\in V_-} t_v(u)v$ lie in a common
$\overline{C_\varsigma}$ on which $\norm{\cdot}_1$ and $f$ are linear. Hence
$h_t(u) := (1-t)u +
t\nu(u)^{-1}\sum_{v\in V_-} t_v(u)v$ ($t \in [0,1]$) deformation-retracts
$\Lk_-(f)$ onto $|\Delta_-(f)|$. Item~3 follows by checking the four sign
pairs.
\end{proof}

\begin{lemma}
\label{lem:pr-injectivity}
For an ordered $q$-tuple $I$ ($q \ge 1$) satisfying local transversality with
$\rank(J_0 Q_\perp) = \dimin - q$, set $\Psi_0(u) := \sum_{a=1}^q (u_a r_a +
(u_a)_+ p_a) : \RR^q \to \mathcal{T}_I^\perp$, $X_\tau :=
[X_1(\tfrac{1+\tau_1}{2})\ \cdots\ X_q(\tfrac{1+\tau_q}{2})]$, and $M_\tau :=
X_\tau\transp X_\tau = M(\tfrac{\mathbf{1}+\tau}{2})$ ($\tau \in \{-,+\}^q$).
Then:
\begin{enumerate}
  \item Positive path rank holds at $I$ if and only if $\Psi_0$ is globally
    injective, in which case $\Psi_0$ is a bi-Lipschitz homeomorphism onto
    $\mathcal{K}_I := \Psi_0(\RR^q)$ inducing the simplicial cone complex
    $\{\Psi_0(\overline{C_\varsigma})\}_{\varsigma\in\{-,0,+\}^q} \cong
    \mathcal{F}_q$.
  \item For $q = 2k$ even, if $r_a \ne 0$ and $r_a + p_a \ne 0$ for all $a$ (as
    when every $M_\tau \succ 0$), then for all $\xi \in \mathcal{T}_I^\perp
    \setminus E_{\mathrm{lin}}$,
    \[
      \ind_0(\xi\cdot\Psi_0)
      = \prod_{a=1}^{2k} \varepsilon_a(\xi).
    \]
\end{enumerate}
\end{lemma}

\begin{proof}
For $u, u' \in \RR^q$, $(u_a)_+ - (u'_a)_+ = \bar{s}_a(u_a - u'_a)$ with
$\bar{s}_a := \int_0^1 \mathbf{1}\{t u_a + (1-t)u'_a > 0\}\,\dd t \in [0,1]$,
so $\Psi_0(u) - \Psi_0(u') = X(\bar{s})(u - u')$. Positive path rank gives
$\sqrt{\lambda_0}\norm{u-u'}_2 \le \norm{\Psi_0(u)-\Psi_0(u')}_2 \le
\sqrt{\overline{\lambda}_0}\norm{u-u'}_2$. Conversely, if $X(s)c = 0$ for
$s \in [0,1]^q$ and $c \ne 0$, setting $(u_a, u'_a) := (s_a c_a, -(1-s_a)c_a)$
for $c_a \ge 0$ and $((1-s_a)c_a, -s_a c_a)$ for $c_a < 0$ gives $u - u' = c$
and $\Psi_0(u) = \Psi_0(u')$. Lemma~\ref{lem:lower-link-index}(1),(3) applied
to $f = \xi\cdot\Psi_0$ gives Item~2.
\end{proof}

\subsection{Banchoff curvature, orthant probabilities, and Schur complements}\label{subsec:banchoff-schur}

Define the \emph{Banchoff curvature} of $\Phi_I$ by
\[
  K(\Phi_I)
  := \EE_{\xi\sim\Norm(0,P_I)}\bigl[\ind_0(\xi\cdot\Psi_0)\bigr]
   = \frac{1}{\abs{\Sph{2k+m-1}}}\int_{\Sph{2k+m-1}(\mathcal{T}_I^\perp)}
     \ind_0(\xi\cdot\Psi_0)\,\dd\omega(\xi),
\]
which equals $\abs{\Sph{n-1}}^{-1}\int_{\Sph{n-1}}\ind_0(\xi\cdot\Psi_0)\,
\dd\omega_{n-1}(\xi)$ since $\ind_0(\xi\cdot\Psi_0)$ depends only on the
direction of $P_I\xi$. Moreover, $K(\Phi_I)$ is constant on each connected
component of $\Phi_I\setminus\bigcup_{h\notin I}\wall{h}$. When $2k = \dimin$
this is Banchoff's vertex curvature~\cite{banchoff1967}. For flats of positive
dimension it extends that curvature transversally, in the spirit
of~\cite{banchoff1983}.

\begin{theorem}[Identification with Banchoff curvature and orthant sum]
\label{thm:banchoff-orthant}
In the setting of Section~\ref{subsec:kac-rice-links}:
\begin{equation}\label{eq:banchoff-orthant}
  \int_{\RR^{2k}} e_{2k}(g_{\Psi,\beta})
  = \int_{[0,1]^{2k}} h_{2k}(s)\,\dd s
  = \EE_\xi\Bigl[\prod_{a=1}^{2k} \varepsilon_a(\xi)\Bigr]
  = K(\Phi_I)
  = \sum_{\tau\in\{\pm\}^{2k}} (-1)^{\#\tau_+}\PP(Y^\tau > 0),
\end{equation}
where $\#\tau_+ := \#\{a : \tau_a = +\}$ and $Y^\tau \sim \Norm(0, M_\tau)$
with $M_\tau = M(\tfrac{\mathbf{1}+\tau}{2})$. For $k = 1$, writing $\theta_\tau
:= \arccos((-1)^{\#\tau_+}(M_\tau)_{12}/\sqrt{(M_\tau)_{11}(M_\tau)_{22}}) \in
(0,\pi)$ for the interior angle of $\Psi_0(C_\tau)$,
\begin{equation}\label{eq:codim2-arcsin}
  K(\Phi_{\{h_1,h_2\}})
  = 1 - \frac{1}{2\pi}\sum_{\tau\in\{-,+\}^2} \theta_\tau
  = \frac{1}{2\pi}\sum_{\tau\in\{-,+\}^2} (-1)^{\#\tau_+}
    \arcsin\!\left(\frac{(M_\tau)_{12}}{\sqrt{(M_\tau)_{11}(M_\tau)_{22}}}\right).
\end{equation}
\end{theorem}

\begin{proof}
By Proposition~\ref{prop:hafnian-euler}(3), $e_{2k}(g_{\Psi,\beta})$ is
supported on $\Rtube_\beta$ and pushes forward under the diffeomorphism
$s_\beta(u) := (\rho'(\beta u_1),\dots,\rho'(\beta u_{2k})) :
(-\beta^{-1},\beta^{-1})^{2k} \to (0,1)^{2k}$ to $h_{2k}(s)\,\dd s$ on
$[0,1]^{2k}$. Thus $\int_{\Rtube_\beta}(f\circ s_\beta)e_{2k}(g_{\Psi,\beta}) =
\int_{[0,1]^{2k}} f h_{2k}$ for bounded Borel $f$. Taking $f \equiv 1$ and
applying Lemmas~\ref{lem:separable-crit}, \ref{lem:kac-rice-hafnian},
and~\ref{lem:pr-injectivity}(2) gives the first three equalities
in~\eqref{eq:banchoff-orthant}.

For the orthant sum, Lemma~\ref{lem:lower-link-index}(1),(3) gives, off
$E_{\mathrm{lin}}$, $\varepsilon_a(\xi) =
\mathbf{1}\{-\ell_a^-(\xi)<0\} - \mathbf{1}\{\ell_a^+(\xi)<0\} =
\sum_{\tau_a\in\{-,+\}}(-1)^{\mathbf{1}\{\tau_a=+\}}
\mathbf{1}\{-\ip{\xi}{X_a(\tfrac{1+\tau_a}{2})} > 0\}$. Multiplying over $a$,
setting $Y^\tau := -X_\tau\transp\xi \sim \Norm(0, M_\tau)$, and taking
expectations yields the orthant sum in~\eqref{eq:banchoff-orthant}. On the
product fan $\mathcal{F}_{2k}$, this expansion is the coordinate-factored form
of the Brianchon--Gram relation~\cite{mcmullen1975, schneiderweil2008}.

For $k = 1$, the ray images of $C_\tau$ are $\tau_a X_a(\tfrac{1+\tau_a}{2})$,
and $\tau_1\tau_2 = (-1)^{\#\tau_+}$. Setting $\varrho_\tau :=
(M_\tau)_{12}/\sqrt{(M_\tau)_{11}(M_\tau)_{22}}$, we get $\theta_\tau =
\tfrac{\pi}{2} - (-1)^{\#\tau_+}\arcsin \varrho_\tau$. Sheppard's
formula~\cite{sheppard1899}, $\PP(Y^\tau_1>0, Y^\tau_2>0) = \tfrac{1}{4} +
\frac{1}{2\pi}\arcsin \varrho_\tau$, together with $\sum_{\tau}(-1)^{\#\tau_+}
= 0$, then gives~\eqref{eq:codim2-arcsin}.
\end{proof}

\begin{remk}\label{rem:profile-scope}
We record two consequences of the proof.
\begin{enumerate}
  \item The shape of $\rho$ enters~\eqref{eq:banchoff-orthant} only through the
    change of variables $s_a = \rho'(\beta u_a)$. Consequently
    \eqref{eq:banchoff-orthant} holds for the translation submanifold
    $\Psi_\beta$ of any smooth convex ramp with $\rho'' > 0$ on $\RR$,
    $\rho'(-\infty) = 0$, and $\rho'(+\infty) = 1$, such as softplus $\rho(z) =
    \log(1 + e^z)$. Indeed, Lemma~\ref{lem:separable-crit} applies verbatim, and in
    Lemma~\ref{lem:kac-rice-hafnian} the integrals over $\Rtube_\beta$ become
    integrals over $\RR^{2k}$ dominated by a constant multiple of
    $\prod_a\rho''_\beta(u_a)$. For the smoothed map $F_\beta$ itself this
    covers the case $I = \{1,\dots,\ngen\}$, in which no further wall enters
    the local model~\eqref{eq:local-model}. For general $I$, and throughout
    Section~\ref{sec:global}, we need the compact support $\supp\rho'' =
    [-1,1]$ built into $\Pprof$: it keeps the transition tubes of different
    walls disjoint at a finite scale $\beta \ge \beta_\sharp$
    (Lemma~\ref{lem:tube-separation}).
  \item When $n = \dimin$ ($m = 0$), positive path rank forces $\Pi^\perp(s)
    \equiv 0$, hence $N(s) \equiv 0$, $h_{2k} \equiv 0$, and $K(\Phi_I) = 0$: a
    locally injective piecewise-affine map $\RR^{\dimin} \to \RR^{\dimin}$ is a
    local isometry off the walls, and so are its smoothings.
    Remark~\ref{rem:PR-dichotomy} shows that positive path rank cannot be
    dropped from this statement.
\end{enumerate}
\end{remk}

\begin{remk}\label{rem:PR-dichotomy}
We assume $\dimin \ge 2$ in Sections~\ref{sec:banchoff}--\ref{sec:global},
since no even-codimension strata exist for $\dimin = 1$. The algebraic
definitions of Sections~\ref{sec:hafnian} and~\ref{sec:moduli} allow $\dimin
\ge 1$; odd-codimension vanishing, folding bounds, and orthogonal
factorisation of $K(\Phi_I)$ are recorded in
Proposition~\ref{prop:folding-bounds}. By Lemma~\ref{lem:pr-injectivity}(2),
the last two equalities in~\eqref{eq:banchoff-orthant}
and~\eqref{eq:codim2-arcsin} require only $M_\tau \succ 0$ ($\tau \in
\{-,+\}^{2k}$), whereas failure of positive path rank forces $\det M(s) = 0$
somewhere on $[0,1]^{2k}\setminus\{0,1\}^{2k}$.

Already for $\dimin = n = 2$, $W = L = I_2$, $b = 0$, the reflection fold
$F(x_1,x_2) = (-|x_1|, x_2)\transp$ ($a_1 = (-2,0)\transp$, $a_2 = 0$) has
$g_\tau = I_2$ on all quadrants, yet $\det DF_\beta = 0$ along
$\rho'(\beta x_1) = \tfrac{1}{2}$. Now take $a_1 = (-2,1)\transp$ and
$a_2 = (0,-2)\transp$, that is, $F(x_1,x_2) = (x_1 - 2(x_1)_+,\, x_2 + (x_1)_+ -
2(x_2)_+)\transp$. This keeps $\det g_\tau = 1$ on all four quadrants, with
$\theta_{--} = \theta_{-+} = \tfrac{\pi}{2}$ and $\theta_{+-} = \theta_{++} =
\tfrac{3\pi}{4}$. Hence $(2\pi)^{-1}\mathcal{R}^2(0) = K(0) = -\tfrac{1}{4}$ by
Proposition~\ref{prop:schur-cms}(3); the last two items of that proposition
require only $M_\tau \succ 0$. Off $\{\rho'(\beta x_a) = \tfrac{1}{2}\}$,
however, the map $F_\beta : \RR^2 \to \RR^2$ is locally Euclidean with
$e_2(g_\beta) \equiv 0$. Without single-wall
orientation preservation, convex coordinate-wise smoothings fold across walls
and need not recover Cheeger--M\"uller--Schrader curvature.
\end{remk}

For a locally finite $\dimin$-dimensional piecewise-flat polyhedral space $X$
and an open $(\dimin - 2k)$-cell $\face$ with normal cone $\mathsf{N}(\face)$,
the \emph{codimension-$2k$ Lipschitz--Killing curvature}~\cite[\S1 and
\S3]{cheeger1984} is $\mathcal{R}^{2k} = \sum_{\face} \mathcal{R}^{2k}(\face)\,
\Haus{\dimin-2k}\restriction_{\face}$ with
\begin{equation}\label{eq:cms-ext-angle}
  (2\pi)^{-k} \mathcal{R}^{2k}(\face)
  := \sum_{\mathfrak{f} \text{ cell of } \mathsf{N}(\face)} (-1)^{\dim \mathfrak{f}}\,\gamma\bigl(\{0\}, \mathfrak{f}\bigr),
\end{equation}
where $\gamma(\{0\},\mathfrak{f}) \in [0,1]$ is the normalised exterior angle of
$\mathfrak{f}$ at $0$ ($\gamma(\{0\},\{0\}) := 1$). No embedding of $X$ is
involved. Here $\mathsf{N}(\face)$ is the intrinsic normal cone of $X$ along
$\face$: the piecewise-Euclidean cone complex obtained by gluing the normal
cones of the cells containing $\face$ along their common faces. The angle
$\gamma(\{0\},\mathfrak{f})$, like $\gamma(\mathfrak{f}'',\mathfrak{f})$
and $\alpha(\mathfrak{f}',\mathfrak{f}'')$ below, is the angle of the convex
Euclidean cone $\mathfrak{f}$ computed within its own linear span.

We next relate \eqref{eq:cms-ext-angle} to the interior-angle flag sum
of~\cite[\S3]{cheeger1984},
\[
  \sum_{j\ge 0}(-1)^j\sum_{\{0\}=\mathfrak{f}_0<\cdots<\mathfrak{f}_j}
  \prod_{i=1}^j \alpha(\mathfrak{f}_{i-1},\mathfrak{f}_i)
\]
where $\alpha(\mathfrak{f}',\mathfrak{f}'')$ denotes the normalised interior
angle and $\alpha(\mathfrak{f}',\mathfrak{f}') := 1$. Order the cells of
$\mathsf{N}(\face)$ by inclusion and set
$A_{\mathrm{inc}}(\mathfrak{f}',\mathfrak{f}'') :=
\alpha(\mathfrak{f}',\mathfrak{f}'')\mathbf{1}\{\mathfrak{f}' < \mathfrak{f}''\}$
and $\widetilde{\Gamma}(\mathfrak{f}'',\mathfrak{f}) :=
(-1)^{\dim\mathfrak{f}-\dim\mathfrak{f}''}\gamma(\mathfrak{f}'',\mathfrak{f})
\mathbf{1}\{\mathfrak{f}'' \le \mathfrak{f}\}$. Sommerville's relation
$\sum_{\mathfrak{f}'\le\mathfrak{f}''\le\mathfrak{f}}
\alpha(\mathfrak{f}',\mathfrak{f}'')(-1)^{\dim\mathfrak{f}-\dim\mathfrak{f}''}
\gamma(\mathfrak{f}'',\mathfrak{f}) =
\delta_{\mathfrak{f}',\mathfrak{f}}$~\cite{mcmullen1975},
\cite[Thm.~6.5.5]{schneiderweil2008} reads
$(I + A_{\mathrm{inc}})\widetilde{\Gamma} = I$ in the incidence algebra.
Since $A_{\mathrm{inc}}$ is strictly upper triangular, we may expand
$\widetilde{\Gamma} = (I + A_{\mathrm{inc}})^{-1} = \sum_{j\ge 0}(-1)^j
A_{\mathrm{inc}}^j$. Evaluating at $\mathfrak{f}' = \{0\}$ and summing over
$\mathfrak{f}$ recovers the flag sum of~\cite[\S3]{cheeger1984}. This matches
$\mathcal{R}^{2k}(X) = (2\pi)^k\chi(X)$ on closed
$2k$-manifolds~\cite{banchoff1967, cheeger1984}.

\begin{proposition}
\label{prop:schur-cms}
In the setting of Section~\ref{subsec:kac-rice-links}:
\begin{enumerate}
  \item For all $s \in [0,1]^{2k}$, $M(s) = V_I\transp\bigl(g(s) -
    g(s)Q_\perp(Q_\perp\transp g(s)Q_\perp)^{-1}Q_\perp\transp g(s)\bigr)V_I =
    (W_I g(s)^{-1} W_I\transp)^{-1}$ where $g(s) = J(s)\transp J(s)$ and $J(s) =
    J_0 + \sum_{a=1}^{2k} s_a a_{h_a}w_{h_a}\transp$.
  \item At each vertex $\tfrac{\mathbf{1}+\tau}{2} \in \{0,1\}^{2k}$ ($\tau \in
    \{-,+\}^{2k}$), $M_\tau = (W_I g_\tau^{-1} W_I\transp)^{-1}$, so $K(\Phi_I)$
    depends only on $W_I$ and the $2^{2k}$ incident chamber metrics
    $\{g_\tau\}_{\tau\in\{-,+\}^{2k}}$.
  \item Let $g = DF\transp DF$ on $U_0$. The normal cone of $(U_0, g)$ along
    $\Phi_I \cap U_0$ is isometric to $\mathcal{K}_I = \Psi_0(\RR^{2k}) \subset
    \mathcal{T}_I^\perp$ with chamber cones $\mathsf{K}_\tau =
    \Psi_0(\overline{C_\tau})$ of metric $M_\tau$ and cells $\mathsf{K}_\varsigma :=
    \Psi_0(\overline{C_\varsigma})$ ($\varsigma \in \{-,0,+\}^{2k}$), and
    \begin{equation}\label{eq:cms-banchoff-id}
      (2\pi)^{-k} \mathcal{R}^{2k}(\Phi_I)
      = \sum_{\varsigma\in\{-,0,+\}^{2k}} (-1)^{|\varsigma|}\,\gamma\bigl(\{0\}, \mathsf{K}_\varsigma\bigr)
      = K(\Phi_I),
    \end{equation}
    so $(2\pi)^{-k}\mathcal{R}^{2k}\restriction_{\Phi_I\cap U_0} =
    K(\Phi_I)\sqrt{\det g_\perp}\,\Haus{\dimin-2k}\restriction_{\Phi_I\cap U_0}$.
\end{enumerate}
Items~2 and~3 use positive path rank only through $\rank(J_0 Q_\perp) = \dimin
- 2k$ and $M_\tau \succ 0$ for all $\tau \in \{-,+\}^{2k}$. Under these weaker
hypotheses the intrinsic normal cone of $(U_0,g)$ along $\Phi_I \cap U_0$ is
the abstract cone complex glued from the Euclidean cones
$(\overline{C_\tau}, M_\tau)$, and~\eqref{eq:cms-banchoff-id} holds with
$K(\Phi_I) = \EE_\xi[\ind_0(\xi\cdot\Psi_0)]$ given by the orthant sum
of~\eqref{eq:banchoff-orthant}, whether or not $\Psi_0$ is injective.
\end{proposition}

\begin{proof}
Since $J(s)Q_\perp = J_0 Q_\perp$, we have $X(s) = P_I J(s) V_I$ and $M(s) =
V_I\transp J(s)\transp P_I J(s) V_I$, the Schur complement of $Q_\perp\transp
g(s) Q_\perp$ in $T\transp g(s) T$ for $T := [V_I\ Q_\perp]$. The identities
$W_I V_I = I_{2k}$, $W_I Q_\perp = 0$, and $Q_\perp\transp V_I = 0$ give
$T^{-1} = [W_I;\,Q_\perp\transp]$. Hence this Schur complement equals the
inverse of the
upper-left $2k\times 2k$ block $W_I g(s)^{-1} W_I\transp$ of $(T\transp g(s)
T)^{-1} = T^{-1}g(s)^{-1}T^{-\!\top}$, proving Items~1--2 by
Theorem~\ref{thm:banchoff-orthant}.

For Item~3, replace
$\rho_\beta(u_a)$ by $(u_a)_+$ in Lemma~\ref{lem:shear-splitting}. This yields
the piecewise-linear shear $\Theta_0(u,x_\perp) := (u,\, x_\perp +
\sum_a(J_0Q_\perp)\pinv(u_aJ_0v_a + (u_a)_+a_{h_a}))$, under which
$(F\circ\iota_I\circ\Theta_0^{-1})(u,y) = c_0' + J_0Q_\perp y + \Psi_0(u)$ and the metric
on $C_\tau\times\RR^{\dimin-2k}$ splits as $M_\tau \oplus g_\perp$. The map
$\Theta_0$ is affine on each closed orthant cylinder
$\overline{C_\tau}\times\RR^{\dimin-2k}$, and it shears parallel to $T\Phi_I$
while fixing $\Phi_I$ ($u = 0$). It therefore maps each chamber of $(U_0,g)$
adjacent to $\Phi_I$ isometrically onto $\overline{C_\tau}\times\RR^{\dimin-2k}$
with the product metric $M_\tau \oplus g_\perp$, preserving the orthogonal
quotient metric $M_\tau$ on $(T\Phi_I)^{\perp_{g_\tau}}$. Hence the intrinsic
normal cone of $\Phi_I \cap U_0$ is the cone complex
$\{(\overline{C_\varsigma}, M_\tau)\}_{\varsigma}$ glued along common faces,
with tangential measure $\sqrt{\det g_\perp}\,\Haus{\dimin-2k}$. By
subdivision invariance of $\mathcal{R}^{2k}$~\cite{cheeger1984}, no simplicial
subdivision along $\RR^{\dimin-2k}$ affects~\eqref{eq:cms-ext-angle}. Each
cell $(\overline{C_\varsigma}, M_\tau)$ is isometric to the convex cone
$\mathsf{K}_\varsigma = \Psi_0(\overline{C_\varsigma})$ because $M_\tau \succ
0$, and under positive path rank Lemma~\ref{lem:pr-injectivity}(1) identifies
the whole complex with the embedded cone $\mathcal{K}_I$.

Finally, for $\xi \sim \Norm(0,P_I)$,
the projection of $\xi$ onto $\spn \mathsf{K}_\varsigma$ lies in the open polar
cone of $\mathsf{K}_\varsigma$ if and only if $0 =
\arg\max_{\overline{C_\varsigma}} \xi\cdot\Psi_0$. Hence
$\gamma(\{0\},\mathsf{K}_\varsigma) = \PP_\xi\{0 =
\arg\max_{\overline{C_\varsigma}} \xi\cdot\Psi_0\}$ and taking expectations in
\eqref{eq:banchoff-index-def} gives~\eqref{eq:cms-banchoff-id} (cf.~\cite{fu1994, mccoytropp2014}).
\end{proof}
\section{Global Weak-\texorpdfstring{$*$}{*} Convergence and Offset Invariance}\label{sec:global}

We now sum the flatwise curvatures of Section~\ref{sec:banchoff} over
$\arr(b) = \{\wall{1},\dots,\wall{\ngen}\}$. Under global genericity, the
smoothed curvature measures converge weak-$*$ to the polyhedral measure at local
rate $O(\beta^{-1})$ on compact sets, and globally in top degree
(Theorem~\ref{thm:global-weak-star}). In even dimension, a polyhedral Morse
count at infinity then gives $\sum_v K(v) = K_\infty$
(Theorem~\ref{thm:offset-invariance}), provided three conditions hold: the wall
normals span $\RR^{\dimin}$, the recession map vanishes only at $0$, and no
bounded edge contracts to a point.

\subsection{Global weak-\texorpdfstring{$*$}{*} convergence of smoothed GBC measures}\label{subsec:weak-star}

For each $J \subset \{1,\dots,\ngen\}$, define the $\beta$-tube $T_J(\beta)$ and
\emph{strict $J$-core} $T_J^\circ(\beta)$, where the active transition set
$T_\beta(x) := \{h \in \{1,\dots,\ngen\} : |z_h(x)| \le \beta^{-1}\}$ equals $J$:
\[
  T_J(\beta) := \bigl\{x \in \RR^{\dimin} : |z_h(x)| \le \beta^{-1}\ \forall h \in J\bigr\},
  \qquad
  T_J^\circ(\beta) := T_J(\beta) \setminus \bigcup_{h\notin J} T_{\{h\}}(\beta).
\]
The strict cores partition $\RR^{\dimin}$; on the $\beta$-shrunk open chambers
$T_\emptyset^\circ(\beta)$, $DF_\beta$ is locally constant and
$R_{g_\beta} \equiv 0$. Set $\Phi_\emptyset := \RR^{\dimin}$ and
$\mathcal{I}_{\ne\emptyset} := \{S \subset \{1,\dots,\ngen\} : S \ne \emptyset,\
\Phi_S := \bigcap_{h\in S}\wall{h} \ne \emptyset\}$. Under simplicity, every
$S \in \mathcal{I}_{\ne\emptyset}$ has $|S| \le \dimin$ and $\rank W_S = |S|$.
Define $\lambda^* := (\norm{L}_{\mathrm{op}} + \sum_{h=1}^{\ngen}\norm{a_h}_2\norm{w_h}_2)^2$
and, under positive path rank at every face,
\[
  \lambda_* := \min_{\mathsf{F}\in\mathcal{F}(\arr(b))}\ \min_{s\in[0,1]^{I(\mathsf{F})}}
  \lambda_{\min}\bigl(J_{\mathsf{F}}(s)\transp J_{\mathsf{F}}(s)\bigr) > 0.
\]

\begin{lemma}
\label{lem:tube-separation}
Assume $\arr(b)$ is simple. There is $\beta_\sharp = \beta_\sharp(\arr(b)) \ge 1$
such that for every $\beta \ge \beta_\sharp$ and $x \in \RR^{\dimin}$ there exist
$S^* \in \mathcal{I}_{\ne\emptyset} \cup \{\emptyset\}$ with $T_\beta(x) \subseteq S^*$
and $y \in \Phi_{S^*}$ ($y := x$ if $S^* = \emptyset$) satisfying $|z_h(x)| >
\beta^{-1}$ and $\sgn z_h(y) = \sgn z_h(x) \ne 0$ for all $h \notin S^*$. Writing
$\mathsf{F}(y)$ for the face containing $y$ ($I(\mathsf{F}(y)) = S^*$), we have
$DF_\beta(x) = J_{\mathsf{F}(y)}(s)$ with $s_h = \rho'(\beta z_h(x)) \in [0,1]$;
in particular, $|T_\beta(x)| \le \dimin$ with $\{w_h\}_{h\in T_\beta(x)}$ linearly
independent, and $\lambda_* I_{\dimin} \preceq g_\beta(x) \preceq \lambda^*
I_{\dimin}$ on $\RR^{\dimin}$.
\end{lemma}

For $2 \le 2k \le \dimin$ and $\beta \ge \beta_\sharp$, define the intrinsic
degree-$2k$ GBC measure on $(\RR^{\dimin}, g_\beta)$:
\begin{equation}\label{eq:smoothed-gbc-measure}
  \dd\mu_\beta^{(2k)}(x)
  := \frac{1}{(8\pi)^k\,k!}
     \sum_{\substack{i_1,\dots,i_{2k}=1 \\ j_1,\dots,j_{2k}=1}}^{\dimin}
     \delta^{i_1\cdots i_{2k}}_{j_1\cdots j_{2k}}\,
     \Bigl(\prod_{p=1}^k R^{j_{2p-1} j_{2p}}_{\beta,\,i_{2p-1} i_{2p}}(x)\Bigr)\,
     \sqrt{\det g_\beta(x)}\,\dd x,
\end{equation}
where $\delta^{i_1\cdots i_{2k}}_{j_1\cdots j_{2k}} =
\det(\delta^{i_a}_{j_b})_{a,b=1}^{2k}$,
$R^{cd}_{\beta,ab} = g_\beta^{ce} g_\beta^{df} R_{\beta,abef}$, and
$R_{\beta,abef}$ is the Levi-Civita Riemann tensor of $g_\beta$ (for $2k = \dimin$,
$\dd\mu_\beta^{(\dimin)}(x) = e_{\dimin}(g_\beta)(x)$). Working intrinsically on
$(\RR^{\dimin}, g_\beta)$ keeps $\mu_\beta^{(2k)}$ well defined even when
$F(\RR^{\dimin})$ self-intersects.

\begin{theorem}[Global weak-$*$ convergence of smoothed GBC measures]
\label{thm:global-weak-star}
Assume global genericity, let $\rho \in \Pprof$, and fix
$2k \in \{2,4,\dots,2\lfloor\dimin/2\rfloor\}$. For each $\Phi_I \ne \emptyset$
with $|I| = 2k$, $\Phi_I^{\mathrm{reg}} := \Phi_I \setminus \bigcup_{h\notin I}
\wall{h}$ decomposes into finitely many open $(\dimin-2k)$-cells
$\Phi_{I,\alpha}$, on each of which the base Jacobian $J_{0,I,\alpha}$ of
Definition~\ref{def:local-model}(1), the tangential metric $g_{\perp,I,\alpha} :=
(J_{0,I,\alpha} Q_{\perp,I})\transp(J_{0,I,\alpha} Q_{\perp,I})$ ($Q_{\perp,I}$ an
orthonormal basis of $\ker W_I$), and $K(\Phi_{I,\alpha})$ of
Theorem~\ref{thm:banchoff-orthant} are constant. Set
\[
  \mu_0^{(2k)}
  := \sum_{\substack{I \subset \{1,\dots,\ngen\} \\ |I|=2k,\,\Phi_I\ne\emptyset}}
     \sum_\alpha K(\Phi_{I,\alpha})\,\sqrt{\det g_{\perp,I,\alpha}}\,
     \Haus{\dimin-2k}\!\restriction_{\Phi_{I,\alpha}},
\]
which coincides on $(\RR^{\dimin}, g)$ with $(2\pi)^{-k}\mathcal{R}^{2k}$
by~\eqref{eq:cms-banchoff-id}. Then:
\begin{enumerate}
  \item For every $\beta \ge \beta_\sharp$ and $x \in \RR^{\dimin}$,
    $|T_\beta(x)| \le \dimin$ with $\{w_h\}_{h\in T_\beta(x)}$ linearly
    independent, and $\lambda_* I_{\dimin} \preceq g_\beta(x) \preceq
    \lambda^* I_{\dimin}$.
  \item When $\dimin = 2k$, each $\Phi_I \ne \emptyset$ with $|I| = \dimin$ is a
    single vertex $v_I = -W_I^{-1} b_I$. For all $\beta \ge \beta_\sharp$, the
    vertex boxes $T_I(\beta) = v_I + W_I^{-1}[-\beta^{-1},\beta^{-1}]^{\dimin}$
    are pairwise disjoint, lie in $T_I^\circ(\beta)$, contain $\supp\mu_\beta^{(\dimin)}$,
    and satisfy (writing $\int e_{\dimin}(g_\beta)$ for the unoriented density
    integral $\int \dd\mu_\beta^{(\dimin)}$)
    \begin{equation}\label{eq:vertex-box-mass}
      \mu_\beta^{(\dimin)}\bigl(T_I(\beta)\bigr)
      = \int_{T_I(\beta)} e_{\dimin}(g_\beta)
      = K(v_I),
      \qquad
      \int_{\RR^{\dimin}} e_{\dimin}(g_\beta)
      = \sum_{v\in V(b)} K(v),
    \end{equation}
    and for every Lipschitz $\varphi : \RR^{\dimin} \to \RR$,
    \[
      \left|\int_{\RR^{\dimin}} \varphi\,\dd\mu_\beta^{(\dimin)}
            - \sum_{v\in V(b)} K(v)\varphi(v)\right|
      \le \frac{C\,\operatorname{Lip}(\varphi)}{\beta}
    \]
    for a constant $C \ge 0$ independent of $\varphi$ and $\beta$.
  \item For every $R_\varphi \ge 1$ there is $C = C(R_\varphi) \ge 0$ such that
    for all $\varphi \in C_c^{0,1}(\RR^{\dimin})$ with $\supp\varphi \subset
    [-R_\varphi,R_\varphi]^{\dimin}$ and all $\beta \ge \beta_\sharp$,
    \begin{equation}\label{eq:weak-star-rate}
      \left|\int_{\RR^{\dimin}} \varphi\,\dd\mu_\beta^{(2k)}
            - \int_{\RR^{\dimin}} \varphi\,\dd\mu_0^{(2k)}\right|
      \le \frac{C\bigl(\operatorname{Lip}(\varphi) + \norm{\varphi}_\infty\bigr)}{\beta},
    \end{equation}
    where the $\norm{\varphi}_\infty$ term is absent when $2k = \dimin$ by~(2).
\end{enumerate}
\end{theorem}

Lemma~\ref{lem:tube-separation} and Theorem~\ref{thm:global-weak-star} are proved
in Appendix~\ref{app:global-proofs}.

\subsection{PL Morse theory at infinity and top-degree offset invariance}\label{subsec:offset-invariance}

The domain $\RR^{\dimin}$ is unbounded, and $F_b$ need not be globally
embedded; at multi-wall collisions it may fail even to be locally injective.
We therefore adapt polyhedral Morse
theory~\cite{banchoff1967, grigsbylindseymasden2022, grunertkuhnelrote2023} to
$\arr(b)$, using the Euler valuation $\chi_c(\mathsf{F}) = (-1)^{\dim\mathsf{F}}$
of relatively open polyhedra~\cite{mcmullen1975, schneiderweil2008}.
For fixed $(L, c_0, A, W, b)$, define $F_b$ and its positively $1$-homogeneous
\emph{recession map}:
\begin{equation}\label{eq:recession-map}
  F_b(x) := Lx + c_0 + A(Wx+b)_+,\qquad
  F_{\mathrm{rec}}(x) := Lx + A(Wx)_+ = \lim_{t\to+\infty}t^{-1}F_b(t x),
\end{equation}
so $F_b - F_{\mathrm{rec}}$ is bounded on $\RR^{\dimin}$ and $F_{\mathrm{rec}}$
is linear on each closed face $\overline{\face}$ of $\arr(0)$. Write $\mathsf{F} = \{x \in \RR^{\dimin} :
\sgn z_h(x) = \eta_h\ \forall h\} \in \mathcal{F}(\arr(b))$ ($\eta \in
\{-1,0,+1\}^{\ngen}$) for a non-empty face of $\arr(b)$ and $\overline{\mathsf{F}}$
($\dim\overline{\mathsf{F}} = \dim\mathsf{F}$) for its closure. For $\xi \in
\RR^n\setminus\{0\}$, set $\ell_\xi(x) := \ip{\xi}{F_b(x)}$ and
$\ell_\xi^{\mathrm{rec}}(x) := \ip{\xi}{F_{\mathrm{rec}}(x)}$. Banchoff's
critical-point index at $x \in \RR^{\dimin}$ and the open \emph{links at
infinity} are
\[
  \ind_x(\xi\transp F_b)
  := \sum_{\substack{\mathsf{F} \in \mathcal{F}(\arr(b)) :\, x \in \overline{\mathsf{F}} \\ \ell_\xi < \ell_\xi(x)\text{ on }\overline{\mathsf{F}}\setminus\{x\}}}
     (-1)^{\dim\mathsf{F}},
  \qquad
  L_\pm^\infty(\xi)
  := \bigl\{d \in \Sph{\dimin-1} : \pm\ell_\xi^{\mathrm{rec}}(d) > 0\bigr\}
\]
(and similarly $\ind_0(\xi\transp F_{\mathrm{rec}})$ on $\mathcal{F}(\arr(0))$).
We say that $\arr$ has \emph{spanning normals} if $\rank W = \dimin$, that $F_b$
has \emph{non-degenerate recession} if $F_{\mathrm{rec}}(d) \ne 0$ for all $d \in
\RR^{\dimin}\setminus\{0\}$, and that $F_b$ is \emph{edge non-degenerate} if
$F_b(v) \ne F_b(v')$ for the endpoints $v \ne v'$ of every bounded edge
(bounded $1$-face) of $\arr(b)$. Whenever
$F_{\mathrm{rec}} \ne 0$ on every ray ($1$-face) of $\arr(0)$ and $F_b$ is edge
non-degenerate, the finite union of proper hyperplanes
\begin{equation}\label{eq:exceptional-dirs}
  E_{\mathrm{exc}}(b)
  := \bigcup_{\substack{v,v'\text{ endpoints of a}\\ \text{bounded edge of }\arr(b)}} \bigl(F_b(v') - F_b(v)\bigr)^\perp
     \;\cup\;
     \bigcup_{\varrho\text{ ray ($1$-face) of }\arr(0)} F_{\mathrm{rec}}(\varrho)^\perp
\end{equation}
has measure zero in $\RR^n$. For $\mathsf{F} \in \mathcal{F}(\arr(b))$ with sign
vector $\eta$, write $\rec(\mathsf{F})$ for the recession cone of
$\overline{\mathsf{F}}$ and $J_{\mathsf{F}} := J_{\sigma(\eta)}$ with
$\sigma(\eta)_h := \mathbf{1}\{\eta_h = +1\}$ in~\eqref{eq:chamber-jacobian}.
Since $w_h\transp d = 0$ when $\eta_h = 0$ and $d \in \rec(\mathsf{F})$, the
vector $J_{\mathsf{F}} d$ does not depend on the zero-coordinate convention. At $v \in
V(b)$, write $F_v(d) := \lim_{t\downarrow 0}t^{-1}(F_b(v+td)-F_b(v))$ for the
tangent cone map.

\begin{lemma}
\label{lem:pl-morse-infinity}
In the setting of~\eqref{eq:recession-map}--\eqref{eq:exceptional-dirs}:
\begin{enumerate}
  \item Assume spanning normals. For every $\mathsf{F} \in \mathcal{F}(\arr(b))$
    with sign vector $\eta \in \{-1,0,+1\}^{\ngen}$, $\overline{\mathsf{F}}$ is
    pointed, $\rec(\mathsf{F}) = \{d \in \RR^{\dimin} : \eta_h w_h\transp d \ge
    0\ \forall h,\ w_h\transp d = 0\text{ if }\eta_h = 0\}$ is a union of closed
    faces of $\arr(0)$ whose extreme rays are rays of $\arr(0)$, and
    $F_{\mathrm{rec}}(d) = J_{\mathsf{F}} d$ and $\nu_1(d) :=
    \sum_{h=1}^{\ngen}|w_h\transp d| = \sum_{h:\,\eta_h\ne 0}\eta_h w_h\transp d$
    are linear on $\rec(\mathsf{F})$. In particular, $\ind_x(\xi\transp F_b) = 0$
    for every $\xi \in \RR^n$ and $x \notin V(b)$. Suppose further that
    $F_{\mathrm{rec}} \ne 0$ on every ray of $\arr(0)$, edge non-degeneracy holds,
    and $\xi \notin E_{\mathrm{exc}}(b)$. Then $\ell_\xi$ is non-constant on every
    edge of $\arr(b)$, and for each $\mathsf{F} \in \mathcal{F}(\arr(b))$ the
    following are equivalent: (i)~$\ell_\xi$ is bounded above on
    $\overline{\mathsf{F}}$; (ii)~$\ell_\xi^{\mathrm{rec}} < 0$ on
    $\rec(\mathsf{F})\setminus\{0\}$; (iii)~$\ell_\xi$ attains its supremum on
    $\overline{\mathsf{F}}$ at a unique point $v^*(\mathsf{F}) \in V(b)$.
  \item Assume spanning normals, edge non-degeneracy, and $F_{\mathrm{rec}} \ne 0$
    on every ray of $\arr(0)$. Then for every $\xi \in \RR^n \setminus
    E_{\mathrm{exc}}(b)$, $L_\pm^\infty(\xi)$ have the homotopy type of finite CW
    complexes and
    \begin{equation}\label{eq:pl-morse-balance}
      \sum_{v \in V(b)} \ind_v(\xi\transp F_b)
      = (-1)^{\dimin}\bigl(1-\chi(L_+^\infty(\xi))\bigr)
      = 1-\chi(L_-^\infty(\xi))
      = \ind_0(\xi\transp F_{\mathrm{rec}}),
    \end{equation}
    which equals $\sum_{\mathsf{F} \in \mathcal{F}(\arr(b)) : \ell_\xi\text{ bdd.\ above on }\overline{\mathsf{F}}} (-1)^{\dim\mathsf{F}}$
    and depends only on $(L, A, W)$ and $\xi$.
  \item Global genericity at $b$ implies non-degenerate recession and edge
    non-degeneracy. If also $\rank W = \dimin$ and $\dimin = 2k$ is even, then
    for every $v \in V(b)$ and $\xi \notin E_{\mathrm{exc}}(b)$,
    $\ind_v(\xi\transp F_b) = \prod_{a=1}^{\dimin}\varepsilon_a^{(v)}(\xi)$ for
    the transverse cone $\Psi_0^{(v)}(u) = F_v(W_{I(v)}^{-1}u)$
    (Lemma~\ref{lem:lower-link-index}(1),(3)).
\end{enumerate}
\end{lemma}

\begin{proof}
\textit{(1).}\
Since $\ker W = \{0\}$, $\overline{\mathsf{F}} = \operatorname{conv}(\overline{\mathsf{F}}
\cap V(b)) + \rec(\mathsf{F})$ is pointed, and $\rec(\mathsf{F})$ has the stated
form and is a union of closed faces $\overline{\face}$ of $\arr(0)$ on which
$F_{\mathrm{rec}}(d) = J_{\mathsf{F}} d$ and $\nu_1(d)$ are
linear~\cite[Thms.~8.3 and~18.5]{rockafellar1970}. Any $x \notin V(b)$ lies in a face $\mathsf{F}_0
\in \mathcal{F}(\arr(b))$ of dimension $\ge 1$ on which $\ell_\xi$ is affine, so
$\ell_\xi$ cannot have a strict maximum at $x$ on any $\overline{\mathsf{F}} \ni
x$ ($\mathsf{F}_0 \subseteq \overline{\mathsf{F}}$), giving $\ind_x(\xi\transp
F_b) = 0$.

For $\xi \notin E_{\mathrm{exc}}(b)$, $\ell_\xi$ is non-constant on bounded edges
and $\ell_\xi^{\mathrm{rec}} \ne 0$ on rays of $\arr(0)$. Thus $\ell_\xi(p + d) =
\ell_\xi(p) + \ell_\xi^{\mathrm{rec}}(d)$ is bounded above on
$\overline{\mathsf{F}}$ if and only if $\ell_\xi^{\mathrm{rec}} \le 0$ on
$\rec(\mathsf{F})$, which is equivalent to $\ell_\xi^{\mathrm{rec}} < 0$ on
$\rec(\mathsf{F})\setminus\{0\}$. Indeed, any non-zero $d \in \rec(\mathsf{F})$
with $\ell_\xi^{\mathrm{rec}}(d) = 0$ would lie in a closed face $\overline{\face}
\subseteq \rec(\mathsf{F})$ of $\arr(0)$ on which $\ell_\xi^{\mathrm{rec}}$ is
linear and non-zero on extreme rays, forcing $\ell_\xi^{\mathrm{rec}} > 0$ on at
least one of them. When $\ell_\xi$ is bounded above on $\overline{\mathsf{F}}$,
its argmax face $\mathsf{G} \subseteq \overline{\mathsf{F}}$ is non-empty and
bounded since $\rec(\mathsf{G}) = \{0\}$, and it contains no bounded edge because
$\xi \notin E_{\mathrm{exc}}(b)$; hence $\mathsf{G} = \{v^*(\mathsf{F})\}
\subset V(b)$ and $\sum_{v \in V(b)} \ind_v(\xi\transp F_b) =
\sum_{\mathsf{F}:\,\ell_\xi\text{ bdd.\ above on }\overline{\mathsf{F}}}
(-1)^{\dim\mathsf{F}}$.

\ppart{(2).}
To relate the bounded-above cell sum from~(1) to the links at infinity we
proceed in three steps: (i)~replace $\Sph{\dimin-1}$ by a polytopal sphere on
which $\ell_\xi^{\mathrm{rec}}$ is affine on every cell; (ii)~retract each
positive-recession locus onto a finite subcomplex, both for $\RR^{\dimin}$ and
for the recession cones $\rec(\mathsf{F})$; (iii)~exchange the order of
summation and evaluate the inner sums with the Euler valuation.

For~(i), because $\nu_1(d) = \sum_{h=1}^{\ngen}|w_h\transp d|$ is a norm on
$\RR^{\dimin}$ ($\ker W = \{0\}$), radial projection $d \mapsto d/\norm{d}_2$
maps the polytopal sphere $\Sph{\dimin-1}_W := \{d \in \RR^{\dimin} : \nu_1(d)
= 1\}$ homeomorphically onto $\Sph{\dimin-1}$ and identifies
$\{\ell_\xi^{\mathrm{rec}} \gtrless 0\} \cap \Sph{\dimin-1}_W$ with
$L_\pm^\infty(\xi)$. Each non-zero face $\face^\circ \in
\mathcal{F}(\arr(0))\setminus\{\{0\}\}$ meets $\Sph{\dimin-1}_W$ in a compact
polytope $P_\face := \overline{\face}\cap\Sph{\dimin-1}_W$ of dimension
$\dim\face^\circ - 1$, on which $\ell_\xi^{\mathrm{rec}}$ is affine and non-zero
at every vertex since $\xi \notin E_{\mathrm{exc}}(b)$. For any closed cone
$C \subseteq \RR^{\dimin}$ that is a union of closed faces of $\arr(0)$, define
the finite subcomplex
\[
  \Delta_+^\infty(C;\xi)
  := \bigcup\bigl\{P_\face : \overline{\face} \subseteq C,\ \face^\circ \ne \{0\},\
     \ell_\xi^{\mathrm{rec}} > 0\text{ on }P_\face\bigr\},
\]
and write $\mathcal{F}_+^\infty(\xi) := \{\face^\circ \in
\mathcal{F}(\arr(0))\setminus\{\{0\}\} : \ell_\xi^{\mathrm{rec}} > 0\text{ on
}P_\face\}$ for the cell set at $C = \RR^{\dimin}$. Since positivity of
$\ell_\xi^{\mathrm{rec}}$ on a cell implies positivity on its faces,
$\Delta_+^\infty(C;\xi)$ is indeed a subcomplex.

For~(ii), consider any $j$-cell $P_\face \subset C\cap\Sph{\dimin-1}_W$ ($j \ge
1$) not belonging to $\Delta_+^\infty(C;\xi)$ that meets
$\{\ell_\xi^{\mathrm{rec}} \ge 0\}$. Because $\ell_\xi^{\mathrm{rec}}$ is affine
on $P_\face$ and fails to be positive everywhere on $P_\face$, we have
$\min_{P_\face}\ell_\xi^{\mathrm{rec}} < 0$, so we may choose $p_\face \in
\operatorname{ri}(P_\face)$ with $\ell_\xi^{\mathrm{rec}}(p_\face) < 0$. Every
$x \in P_\face\setminus\{p_\face\}$ decomposes uniquely via the Minkowski gauge
$s \in (0,1]$ of $P_\face - p_\face$ as $x = (1-s)p_\face + s q$ with $q \in
\partial P_\face$. The radial homotopy
\[
  R_\face(x,t) := (1-t)(1-s)p_\face + \bigl(s + t(1-s)\bigr)q,
  \qquad t \in [0,1],
\]
is a strong deformation retraction of $P_\face\setminus\{p_\face\}$ onto
$\partial P_\face$ fixing $\partial P_\face$ pointwise. If $s = 1$, then $x =
q \in \partial P_\face$ and $R_\face(q,t) = q$ for all $t \in [0,1]$; if $s \in
(0,1)$, then $\ell_\xi^{\mathrm{rec}}(x) = (1-s)\ell_\xi^{\mathrm{rec}}(p_\face)
+ s\ell_\xi^{\mathrm{rec}}(q)$ with $\ell_\xi^{\mathrm{rec}}(p_\face) < 0$, so
$\ell_\xi^{\mathrm{rec}}(x) \ge 0$ forces $\ell_\xi^{\mathrm{rec}}(q) > 0$, and
$\ell_\xi^{\mathrm{rec}}(R_\face(x,t)) = (1-t)\ell_\xi^{\mathrm{rec}}(x) +
t\ell_\xi^{\mathrm{rec}}(q)$ shows that $R_\face(\cdot,t)$ preserves both
$\{\ell_\xi^{\mathrm{rec}} > 0\}$ and $\{\ell_\xi^{\mathrm{rec}} \ge 0\}$. We
apply these retractions in decreasing dimension $j = \dimin-1,\dots,1$. Within
a fixed dimension the cells to be retracted have pairwise disjoint relative
interiors. They meet only along their boundaries, where each $R_\face(\cdot,t)$
is the identity, so the cellwise retractions glue continuously into a single
homotopy. After stage $j$, the remaining set lies in the union of
$\Delta_+^\infty(C;\xi)$ with the $(j-1)$-skeleton. The composite is therefore
a well-defined strong deformation retraction of both
$\{\ell_\xi^{\mathrm{rec}} > 0\}\cap C\cap\Sph{\dimin-1}_W$ and
$\{\ell_\xi^{\mathrm{rec}} \ge 0\}\cap C\cap\Sph{\dimin-1}_W$ onto the finite CW
complex $\Delta_+^\infty(C;\xi)$, whose Euler characteristic is
\[
  \chi\bigl(\Delta_+^\infty(C;\xi)\bigr)
  = \sum_{\substack{\face^\circ \in \mathcal{F}_+^\infty(\xi) \\ \overline{\face} \subseteq C}}
    (-1)^{\dim\face^\circ - 1}.
\]
At $C = \RR^{\dimin}$, this shows that $L_+^\infty(\xi)$ has the homotopy type
of $\Delta_+^\infty(\RR^{\dimin};\xi)$, so $\chi(L_+^\infty(\xi)) =
\chi(\Delta_+^\infty(\RR^{\dimin};\xi))$. Replacing $\xi$ by $-\xi$ gives
$\ind_0(\xi\transp F_{\mathrm{rec}}) = 1 - \chi(\Delta_+^\infty(\RR^{\dimin};
-\xi)) = 1 - \chi(L_-^\infty(\xi))$. Now take $C = \rec(\mathsf{F})$ for a
face $\mathsf{F} \in \mathcal{F}(\arr(b))$. By part~(1), both
$\ell_\xi^{\mathrm{rec}}$ and $\nu_1$ are linear on $\rec(\mathsf{F})$. The
intersection $\{\ell_\xi^{\mathrm{rec}} \ge 0\}\cap\rec(\mathsf{F})\cap
\Sph{\dimin-1}_W$ is thus a single compact convex polytope: it is empty if
$\ell_\xi^{\mathrm{rec}} < 0$ on $\rec(\mathsf{F})\setminus\{0\}$, and
contractible (hence of Euler characteristic~$1$) otherwise. The retraction onto
$\Delta_+^\infty(\rec(\mathsf{F});\xi)$ then gives the indicator identity
\[
  \mathbf{1}\bigl\{\rec(\mathsf{F})\setminus\{0\} \subset \{\ell_\xi^{\mathrm{rec}} < 0\}\bigr\}
  = 1 - \chi\bigl(\Delta_+^\infty(\rec(\mathsf{F});\xi)\bigr).
\]

For~(iii), substituting this indicator into the bounded-above cell sum
from~(1) and exchanging the finite sums yields
\[
  \sum_{v\in V(b)} \ind_v(\xi\transp F_b)
  = \sum_{\mathsf{F}\in\mathcal{F}(\arr(b))} (-1)^{\dim\mathsf{F}}
    - \sum_{\face^\circ\in\mathcal{F}_+^\infty(\xi)} (-1)^{\dim\face^\circ-1}
      \sum_{\substack{\mathsf{F}\in\mathcal{F}(\arr(b)) \\ \overline{\face}\subseteq\rec(\mathsf{F})}} (-1)^{\dim\mathsf{F}}.
\]
Fix $\face^\circ \in \mathcal{F}_+^\infty(\xi)$ and pick $d \in \face^\circ =
\operatorname{ri}(\overline{\face})\setminus\{0\}$. By~(1), a face $\mathsf{F}
\in \mathcal{F}(\arr(b))$ satisfies $\overline{\face} \subseteq
\rec(\mathsf{F})$ if and only if $d \in \rec(\mathsf{F})$, which holds if and
only if $\mathsf{F}$ is contained in the open polyhedron
\[
  O_d := \bigcap_{h:\,w_h\transp d > 0}\{x \in \RR^{\dimin} : z_h(x) > 0\}
         \;\cap\;
         \bigcap_{h:\,w_h\transp d < 0}\{x \in \RR^{\dimin} : z_h(x) < 0\}.
\]
For all sufficiently large $t > 0$ we have $t d \in O_d$. Hence $O_d$ is a
non-empty open convex polyhedron, partitioned into the faces $\mathsf{F}$ with
$\overline{\face} \subseteq \rec(\mathsf{F})$. Finite additivity of the Euler
valuation $\chi_c(\mathsf{F}) = (-1)^{\dim\mathsf{F}}$~\cite{mcmullen1975,
schneiderweil2008} gives $\sum_{\mathsf{F}:\,\overline{\face}\subseteq
\rec(\mathsf{F})}(-1)^{\dim\mathsf{F}} = \chi_c(O_d) = (-1)^{\dimin}$ and
$\sum_{\mathsf{F}\in\mathcal{F}(\arr(b))}(-1)^{\dim\mathsf{F}} =
\chi_c(\RR^{\dimin}) = (-1)^{\dimin}$, whence $\sum_{v\in V(b)}\ind_v(\xi\transp
F_b) = (-1)^{\dimin}(1 - \chi(L_+^\infty(\xi)))$. Now specialise this identity
to $b = 0$, where $E_{\mathrm{exc}}(0) \subseteq E_{\mathrm{exc}}(b)$ and $V(0) =
\{0\}$. It equates $(-1)^{\dimin}(1 - \chi(L_+^\infty(\xi)))$ with
$\ind_0(\xi\transp F_{\mathrm{rec}}) = 1 - \chi(L_-^\infty(\xi))$, which
completes~\eqref{eq:pl-morse-balance}. (The resulting relation
$(-1)^{\dimin}(1-\chi(L_+^\infty(\xi))) = 1 - \chi(L_-^\infty(\xi))$ is an
instance of Alexander duality on $\Sph{\dimin-1}$, the two links being
separated by the zero set of $\ell_\xi^{\mathrm{rec}}$.)

\ppart{(3).}
Given $d \in \RR^{\dimin}\setminus\{0\}$, pick $\mathsf{F} \in
\mathcal{F}(\arr(b))$ with $d \in \rec(\mathsf{F})$ and $x \in \mathsf{F}$. Since
$x + td \in \mathsf{F}$ for all $t \ge 0$, $w_h\transp d = 0$ for $h \in
I(\mathsf{F})$ and $\mathbf{1}\{z_h(x) > 0\} = \mathbf{1}\{w_h\transp d > 0\}$
whenever $w_h\transp d \ne 0$, so $J_{\mathsf{F}}(s)d = F_{\mathrm{rec}}(d) \ne 0$
by positive path rank at $\mathsf{F}$. Along a bounded edge $\mathsf{F}$ with
endpoints $v, v'$, $W_{I(\mathsf{F})}(v'-v) = 0$ likewise gives $F_b(v') - F_b(v) =
J_{\mathsf{F}}(s)(v'-v) \ne 0$. At $v \in V(b)$ with $u = W_{I(v)}(x - v)$,
$\xi\transp(F_b(x) - F_b(v)) = \xi\cdot\Psi_0^{(v)}(u)$ has edge directions $\pm
W_{I(v)}^{-1}e_a$, so $\{\ell_a^\pm(\xi) = 0\} \subset E_{\mathrm{exc}}(b)$ and
Lemma~\ref{lem:lower-link-index}(1),(3) yields~(3).
\end{proof}

At any vertex $v \in V(b)$, including non-simple multi-wall collisions, we
define Banchoff's vertex curvature on $\arr(b)$ as
\[
  K(v)
  := \EE_{\xi\sim\Norm(0,I_n)}\bigl[\ind_v(\xi\transp F_b)\bigr]
   = \EE_{\xi\sim\mathrm{Unif}(\Sph{n-1})}\bigl[\ind_v(\xi\transp F_b)\bigr],
\]
By Proposition~\ref{prop:schur-cms}(3), this agrees with the intrinsic density
$(2\pi)^{-\dimin/2}\mathcal{R}^{\dimin}(v)$ at every locally simple vertex $v$
with $\rank J_\sigma = \dimin$ on the incident chambers. This covers every
vertex at which $F_v$ embeds a neighbourhood of $v$
(cf.~\eqref{eq:cms-ext-angle}).

\begin{theorem}[Top-degree offset invariance]
\label{thm:offset-invariance}
Fix $(L, A, W)$ with spanning normals and non-degenerate recession, and set
\[
  K_\infty
  := \EE_{\xi\sim\Norm(0,I_n)}\bigl[\ind_0(\xi\transp F_{\mathrm{rec}})\bigr]
   = \sum_{\face^\circ\in\mathcal{F}(\arr(0))} (-1)^{\dim\face^\circ}\,\gamma_{\face}
\]
(also equal to $\EE_{\xi\sim\Norm(0,I_n)}[1 - \chi(L_-^\infty(\xi))]$
by~\eqref{eq:pl-morse-balance}), where $\gamma_0 := 1$ for $\face^\circ = \{0\}$,
and for $\face^\circ \ne \{0\}$ with extreme-ray generators
$d_{\face,1},\dots,d_{\face,m_{\face}}$ and unit image rays $v_{\face,j} :=
F_{\mathrm{rec}}(d_{\face,j})/\norm{F_{\mathrm{rec}}(d_{\face,j})}_2$,
$\gamma_{\face} := \PP_{\xi\sim\Norm(0,I_n)}(\xi\transp v_{\face,j} < 0\
\forall\,j \le m_{\face})$ is the normalised exterior solid angle of
$F_{\mathrm{rec}}(\overline{\face})$~\cite{mccoytropp2014, schneiderweil2008}. Then:
\begin{enumerate}
  \item In even dimension $\dimin = 2k$, for every $b \in \RR^{\ngen}$
    satisfying edge non-degeneracy (including collision offsets on $\mathcal{D}$
    in Remark~\ref{rem:circuit}), $\sum_{v \in V(b)} K(v) = K_\infty$ depends only
    on $(L, A, W)$, and $K(v)$ is given by the orthant-probability sum
    in~\eqref{eq:banchoff-orthant} at every locally simple vertex $v$
    ($\dimin$ active walls and $\rank J_\sigma = \dimin$
    on the $2^{\dimin}$ incident chambers). Under global genericity at $b$, every
    vertex is locally simple and, for all $\rho \in \Pprof$ and $\beta \ge
    \beta_\sharp(b)$,
    \begin{equation}\label{eq:offset-invariance}
      \sum_{v \in V(b)} K(v)
      = \int_{\RR^{\dimin}} e_{\dimin}\bigl(g_\beta(\cdot;b)\bigr)
      = \lim_{\beta\to\infty} \int_{\RR^{\dimin}} e_{\dimin}\bigl(g_\beta(\cdot;b)\bigr)
      = K_\infty.
    \end{equation}
  \item If $\dimin$ is odd, $\ind_0(-\xi\transp F_{\mathrm{rec}}) =
    -\ind_0(\xi\transp F_{\mathrm{rec}})$ a.e., so $K_\infty = 0$; if $b$ also
    satisfies edge non-degeneracy, $\ind_v(-\xi\transp F_b) = -\ind_v(\xi\transp
    F_b)$ for all $v \in V(b)$ and $\xi \notin E_{\mathrm{exc}}(b)$, so $K(v) = 0$
    for all $v \in V(b)$ (see Remark~\ref{rem:morse-hyp-sharpness}).
  \item When $\dimin = 2$, let $\arr(0)$ have $m_{\mathrm{rays}} \in [4, 2\ngen]$
    distinct rays dividing $\RR^2$ into pointed sectors
    $\psector_1,\dots,\psector_{m_{\mathrm{rays}}}$ of image angle
    $\theta(\psector_j) \in [0,\pi)$ under $F_{\mathrm{rec}}$. Then
    \begin{equation}\label{eq:planar-defect}
      K_\infty
      = 1 - \frac{1}{2\pi}\sum_{j=1}^{m_{\mathrm{rays}}} \theta(\psector_j)
      = 1 - \frac{m_{\mathrm{rays}}}{2}
        + \sum_{j=1}^{m_{\mathrm{rays}}}\Bigl(\frac{1}{2} - \frac{\theta(\psector_j)}{2\pi}\Bigr).
    \end{equation}
    For every edge-non-degenerate $b$, $2\pi\sum_{v\in V(b)} K(v) = 2\pi K_\infty
    = 2\pi\chi(\RR^2) - \Theta_\infty$ with total cone angle $\Theta_\infty :=
    \sum_{j} \theta(\psector_j)$ at infinity~\cite{cohnvossen1935, huber1957}, and
    each $v \in V(b)$ with incident sectors $\psector_{v,1},\dots,\psector_{v,m_v}$
    of tangent angles $\theta(\psector_{v,j}) \in [0,\pi]$ under $F_v$ has
    $K(v) = 1 - \frac{1}{2\pi}\sum_{j=1}^{m_v}\theta(\psector_{v,j})$.
\end{enumerate}
\end{theorem}

\begin{proof}
\textit{(1).}\
Taking $\EE_{\xi\sim\Norm(0,I_n)}$ in~\eqref{eq:pl-morse-balance} gives
$\sum_{v\in V(b)} K(v) = K_\infty$. At a locally simple vertex $v$, setting $u =
W_{I(v)}(x - v)$ identifies $\ind_v(\xi\transp F_b) = \ind_0(\xi\cdot\Psi_0)$ for
$\Psi_0(u) := F_v(W_{I(v)}^{-1}u)$, and Lemma~\ref{lem:pr-injectivity}(2) with
Theorem~\ref{thm:banchoff-orthant} (whose last two equalities require only
$M_\tau = W_{I(v)}^{-\!\top} g_\tau W_{I(v)}^{-1} \succ 0$)
yields~\eqref{eq:banchoff-orthant}. Under global genericity,
Lemma~\ref{lem:pl-morse-infinity}(3) and Theorem~\ref{thm:global-weak-star}(2)
give~\eqref{eq:offset-invariance}.

\ppart{(2).}
At $b = 0$ ($\arr(0)$ has no bounded edges and $V(0) = \{0\}$),
Lemma~\ref{lem:pl-morse-infinity}(2) for $\xi \notin E_{\mathrm{exc}}(0) =
-E_{\mathrm{exc}}(0)$ gives $\ind_0(\xi\transp F_{\mathrm{rec}}) = 1 -
\chi(L^\infty_-(\xi)) = (-1)^{\dimin}(1 - \chi(L^\infty_+(\xi)))$, and since
$L^\infty_\mp(\xi) = L^\infty_\pm(-\xi)$, $\ind_0(-\xi\transp F_{\mathrm{rec}}) =
(-1)^{\dimin}\ind_0(\xi\transp F_{\mathrm{rec}})$, forcing $K_\infty = 0$ for odd
$\dimin$.

The vertex case is similar. Assume edge non-degeneracy, fix $v \in V(b)$, and
let $\mathcal{A}_v(0)$ be the central arrangement of walls through $v$
translated to $0$. Then $\ind_v(\xi\transp F_b) = \ind_0(\xi\transp F_v)$, and
$(\mathcal{A}_v(0), F_v)$ has spanning normals and no bounded edges. Moreover
$F_v \ne 0$ on every ray $\varrho$ of $\mathcal{A}_v(0)$: either $\varrho$
directs an unbounded edge, where $F_v = F_{\mathrm{rec}} \ne 0$, or it directs a
bounded edge $[v,v']$, where $F_v(v'-v) = F_b(v') - F_b(v) \ne 0$. Applying
Lemma~\ref{lem:pl-morse-infinity}(2) to $(\mathcal{A}_v(0),F_v)$ for $\xi \notin
\bigcup_\varrho F_v(\varrho)^\perp \subseteq E_{\mathrm{exc}}(b)$ gives
$\ind_v(-\xi\transp F_b) = (-1)^{\dimin}\ind_v(\xi\transp F_b)$, hence $K(v) = 0$
for odd $\dimin$.

\ppart{(3).}
When $\dimin = 2$, $\gamma_0 = 1$, $\gamma_{\varrho_j} = \tfrac{1}{2}$ on each ray,
and $\gamma_{\psector_j} = \tfrac{1}{2} - \theta(\psector_j)/(2\pi)$ by Sheppard's
formula, giving~\eqref{eq:planar-defect}. Applying the same count to
$(\mathcal{A}_v(0),F_v)$, where Sheppard's formula remains valid at the boundary
angles $\theta = 0$ ($\gamma = \tfrac12$) and $\theta = \pi$ ($\gamma = 0$), yields
the formula for $K(v)$.
\end{proof}

\begin{remk}\label{rem:morse-hyp-sharpness}
By Lemma~\ref{lem:pl-morse-infinity}(2), the balance $\sum_{v\in V(b)} K(v) =
\EE_\xi[\ind_0(\xi\transp F_{\mathrm{rec}})]$ holds already whenever
$F_{\mathrm{rec}} \ne 0$ on every ray of $\arr(0)$ and edge non-degeneracy
holds. In~\eqref{eq:planar-defect}, parallel walls contribute coincident rays of
$\arr(0)$ once. Full non-degenerate recession is used only to exclude antipodal
ray images on $\psector_j$, which would force an interior zero; this keeps
$\theta(\psector_j) \in [0,\pi)$. At $v \in V(b)$, by contrast, an interior
zero of $F_v$ is permitted and $\theta(\psector_{v,j})$ can attain $\pi$.

Edge non-degeneracy matters. If it fails and a bounded edge with endpoints
$v, v'$ contracts to a point, $F_v$ vanishes along that ray, and both
Theorem~\ref{thm:offset-invariance}(1) and~(2) can fail. For instance,
$g(t) := t - (t)_+ + (t-1)_+$ maps $[0,1]$ to $0$ with $g_{\mathrm{rec}}(t) = t$,
giving $K(0) = K(1) = \tfrac{1}{2} \ne 0$ in $\dimin = 1$. For $F_b(x_1,x_2) =
(g(x_1), g(x_2))$ on $\RR^2$ ($F_{\mathrm{rec}} = \mathrm{id}_{\RR^2}$,
$K_\infty = 0$) we get $K(v) = \tfrac{1}{4}$ on $\{0,1\}^2$, so
$\sum_{v\in V(b)} K(v) = 1 \ne K_\infty$.
\end{remk}

\begin{figure}[t]
\centering
\includegraphics[width=\textwidth]{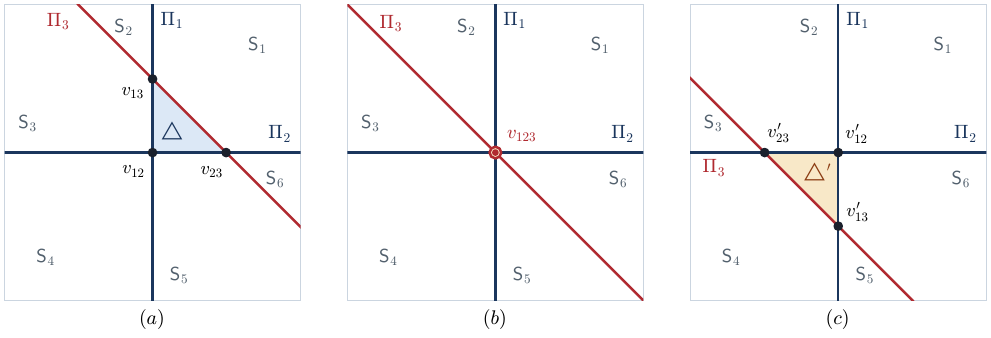}
\caption{Crossing the circuit wall $\mathcal{D} = \{b_1+b_2-b_3=0\}$ for $\dimin = 2$, $\ngen = 3$, $w_1 = (1,0)$, $w_2 = (0,1)$, $w_3 = (1,1)$ (Remark~\ref{rem:circuit}). \emph{(a)}~$b_1+b_2-b_3>0$ and \emph{(c)}~$b_1+b_2-b_3<0$: seven-chamber offsets around a bounded triangle ($\triangle$, $\triangle'$). \emph{(b)}~$b_1+b_2-b_3=0$: triple-line collision $v_{123}$. In all panels the unbounded sectors $\psector_1,\dots,\psector_6$ are fixed and $\sum_v K(v) = K_\infty = 1 - \Theta_\infty/(2\pi)$.}
\label{fig:circuit}
\end{figure}

\begin{remk}\label{rem:circuit}
Theorem~\ref{thm:offset-invariance}(1) holds pointwise in $b$ without deformation,
even though the globally generic locus is disconnected. Indeed, $\arr(b)$ is
non-simple if and only if $b_C \in \im W_C =: \Delta_C$ for some circuit $C$ of
the row matroid of $W$ ($\rank W_C = |C| - 1$). Simple offsets therefore form
the complement of the \emph{circuit} (or \emph{discriminantal})
\emph{arrangement} $\mathcal{D} := \{\Delta_C\}$~\cite{maninschechtman1989},
and Theorem~\ref{thm:offset-invariance}(1) remains valid on its collision walls.
In Figure~\ref{fig:circuit} we have $w_1 + w_2 - w_3 = 0$ and $\mathcal{D} =
\{b_1 + b_2 - b_3 = 0\}$. Suppose $\rank DF_b = 2$ on all seven chambers of
Panel~(a). The interior angles of $F_b(\triangle)$ sum to $\pi$, the three
chambers sharing an edge with $\triangle$ contribute
$\pi + \theta(\psector_{2j-1})$, and the three corner chambers contribute
$\theta(\psector_{2j})$. Around $\triangle'$ in Panel~(c) the roles of
$\psector_{2j-1}$ and $\psector_{2j}$ are reversed. In both panels the twelve
vertex angles thus sum to $4\pi + \Theta_\infty$, and
$\sum_{v\in V(b)} K(v) = 3 - (4\pi + \Theta_\infty)/(2\pi) =
1 - \Theta_\infty/(2\pi) = K(v_{123}) = K_\infty$.
\end{remk}

\subsection{Convex piecewise-affine graphs}
\label{subsec:convex-graphs}

For the graph of a convex piecewise-affine function in even dimension with
spanning normals, $h_{2k}(s) > 0$ pointwise and Theorem~\ref{thm:offset-invariance}(1)
recovers the Gauss--Bonnet evaluation induced by the subdifferential tiling of
the Newton zonotope, a regular zonotopal tiling in the sense of the
Bohne--Dress theorem~\cite{richtergebertziegler1994}, \cite[\S7.5]{ziegler1995}.

\begin{corollary}[Zonotopal subdifferential tilings of convex graphs]
\label{cor:conic-zonotope}
In even dimension $\dimin = 2k \ge 2$, let $f(x) = \ell_0\transp x + c_0 +
\sum_{h=1}^{\ngen} \alpha_h(w_h\transp x + b_h)_+$ ($\alpha_h > 0$) have $\rank W
= \dimin$ and simple $\arr(b)$, and let $F(x) := (x, f(x)) \in \RR^{\dimin+1}$
(positive path rank holds at every face since $DF = [I_{\dimin};\, \nabla f\transp]$).
At $v \in V(b)$ with $I(v) = (h_1,\dots,h_{2k})$, write $I_+(v) := \{h \notin
I(v) : w_h\transp v + b_h > 0\}$ and $g_v(s) := \ell_0 +
\sum_{h\in I_+(v)}\alpha_h w_h + \sum_{a=1}^{2k} s_a \alpha_{h_a} w_{h_a} \in
\partial f(v)$. Then for all $s \in [0,1]^{2k}$,
\[
  h_{2k}(s)
  = \frac{(2k-1)!!\Bigl(\prod_{a=1}^{2k}\alpha_{h_a}\Bigr)\abs{\det W_{I(v)}}}
         {(2\pi)^k\,\bigl(1 + \norm{g_v(s)}_2^2\bigr)^{k+1/2}}
  > 0,
\]
and under $\nu_+(g) := (-g, 1)/\sqrt{1+\norm{g}_2^2} \in S_+^{2k}$, the
parallelepipeds $\{\partial f(v)\}_{v\in V(b)}$ tile the Newton zonotope
$\mathcal{Z}(f) := \ell_0 + \sum_{h=1}^{\ngen}[0,\alpha_h]w_h =
\partial f_{\mathrm{rec}}(0)$ with disjoint interiors and
\begin{align}\label{eq:convex-gnomonic}
  K(v)
  &= \frac{2\,\Haus{2k}\bigl(\nu_+(\partial f(v))\bigr)}{\abs{\Sph{2k}}} > 0,\\
  \label{eq:zonotope-conservation}
  \sum_{v\in V(b)} K(v)
  &= K_\infty
   = \frac{2\,\Haus{\dimin}(\nu_+(\mathcal{Z}(f)))}{\abs{\Sph{\dimin}}}
   \in (0, 1).
\end{align}
\end{corollary}

\begin{proof}
With $J(s) = [I_{2k};\,g_v(s)\transp]$, $V_I = W_{I(v)}^{-1}$, and upward unit
normal $\nu(s) = (1+\norm{g_v(s)}_2^2)^{-1/2}(-g_v(s), 1)\transp$, one has
$(\det M(s))^{-1/2} = \abs{\det W_{I(v)}}/\sqrt{1+\norm{g_v(s)}_2^2}$ and
$N_{ab}(s) = \alpha_{h_a}\alpha_{h_b}/(1+\norm{g_v(s)}_2^2)$, giving $h_{2k}(s)$
and, via $s \mapsto g_v(s)$ and $\abs{\Sph{2k}} = 2(2\pi)^k/(2k-1)!!$,
\eqref{eq:convex-gnomonic}.

For the tiling, Moreau--Rockafellar
additivity~\cite[Theorem~23.8]{rockafellar1970} gives $\partial f(x) \subseteq
\mathcal{Z}(f)$. Every $g_0 \in \operatorname{int}\mathcal{Z}(f) =
\operatorname{int}\partial f_{\mathrm{rec}}(0)$ satisfies $f_{\mathrm{rec}}(d) -
g_0\transp d > 0$ on $\Sph{\dimin-1}$, so $f(x) - g_0\transp x$ attains its
minimum at some $x_*$ with $g_0 \in \partial f(x_*)$; since simplicity gives
$\dim\partial f(x) \le |I(\mathsf{F}(x))| < \dimin$ off $V(b)$, the compact
parallelepipeds $\{\partial f(v)\}_{v\in V(b)}$ cover $\mathcal{Z}(f)$. Finally,
for $g = g_v(s) \in \operatorname{int}\partial f(v)$ ($s \in (0,1)^{2k}$) and
$v' \ne v$, invertibility of $W_{I(v)}$ gives $f(v') - f(v) \ge f'(v;v'-v) >
g\transp(v'-v)$, so $g \notin \partial f(v')$ and $\{\partial f(v)\}_{v\in V(b)}$
tiles $\mathcal{Z}(f)$~\cite{richtergebertziegler1994}, \cite[\S7.5]{ziegler1995},
yielding~\eqref{eq:zonotope-conservation} (cf.~\cite{brandenburggrillohertrich2026}).
\end{proof}
\section{Spectrahedral Realisation Moduli and Polyhedral Rigidity}\label{sec:moduli}

\subsection{Chamber-metric equivalence and the Gram spectrahedron}\label{subsec:gram-spectrahedron}

With the arrangement $\arr = \{\wall{1},\dots,\wall{\ngen}\}$ and hence $(W,b)$
fixed, we classify all parameters $(L, A)$ for which~\eqref{eq:arr-map} induces
prescribed chamber metrics $\{g_\sigma\}_{\sigma\in\Sigma(\arr)}$ modulo
$O(n)\ltimes\RR^n$. Theorem~\ref{thm:gram-spectrahedron} identifies these
\emph{marked} realisations on $(\RR^{\dimin},\arr)$ with a rank-constrained
spectrahedron $\mathcal{M}(G_0,n)$. Isometric piecewise-linear embeddings
may re-triangulate the domain~\cite{buragozalgaller1996, demaineorourke2007};
$\mathcal{M}(G_0,n)$ instead keeps the crease walls fixed. Piecewise-linear
self-isometries that permute the walls act on it as discrete symmetries, for
example $t \mapsto 1/(3t)$ in Proposition~\ref{prop:vertex-flex-rigid}. Under
codimension-$2$ genericity, Corollary~\ref{cor:spline-exhaustion} shows
that~\eqref{eq:arr-map} exhausts all continuous piecewise-affine maps
subordinate to $\arr$.

Fix pairwise distinct hyperplanes $\arr = \{\wall{1},\dots,\wall{\ngen}\}$ in
$\RR^{\dimin}$ with normals $w_1,\dots,w_{\ngen} \in \RR^{\dimin}\setminus\{0\}$,
and write $\Sigma(\arr) \subset \{0,1\}^{\ngen}$ for the sign patterns of its
non-empty open chambers $\cell{\sigma}$. Collect $(L, A)$ into
\[
  B := [L\ \ A] \in \RR^{n\times(\dimin+\ngen)},\qquad
  G := B\transp B = \begin{pmatrix} L\transp L & L\transp A \\ A\transp L & A\transp A \end{pmatrix}
  \in \Spec_+^{\dimin+\ngen}.
\]
On chamber $\cell{\sigma}$ ($\sigma \in \Sigma(\arr)$),
\[
  J_\sigma = B \Upsilon_\sigma,\qquad
  g_\sigma = \Upsilon_\sigma\transp G \Upsilon_\sigma,\qquad
  \Upsilon_\sigma := \begin{pmatrix} I_{\dimin} \\ \diag(\sigma)\,W \end{pmatrix}
  \in \RR^{(\dimin+\ngen)\times\dimin}.
\]
For $h \in \{1,\dots,\ngen\}$, with standard basis vector $e_h \in \RR^{\ngen}$,
define the rank-two \emph{wall matrix}
\begin{equation}\label{eq:wall-matrix}
  E_h
  := \begin{pmatrix} 0 & w_h e_h\transp \\ e_h w_h\transp & -2\,e_h e_h\transp \end{pmatrix}
  \in \Spec^{\dimin+\ngen}.
\end{equation}
We call $\arr$ \emph{codimension-$2$ generic} if no two normals $w_h, w_{h'}$
with $h \ne h'$ are parallel and no intersection $\wall{h}\cap\wall{h'}$ lies in
a third wall $\wall{h''}$, $h'' \notin \{h,h'\}$. When $\dimin = 1$ we require
$\ngen = 1$.
Every simple arrangement with pairwise non-parallel normals is codimension-$2$
generic.

\begin{theorem}[Spectrahedral realisation slice]
\label{thm:gram-spectrahedron}
Let $W = [w_1\ \cdots\ w_{\ngen}]\transp \in \RR^{\ngen\times\dimin}$ have
$w_h \ne 0$ for all $h$.
\begin{enumerate}
  \item For every $h \in \{1,\dots,\ngen\}$ and $\sigma \in \{0,1\}^{\ngen}$,
    $\sigma_h^2 = \sigma_h$ gives
    \begin{equation}\label{eq:boolean-annihilation}
      \Upsilon_\sigma\transp E_h \Upsilon_\sigma
      = 2(\sigma_h - \sigma_h^2)\,w_h w_h\transp
      = 0,
    \end{equation}
    and $E_1,\dots,E_{\ngen}$ are linearly independent in $\Spec^{\dimin+\ngen}$.
  \item $\mathcal{N}_{\mathrm{all}}(W) := \{X \in \Spec^{\dimin+\ngen} :
    \Upsilon_\sigma\transp X \Upsilon_\sigma = 0\ \forall\sigma \in
    \{0,1\}^{\ngen}\} = \spn\{E_1,\dots,E_{\ngen}\}$ has dimension $\ngen$.
  \item Define the \emph{realisable chamber nullspace}
    $\mathcal{N}(\arr) := \{X \in \Spec^{\dimin+\ngen} : \Upsilon_\sigma\transp X
    \Upsilon_\sigma = 0\ \forall\sigma \in \Sigma(\arr)\}$. Each wall $\wall{h}$
    bounds an adjacent chamber pair in $\Sigma(\arr)$ differing only in~$h$, and
    every pair $\{\wall{h},\wall{h'}\}$ ($h \ne h'$) meets in a codimension-$2$
    flat surrounded by a $4$-cycle in $\Sigma(\arr)$ flipping $(h,h')$ alone if
    and only if $\arr$ is codimension-$2$ generic, in which case
    $\mathcal{N}(\arr) = \spn\{E_1,\dots,E_{\ngen}\}$ of dimension $\ngen$.
  \item Assume $\arr$ is codimension-$2$ generic. Given $B_0 = [L_0\ A_0] \in
    \RR^{n\times(\dimin+\ngen)}$ with $G_0 = B_0\transp B_0$, a matrix $G \in
    \Spec_+^{\dimin+\ngen}$ satisfies $\Upsilon_\sigma\transp G \Upsilon_\sigma =
    \Upsilon_\sigma\transp G_0 \Upsilon_\sigma$ for all $\sigma \in \Sigma(\arr)$
    if and only if $G = G(\mu) := G_0 + \sum_{h=1}^{\ngen} \mu_h E_h$ with $\mu$
    in $\mathcal{S}(G_0) := \{\mu \in \RR^{\ngen} : G(\mu) \succeq 0\}$, and
    $\widetilde{B} \mapsto \widetilde{B}\transp\widetilde{B}$ gives a bijection
    from $\{\widetilde{B} \in \RR^{n\times(\dimin+\ngen)} :
    \Upsilon_\sigma\transp\widetilde{B}\transp\widetilde{B}\Upsilon_\sigma =
    g_\sigma\ \forall\sigma \in \Sigma(\arr)\}/O(n)$ onto
    \[
      \mathcal{M}(G_0, n)
      := \bigl\{\mu \in \RR^{\ngen} : G(\mu) \succeq 0\text{ and }\rank G(\mu) \le n\bigr\}.
    \]
\end{enumerate}
\end{theorem}

\begin{proof}
\textit{(1)--(3).}\
Since the walls are pairwise distinct, each $\wall{h}\cap\wall{h'}$ ($h' \ne h$)
is empty or proper in $\wall{h}$, so a small ball at a generic point of
$\wall{h}$ meets an adjacent chamber pair flipping $\sigma_h$ alone. If $w_h
\parallel w_{h'}$, at most three sign patterns on $(h,h')$ occur. If instead
$\Phi_{\{h,h'\}} := \wall{h}\cap\wall{h'} \subseteq \wall{h''}$ ($h'' \notin
\{h,h'\}$), then fixing $\sigma_{h''} \in \{0,1\}$ leaves only three of the six
transverse sectors, so no $4$-cycle exists. Conversely, under codimension-$2$
genericity, $\dim\Phi_{\{h,h'\}} = \dimin - 2$ and every $\wall{h''}$ ($h''
\notin \{h,h'\}$) cuts $\Phi_{\{h,h'\}}$ properly, so a small ball around a
generic point of $\Phi_{\{h,h'\}}$ meets a $4$-cycle flipping $(h,h')$.

Writing $X = \bigl(\begin{smallmatrix} X_{LL} & X_{LA} \\ X_{LA}\transp & X_{AA}
\end{smallmatrix}\bigr) \in \Spec^{\dimin+\ngen}$ with $X_{LA} =
[u_1\ \cdots\ u_{\ngen}]$ and $X_{AA} = (x_{hh'})$ gives
\begin{equation}\label{eq:nullspace-quadratic}
  \Upsilon_\sigma\transp X \Upsilon_\sigma
  = X_{LL}
    + \sum_{h=1}^{\ngen} \sigma_h\bigl(u_h w_h\transp + w_h u_h\transp + x_{hh} w_h w_h\transp\bigr)
    + \sum_{1\le h < h'\le\ngen} \sigma_h\sigma_{h'}\,x_{hh'}
      \bigl(w_h w_{h'}\transp + w_{h'} w_h\transp\bigr).
\end{equation}
Taking $X = E_h$ gives~\eqref{eq:boolean-annihilation}, linear independence of
$\{E_h\}$ (from their $-2\,e_he_h\transp$ diagonal entries), and
$\spn\{E_1,\dots,E_{\ngen}\} \subseteq \mathcal{N}(\arr)$. Conversely, on
$\{0,1\}^{\ngen}$ or under codimension-$2$ genericity, the second difference
of~\eqref{eq:nullspace-quadratic} across the $(h,h')$ $4$-cycle gives
$x_{hh'}(w_h w_{h'}\transp + w_{h'} w_h\transp) = 0$, hence $x_{hh'} = 0$ for
$h \ne h'$. The first difference across an adjacent pair flipping $\sigma_h$
then gives $u_h = -\tfrac{1}{2}x_{hh} w_h =: \mu_h w_h$. Evaluating on any
$\sigma \in \Sigma(\arr)$ gives $X_{LL} = 0$, so $X = \sum_{h=1}^{\ngen} \mu_h
E_h$.

\ppart{(4).}
A matrix $G \in \Spec_+^{\dimin+\ngen}$ factors as
$\widetilde{B}\transp\widetilde{B}$ with $\widetilde{B} \in
\RR^{n\times(\dimin+\ngen)}$, unique modulo $O(n)$, if and only if $G \succeq 0$
and $\rank G \le n$. Hence~(4) follows from~(3).
\end{proof}

\begin{corollary}[Spline exhaustion, after~\cite{billera1988, billerarose1991}]
\label{cor:spline-exhaustion}
Call a continuous map $F : \RR^{\dimin} \to \RR^n$ \emph{piecewise-affine
subordinate to $\arr$} if it is affine on each closed chamber $\overline{\cell{\sigma}}$.
When $\arr$ is codimension-$2$ generic, every such map has the
form~\eqref{eq:arr-map} for a unique $(L, c_0, A)$, so $\mathcal{M}(G_0,n)$
classifies modulo $O(n)\ltimes\RR^n$ all continuous piecewise-affine maps
subordinate to $\arr$ with chamber metrics $\{g_\sigma\}_{\sigma\in\Sigma(\arr)}$.
When codimension-$2$ genericity fails (e.g.\ central arrangements with $\ngen >
\dimin$),~\eqref{eq:arr-map} is in general a proper subfamily.
\end{corollary}

\begin{proof}
For $\dimin = 1$, codimension-$2$ genericity gives $\ngen = 1$, and the claim
is immediate. In fact, for $\dimin = 1$ exhaustion holds for every $\ngen \ge
1$, since a continuous piecewise-affine function of one variable is determined
by its slope jumps. Assume now that $\dimin \ge 2$, and write $D_\sigma$ for the
chamber differentials of $F$. Across any facet $f \subset \wall{h}$ separating
$\sigma,\sigma'$ ($\sigma'_h = 1, \sigma_h = 0$), continuity along $w_h^\perp$
gives $D_{\sigma'} - D_\sigma = a_h^{(f)}w_h\transp$. Across any $(\dimin-2)$-face
$e \subset \wall{h}\cap\wall{h'}$ between adjacent facets $f, f' \subset \wall{h}$,
the $4$-cycle around $e$ (Theorem~\ref{thm:gram-spectrahedron}(3)) gives
$(a_h^{(f)} - a_h^{(f')})w_h\transp = (a_{h'}^{(g')} - a_{h'}^{(g)})w_{h'}\transp$,
forcing $a_h^{(f)} = a_h^{(f')}$ since $w_h \nparallel w_{h'}$. Connectedness of
the chamber graph on $\wall{h} \cong \RR^{\dimin-1}$ and on
$\RR^{\dimin}$~\cite{orlikterao1992, ziegler1995} makes $a_h^{(f)} =: a_h$
constant along $\wall{h}$. Starting from a base chamber $\sigma_0$, it yields
$D_\sigma = L + \sum_h\sigma_h a_hw_h\transp$, so $F(x) = Lx + c_0 + A(Wx+b)_+$
uniquely.

Conversely, three central lines in $\RR^2$ with normals $e_1, e_2, e_1+e_2$ admit
six ray jumps $\alpha_r$ subject to $\sum_r\alpha_rw_{h(r)} = 0$, giving $6 - 2 =
4 > 3$ parameters per coordinate modulo $L$.
\end{proof}

\begin{remk}\label{rem:q2-counterexamples}
Without codimension-$2$ genericity, $\mathcal{N}(\arr)$ can strictly contain
$\spn\{E_1,\dots,E_{\ngen}\}$ (where $x_{hh'} = 0$ for $h \ne h'$). First take
a parallel pair: $\dimin=2$, $\ngen=3$, $w_1=w_2=e_1$, $w_3=e_2$,
$b=(0,-1,0)\transp$. The pattern $(\sigma_1,\sigma_2) = (0,1)$ is absent, so
$\sigma_1\sigma_2 = \sigma_2$ on $\Sigma(\arr)$. Setting $u_2 = -e_1$, $x_{12}
= 1$ and all other entries $0$ gives $\Upsilon_\sigma\transp X\Upsilon_\sigma =
0$ with $X \notin \spn\{E_1,E_2,E_3\}$; here $\dim\mathcal{N}(\arr) = 4$.

Four central lines in $\RR^2$ behave the same way. Take $\ngen = 4$, $b = 0$,
$w_1 = e_1$, $w_2 = e_2$, $w_3 = e_1+e_2$, $w_4 = e_1-e_2$. The $8$ chambers
satisfy $(\sigma_3-\sigma_1)(\sigma_3-\sigma_2) = 0$ and
$(\sigma_4-\sigma_1)(\sigma_4+\sigma_2-1) = 0$, so the matrix $X$ with $X_{LL} =
0$, $x_{hh} = 0$, $u_1 = (1,-1)\transp$, $u_3 = (0,1)\transp$, $u_2 = u_4 = 0$,
and $(x_{12},x_{13},x_{14},x_{23},x_{24},x_{34}) = (2,-1,-1,-1,-1,1)$ lies in
$\mathcal{N}(\arr)\setminus\spn\{E_1,\dots,E_4\}$. In this case
$\dim\mathcal{N}(\arr) = 5$.

By contrast, any $\ngen = 3$ pairwise non-parallel central lines in $\RR^2$ still
satisfy $\mathcal{N}(\arr) = \spn\{E_1,E_2,E_3\}$. To see this, subtract the
first differences of~\eqref{eq:nullspace-quadratic} across the two opposite
chamber pairs flipping $\sigma_h$ alone. This gives $w_h v\transp + v w_h\transp
= 0$ with $v := \pm x_{hh'}w_{h'} \pm x_{hh''}w_{h''}$, which forces $v \in \RR
w_h$ and hence $v = 0$. Independence of $w_{h'}, w_{h''}$ gives $x_{hh'} =
x_{hh''} = 0$, and so $u_h = -\tfrac12 x_{hh}w_h$ and $X_{LL} = 0$.
\end{remk}

\begin{remk}\label{rem:gram-PR}
The spectrahedron $\mathcal{S}(G_0)$ and its rank strata $\mathcal{M}(G_0,n)$ are
a piecewise-linear counterpart of sum-of-squares Gram
spectrahedra~\cite{chuaplaumannsinnvinzant2017, blekherman2013}. We say that
$(B, W)$, where $B = [L\ A]$, has
\emph{full-cube positive path rank} if
\begin{equation}\label{eq:PR-global}
  \rank\bigl(L + A\diag(s)W\bigr) = \dimin
  \qquad\text{for all } s \in [0,1]^{\ngen}.
\end{equation}
This implies positive path rank at every face for every $b$. Suppose $n =
\dimin$ and $\det\widetilde{L} \ne 0$. Restricting~\eqref{eq:PR-global} to the
edges of $[0,1]^{\ngen}$ at $0$ then gives \emph{single-wall orientation
preservation}:
$\det(\widetilde{L} + s_h\widetilde{a}_h w_h\transp)/\det\widetilde{L} = 1 + s_h
w_h\transp\widetilde{L}^{-1}\widetilde{a}_h > 0$ for all $h$ and $s_h \in [0,1]$.
In this square case the whole cube can be treated in the same way. Put
$M_{\mathrm{P}} := I_{\ngen} + W\widetilde{L}^{-1}\widetilde{A}$ and write
$(M_{\mathrm{P}})_{SS}$ for its principal submatrix on $S \subseteq
\{1,\dots,\ngen\}$, with $\det(M_{\mathrm{P}})_{\emptyset\emptyset} := 1$.
Sylvester's determinant identity and multilinearity in the rows give
\begin{equation}\label{eq:PR-pmatrix}
\begin{aligned}
  \frac{\det\bigl(\widetilde{L} + \widetilde{A}\diag(s)W\bigr)}{\det\widetilde{L}}
  &= \det\bigl(I_{\ngen} - \diag(s) + \diag(s)\,M_{\mathrm{P}}\bigr)\\
  &= \sum_{S\subseteq\{1,\dots,\ngen\}}\,
     \prod_{h\in S}s_h\prod_{h\notin S}(1-s_h)\,\det(M_{\mathrm{P}})_{SS}.
\end{aligned}
\end{equation}
For $s \in [0,1]^{\ngen}$ the right-hand side is a convex combination of the
principal minors of $M_{\mathrm{P}}$. It equals $1$ at $s = 0$ and
$\det(M_{\mathrm{P}})_{SS}$ at the vertex $s = \mathbf{1}_S$, and $[0,1]^{\ngen}$
is connected. Hence full-cube positive path rank holds if and only if
$M_{\mathrm{P}}$ is a \emph{P-matrix}, that is, all its principal minors are
positive. Deciding whether a rational matrix is a P-matrix is
co-NP-complete~\cite{coxson1994}, and $W = \widetilde{L} = I_{\dimin}$ makes
$M_{\mathrm{P}} = I_{\dimin} + \widetilde{A}$ arbitrary. Deciding~\eqref{eq:PR-global}
is therefore co-NP-hard when $\dimin$ and $\ngen$ are part of the input.
With $\Upsilon(s) := [I_{\dimin};\,\diag(s)W]$,
\[
  \mathcal{M}_{\mathrm{PR}}(G_0,n)
  := \bigl\{\mu \in \mathcal{M}(G_0,n) :
     \Upsilon(s)\transp G(\mu)\Upsilon(s) \succ 0\ \ \forall s \in [0,1]^{\ngen}\bigr\}
\]
is an open subset of $\mathcal{M}(G_0,n)$ that excludes orientation-reversing
fold-overs (Theorem~\ref{thm:compactness-stress}(2),
Remark~\ref{rem:compactness-examples}).

Both conditions are uniform in the offset $b$ and refer to the formal cube
$[0,1]^{\ngen}$ rather than to the faces that $\arr(b)$ actually has. For a
face $\mathsf{F} \in \mathcal{F}(\arr(b))$ and $s \in [0,1]^{I(\mathsf{F})}$,
set $\Upsilon_{\mathsf{F}}(s) := \Upsilon(s^{\mathsf{F}}(s))$, where
$s^{\mathsf{F}}_h := s_h$ for $h \in I(\mathsf{F})$ and $s^{\mathsf{F}}_h :=
\mathbf{1}\{z_h > 0\text{ on }\mathsf{F}\}$ for $h \notin I(\mathsf{F})$. Then
$\widetilde{J}_{\mathsf{F}}(s)\transp\widetilde{J}_{\mathsf{F}}(s) =
\Upsilon_{\mathsf{F}}(s)\transp G(\mu)\Upsilon_{\mathsf{F}}(s)$ does not depend
on the $O(n)$-choice of factor $\widetilde{B}$ with $G(\mu) =
\widetilde{B}\transp\widetilde{B}$. The offset-specific \emph{face-immersive
locus}
\[
  \mathcal{M}_{\mathrm{face}}(G_0,n;b)
  := \bigl\{\mu \in \mathcal{M}(G_0,n) :
     \Upsilon_{\mathsf{F}}(s)\transp G(\mu)\Upsilon_{\mathsf{F}}(s) \succ 0
     \ \ \forall\mathsf{F} \in \mathcal{F}(\arr(b)),\
     \forall s \in [0,1]^{I(\mathsf{F})}\bigr\}
\]
encodes positive path rank at every face of $\arr(b)$. It contains
$\mathcal{M}_{\mathrm{PR}}(G_0,n)$, and it can be strictly larger when $\arr(b)$
has fewer than $2^{\ngen}$ chambers
(Proposition~\ref{prop:vertex-flex-rigid}(4)). When $\arr(b)$ is simple, every
realisation in this locus is a globally generic arrangement immersion
(Section~\ref{subsec:arr-local-models}).
\end{remk}

\subsection{Local realisation dimension, flexing, and rigidification}
\label{subsec:local-moduli}

Near $\mu = 0$, the Schur-complement chart~\cite{blekherman2013} on the
rank-$\le n_0$ determinantal locus linearises the rank constraint at $G_0$ into
the congruence map $X \mapsto Z\transp X Z$.

\begin{lemma}\label{lem:psd-schur-chart}
Let $G_0 \in \Spec_+^{\dimin+\ngen}$ have rank $n_0 \le \dimin+\ngen$, put
$r := \dimin+\ngen - n_0$, and choose $Q = [Q_+\ Z] \in O(\dimin+\ngen)$ with
$\im Q_+ = \im G_0$ and $\im Z = \ker G_0$. On $\mathcal{U} := \{G \in
\Spec^{\dimin+\ngen} : G_{++}(G) := Q_+\transp G Q_+ \succ 0\}$, set
$\Sigma_{\ker}(G) := Z\transp G Z - (Q_+\transp G Z)\transp G_{++}(G)^{-1}
(Q_+\transp G Z) \in \Spec^r$.

Then $G \in \mathcal{U}$ satisfies $G \succeq 0$
and $\rank G \le n_0$ if and only if $\rank G = n_0$ and $\Sigma_{\ker}(G) = 0$,
and
\begin{equation}\label{eq:schur-chart-diff}
  D\Sigma_{\ker}(G_0)[X] = Z\transp X Z \in \Spec^r,
  \qquad X \in \Spec^{\dimin+\ngen}.
\end{equation}
\end{lemma}

\begin{proof}
By block Gaussian elimination, $Q\transp G Q$ is block congruent to
$\diag(G_{++}(G), \Sigma_{\ker}(G))$. By Haynsworth's inertia additivity
formula~\cite{haynsworth1968}, this gives $\rank G = n_0 +
\rank\Sigma_{\ker}(G)$ and $G \succeq 0 \iff \Sigma_{\ker}(G) \succeq 0$ on
$\mathcal{U}$. The identity $Q_+\transp G_0 Z = 0$
gives~\eqref{eq:schur-chart-diff}.
\end{proof}

\begin{lemma}\label{lem:wall-congruence}
Partition $Z = [Z_L;\, Z_A] \in \RR^{(\dimin+\ngen)\times r}$ into $Z_L \in
\RR^{\dimin\times r}$ and $Z_A = [\eta_1\ \cdots\ \eta_{\ngen}]\transp \in
\RR^{\ngen\times r}$, and set $\zeta_h := Z_L\transp w_h - \eta_h \in \RR^r$. Then
\begin{equation}\label{eq:wall-congruence}
  Z\transp E_h Z
  = \zeta_h \eta_h\transp + \eta_h \zeta_h\transp \in \Spec^r
  \qquad (1 \le h \le \ngen),
\end{equation}
and $\dim\spn\{Z\transp E_1 Z,\dots,Z\transp E_{\ngen} Z\}$ depends only on
$\im Z \subset \RR^{\dimin+\ngen}$.
\end{lemma}

\begin{proof}
Expanding~\eqref{eq:wall-matrix} gives $Z\transp E_h Z = Z_L\transp w_h\eta_h\transp
+ \eta_h w_h\transp Z_L - 2\eta_h\eta_h\transp = \zeta_h\eta_h\transp +
\eta_h\zeta_h\transp$, and $Z \mapsto ZR$ ($R \in \operatorname{GL}_r(\RR)$) acts
on $\Spec^r$ by the automorphism $S \mapsto R\transp S R$.
\end{proof}

Fix $\dimin \ge 2$, $\ngen \ge 1$, $n \ge \dimin$, and $r := \max(0,\,
\dimin+\ngen-n)$. For $B = [L\ A] \in \RR^{n\times(\dimin+\ngen)}$ of full rank
$\min(n,\dimin+\ngen)$ and $W \in \RR^{\ngen\times\dimin}$ with kernel basis
$Z(B) \in \RR^{(\dimin+\ngen)\times r}$, define the linear obstruction map
\[
  \Lambda_{(B,W)} : \RR^{\ngen} \to \Spec^r,\qquad
  \Lambda_{(B,W)}(\mu) := \sum_{h=1}^{\ngen} \mu_h\,Z(B)\transp E_h(W)\,Z(B)
\]
with $\Spec^0 := \{0\}$ and $\Lambda_{(B,W)} := 0$ for $r = 0$; we abbreviate
$\Lambda := \Lambda_{(B_0,W)}$. Here $\binom{r+1}{2} = \dim\Spec^r$ is the
codimension of the rank-$\min(n,\dimin+\ngen)$ stratum in
$\Spec_+^{\dimin+\ngen}$~\cite{barvinok1995, blekherman2013, pataki1998}.
Lemma~\ref{lem:obstruction-genericity} in Appendix~\ref{app:moduli-proofs}
constructs a non-empty Zariski-open locus $\Gamma$ on which the sparse rank-two
wall pencil $G_0 + \spn\{E_h(W)\}$ meets that stratum transversely, that is,
$\rank\Lambda_{(B,W)} = \min\bigl(\ngen,\binom{r+1}{2}\bigr)$. The appendix
then builds the local real-analytic chart (Lemma~\ref{lem:local-moduli-manifold})
and proves joint Lebesgue genericity (Lemma~\ref{lem:joint-genericity}).

\begin{theorem}[Generic local dimension of fixed-$n$ realisation moduli]
\label{thm:generic-local-dim}
Fix $\dimin \ge 2$, $\ngen \ge 1$, $n \ge \dimin$, and $r := \max(0,\,
\dimin + \ngen - n)$. For every $(L_0, A_0, W) \in \Gamma$
(Lemma~\ref{lem:obstruction-genericity}), $B_0 = [L_0\ A_0]$ has rank
$\min(n, \dimin+\ngen)$ and the germ at $\mu = 0$ of $\mathcal{M}(G_0, n)$
($G_0 = B_0\transp B_0$) is a real-analytic embedded manifold germ of rank-$\min(n,
\dimin+\ngen)$ matrices of dimension
\begin{equation}\label{eq:generic-local-dim}
  \dim_0 \mathcal{M}(G_0, n)
  = \delta(\dimin,\ngen,n)
  := \max\!\left(0,\, \ngen - \binom{r+1}{2}\right).
\end{equation}
Suppose further that $\arr$ is codimension-$2$ generic. By
Lemma~\ref{lem:joint-genericity}(3), this holds jointly with $(L_0, A_0, W) \in
\Gamma$ and simplicity of $\arr$ for Lebesgue-almost every $(L_0, A_0, W, b)$.
Set
$\mathcal{R}_n := \{\widetilde{B} \in \RR^{n\times(\dimin+\ngen)} :
\Upsilon_\sigma\transp\widetilde{B}\transp\widetilde{B}\Upsilon_\sigma = g_\sigma\
\forall\sigma\in\Sigma(\arr)\}$. The map $\widetilde{B} \mapsto
\widetilde{B}\transp\widetilde{B}$ induces a homeomorphism
$\mathcal{R}_n/O(n) \to \mathcal{M}(G_0,n)$, which near $[B_0]$ identifies
$\mathcal{R}_n/O(n)$ with this manifold germ; by
Corollary~\ref{cor:spline-exhaustion}, $\mathcal{M}(G_0,n)$ therefore
parametrises all isometric continuous piecewise-affine maps subordinate to
$\arr$ modulo $O(n)\ltimes\RR^n$. In particular, $[B_0]$ is locally rigid
modulo $O(n)$ if and only if $\ngen \le \binom{r+1}{2}$. For any $(L_0,A_0,W)$ with $\rank B_0 =
\min(n,\dimin+\ngen)$, $\dim\ker\Lambda = \ngen - \rank\Lambda \ge
\delta(\dimin,\ngen,n)$ bounds the local dimension at $\mu = 0$ from above.
Finally, the locus of $(L_0, A_0, W) \in \Gamma$ with full-cube positive path
rank is non-empty and open.
\end{theorem}

\begin{proof}
On $\Gamma$, Lemma~\ref{lem:obstruction-genericity} gives $\rank B_0 =
\min(n,\dimin+\ngen)$ and $\rank\Lambda = \min(\ngen, \binom{r+1}{2})$, so
Lemma~\ref{lem:local-moduli-manifold}(1)--(3) gives~\eqref{eq:generic-local-dim},
and Lemma~\ref{lem:local-moduli-manifold}(4) gives the general upper bound
$\ngen - \rank\Lambda$. Under codimension-$2$ genericity,
Theorem~\ref{thm:gram-spectrahedron}(4) makes $[\widetilde{B}] \mapsto
\widetilde{B}\transp\widetilde{B}$ a continuous bijection $\mathcal{R}_n/O(n) \to
\mathcal{M}(G_0,n)$; it is proper since $\norm{\widetilde{B}}_F^2 = \tr G(\mu)$
is bounded on compact subsets of $\mathcal{M}(G_0,n)$, hence a homeomorphism.
Full-cube positive path rank is open by compactness of $[0,1]^{\ngen}$ and holds
at $A_0 = 0$, $\rank L_0 = \dimin$. It therefore meets the dense open set
$\Gamma$ in a non-empty open set.
\end{proof}

Unlike polyhedral surfaces hinged independently along
edges~\cite{demaineorourke2007},~\eqref{eq:arr-map} couples all chambers
across $\wall{h}$ by a single rank-$1$ jump $a_h w_h\transp$. This imposes
$\binom{r+1}{2}$ corank equations on $\mu \in \RR^{\ngen}$. At $(\dimin, \ngen,
n) = (2, 2, 3)$, where $r = 1$ and $\delta = 1$, we recover the flex of the
spherical four-bar linkage. A third generic line ($\ngen = 3$, $r = 2$, $\delta
= 0$) rigidifies the surface in $\RR^3$.

To see both the content and the limits of Theorem~\ref{thm:generic-local-dim},
it is instructive to work through this smallest case in detail
(Proposition~\ref{prop:vertex-flex-rigid}): the generic dimension $\delta = 1$
is realised by an explicit flex, and adjoining a third wall gives local
rigidity but, as we shall see, not global uniqueness. Explicitly, fix $L =
[e_1\ e_2]$, so
that $g_{00} = I_2$, and rotate the crease rays $c_1(\phi) = J_{10}e_1$
(length~$2$) and $c_2(\psi) = J_{01}e_2$ (length~$1$) at cosine $-\frac{3}{5}$
with $e_2$ and $e_1$. Orthogonality $c_1\transp c_2 = 0$ ($g_{11} =
\diag(4,1)$) reads $4\sin\phi\sin\psi = 3(\cos\phi+\cos\psi)$. In half-angle
tangents $t = \tan(\phi/2) > 0$, $u = \tan(\psi/2)$, this factors as $(3tu -
1)(tu + 3) = 0$, the coupling of a flexible Kokotsakis
vertex~\cite{kokotsakis1933, izmestiev2017}. The curve splits into the
admissible arc $u = 1/(3t) > 0$ and the fold-over arc $u = -3/t < 0$ (verified
in Appendix~\ref{app:moduli-proofs}; Figure~\ref{fig:flex-rigid}).

\begin{figure}[t]
\centering
\includegraphics[width=\textwidth]{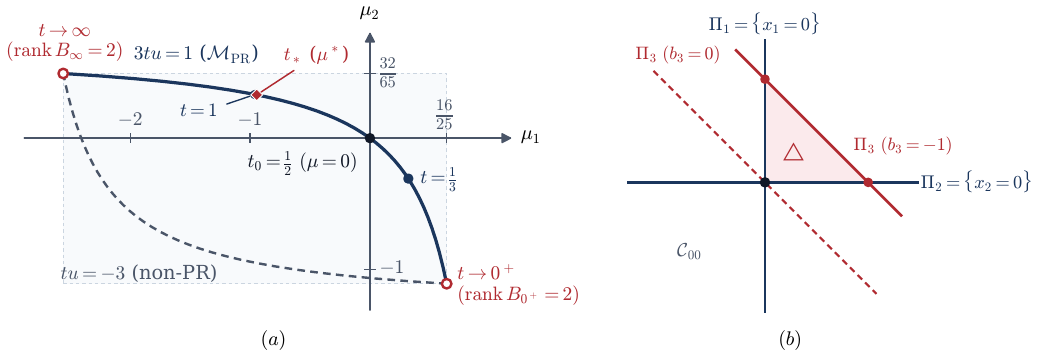}
\caption{Degree-four flex and three-wall rigidification (Proposition~\ref{prop:vertex-flex-rigid}). \emph{(a)}~At $(2,2,3)$, $\mathcal{M}(G_{t_0},3) = \partial\mathcal{S}(G_{t_0})$ splits at the flat-folded endpoints ($t\to 0^+,\infty$) into the solid admissible arc $\mathcal{M}_{\mathrm{PR}}(G_{t_0},3)$ ($3tu=1$, with $t_*$ marking the $(\mu_1,\mu_2)$-projection of $\mu^*$) and the dashed fold-over arc ($tu=-3$). \emph{(b)}~Adjoining $\wall{3}$ ($b_3=-1$ or $b_3=0$) isolates $\mu=0$ locally modulo $O(3)$.}
\label{fig:flex-rigid}
\end{figure}

\begin{proposition}\label{prop:vertex-flex-rigid}
Let $\dimin = 2$, $\ngen = 2$, $n = 3$, $W = I_2$, $b = 0$, and $L = [e_1\ e_2]
\in \RR^{3\times 2}$. For $t \in (0, \infty)$, define $a_1(t) := c_1(t) - e_1$
and $a_2(t) := c_2(t) - e_2$ via
\begin{gather*}
  \cos\phi(t) := \frac{1-t^2}{1+t^2},\quad
  \sin\phi(t) := \frac{2t}{1+t^2},\qquad
  \cos\psi(t) := \frac{9t^2-1}{9t^2+1},\quad
  \sin\psi(t) := \frac{6t}{9t^2+1},\\
  c_1(t) := \bigl(\tfrac{8}{5}\cos\phi(t),\ -\tfrac{6}{5},\ \tfrac{8}{5}\sin\phi(t)\bigr)\transp,
  \qquad
  c_2(t) := \bigl(-\tfrac{3}{5},\ \tfrac{4}{5}\cos\psi(t),\ \tfrac{4}{5}\sin\psi(t)\bigr)\transp.
\end{gather*}
Then:
\begin{enumerate}
  \item For all $t > 0$, $4\sin\phi(t)\sin\psi(t) = 3(\cos\phi(t)+\cos\psi(t))$,
    so $c_1(t)\transp c_2(t) \equiv 0$, $F_t(x) := Lx + A(t)(x)_+$ has constant
    chamber metrics $g_{00} = I_2$,
    $g_{10} = \bigl(\begin{smallmatrix} 4 & -6/5 \\ -6/5 & 1 \end{smallmatrix}\bigr)$,
    $g_{01} = \bigl(\begin{smallmatrix} 1 & -3/5 \\ -3/5 & 1 \end{smallmatrix}\bigr)$,
    $g_{11} = \diag(4,1)$, and relative to
    $t_0 := \frac{1}{2}$, $G_t = B_t\transp B_t = G_{t_0} + \mu_1(t)E_1 +
    \mu_2(t)E_2$ with $\mu_1(t) = \frac{8}{5}(\cos\phi(t) - \frac{3}{5})$ and
    $\mu_2(t) = \frac{4}{5}(\cos\psi(t) - \frac{5}{13})$.
  \item For all $t > 0$, $\rank B_t = 3$, global genericity, codimension-$2$
    genericity, and $(L, A(t), W) \in \Gamma$ hold, and $t \mapsto \mu(t)$ is a
    real-analytic homeomorphism of $(0,\infty)$ onto $\mathcal{M}_{\mathrm{PR}}(G_{t_0},
    3)$ (with $B_t \not\sim_{O(3)} B_{t'}$ for $t \ne t'$, although $F_t(\RR^2)
    \cong_{O(3)} F_{1/(3t)}(\RR^2)$ after a piecewise-linear reparametrisation
    of $\RR^2$ that swaps the two walls). The
    full space $\mathcal{M}(G_{t_0}, 3) = \partial\mathcal{S}(G_{t_0}) \cong S^1$
    is a nodal quartic oval comprising $\mathcal{M}_{\mathrm{PR}}(G_{t_0}, 3)$,
    two rank-$2$ flat-folded endpoints $(\frac{16}{25},-\frac{72}{65})$ ($t \to
    0^+$) and $(-\frac{64}{25},\frac{32}{65})$ ($t \to \infty$), and the fold-over
    arc $tu = -3$ ($\rank B = 3$, $\rank J(\frac{1}{3},\frac{1}{2}) = 1$).
  \item Along the flex, $K(0) = \frac{1}{\pi}\arcsin(3/5)$ is constant whereas
    $h_{2,1/3}(\tfrac{1}{2},\tfrac{1}{2}) > h_{2,1}(\tfrac{1}{2},\tfrac{1}{2})$.
  \item Adjoining $w_3 = (1,1)\transp$, $a_3 = e_3$, and $b_3 \in \{-1, 0\}$ at
    $t_0 = \frac{1}{2}$ gives $\rank B = 3$, full-cube positive path rank, and
    bijective $\Lambda : \RR^3 \to \Spec^2$. At $b_3 = -1$, global genericity,
    codimension-$2$ genericity, and $(L,A,W) \in \Gamma$ hold with
    $\dim_0\mathcal{M}(G_0,3) = 0$ and $\mathcal{M}_{\mathrm{PR}}(G_0,3) = \{0\}$.
    Nevertheless $\mathcal{M}(G_0,3)$ contains a second point $\mu^* =
    (\mu_1(t_*), \mu_2(t_*), \tfrac{16}{13}\mu_1(t_*)) \approx
    (-0.9449,\,0.3296,\,-1.1629)$, where $t_*^2 \approx 0.9813$ is the unique
    positive root of $294{,}849\,x^3 + 203{,}198\,x^2 - 451{,}015\,x -
    31{,}684$. This point satisfies $\mu^* \in \mathcal{M}_{\mathrm{face}}(G_0,3;b)
    \setminus \mathcal{M}_{\mathrm{PR}}(G_0,3)$: it has positive path rank on
    every non-empty face of $\arr(b)$ (the chamber $\cell{001}$ is empty).
    Hence $\mu = 0$ and $\mu^*$ are two non-congruent globally generic
    arrangement immersions with the same chamber metrics. Along a curve in
    $(0,1)^3$ ending on the unrealised face $s_1 = 0$ at $s \approx
    (0,\,0.9290,\,0.2568)$, however, $\rank J(s)$ drops to $1$. At $b_3 = 0$,
    simplicity and
    codimension-$2$ genericity fail, yet $\mathcal{N}(\arr) =
    \operatorname{span}\{E_1,E_2,E_3\}$ gives local rigidity modulo $O(3)$
    within~\eqref{eq:arr-map}.
\end{enumerate}
\end{proposition}

\subsection{Global compactness, minimal-codimension uniqueness, and stress matrices}
\label{subsec:global-spectrahedron}

Whenever $w_h \ne 0$ for all $h$, $\mathcal{S}(G_0)$ is compact. At minimal
target dimension $n = \dimin$, single-wall orientation preservation isolates
$\mu = 0$ globally (Theorem~\ref{thm:compactness-stress}). In higher
codimension, consider a PSD multiplier $\Omega = Z S_\Omega Z\transp$
($S_\Omega \succ 0$) orthogonal to $\{E_h\}$. It is an arrangement counterpart
of a Connelly--Gortler stress matrix~\cite{connelly1982, gortlerthurston2014},
and we write it in the Gale-matrix form of Alfakih~\cite{alfakih2007,
alfakih2018}, with $\rank\Omega = r = \dimin+\ngen-n$. Such an $\Omega$
certifies universal rigidity $\mathcal{S}(G_0) = \{0\}$
(Corollary~\ref{cor:stress-rigidity}).

\begin{theorem}[Spectrahedral compactness and minimal-codimension uniqueness]
\label{thm:compactness-stress}
Fix $n \ge \dimin \ge 1$, $\ngen \ge 1$, and $w_h \ne 0$ ($1 \le h \le \ngen$).
For every $G_0 \in \Spec_+^{\dimin+\ngen}$, $\mathcal{S}(G_0) = \{\mu \in
\RR^{\ngen} : G(\mu) \succeq 0\}$ (hence every $\mathcal{M}(G_0,n)$) has
$\operatorname{rec}(\mathcal{S}(G_0)) = \{0\}$ and is compact. Now let $G_0 =
B_0\transp B_0$ with $B_0 = [L\ A] \in \RR^{n\times(\dimin+\ngen)}$ and
$\rank L = \dimin$, and set $P_L := L(L\transp L)^{-1}L\transp$, $P_L^\perp :=
I_n - P_L$, $q_h := w_h\transp(L\transp L)^{-1}w_h > 0$, and $\ell_h := 1 +
w_h\transp(L\transp L)^{-1}L\transp a_h$. Then:
\begin{enumerate}
  \item $G(\mu) \succeq 0$ if and only if the Schur complement of $L\transp L$ in
    $G(\mu)$,
    \[
      S_L(\mu)
      := A\transp A - 2\diag(\mu)
         - \bigl(A\transp L + \diag(\mu)W\bigr)(L\transp L)^{-1}
           \bigl(L\transp A + W\transp\diag(\mu)\bigr),
    \]
    lies in $\Spec_+^{\ngen}$, in which case
    $\rank G(\mu) = \dimin + \rank S_L(\mu)$. For any realisation $\widetilde{B}
    = [\widetilde{L}\ \widetilde{A}] \in \RR^{n'\times(\dimin+\ngen)}$ ($n' :=
    \max(n,\rank G(\mu))$) of $G(\mu)$ in gauge $\widetilde{L} =
    [L;\,0_{(n'-n)\times\dimin}]$, $P_{L,\RR^{n'}}\widetilde{a}_h = [P_L a_h +
    \mu_h L(L\transp L)^{-1}w_h;\,0]$ and
    \begin{equation}\label{eq:schur-entries}
      \begin{aligned}
        S_L(\mu)_{hh'}
        &= \ip{P_{L,\RR^{n'}}^\perp\widetilde{a}_h}{P_{L,\RR^{n'}}^\perp\widetilde{a}_{h'}}
         = \ip{P_L^\perp a_h}{P_L^\perp a_{h'}} - 2\delta_{hh'}\mu_h \\
        &\quad - \mu_h w_h\transp(L\transp L)^{-1}L\transp a_{h'}
           - \mu_{h'} a_h\transp L(L\transp L)^{-1}w_{h'}
           - \mu_h\mu_{h'} w_h\transp(L\transp L)^{-1}w_{h'}.
      \end{aligned}
    \end{equation}
    Each diagonal entry
    \[
      S_L(\mu)_{hh}
      = \norm{P_{L,\RR^{n'}}^\perp\widetilde{a}_h}_2^2
      = \norm{P_L^\perp a_h}_2^2 - 2\ell_h\mu_h - q_h\mu_h^2
    \]
    is strictly concave in $\mu_h$, giving the bounding box (with
    $0 \in [\mu_h^-,\mu_h^+]$ for each $h$)
    \begin{equation}\label{eq:spectrahedral-box}
      \mathcal{S}(G_0)
      \subset \prod_{h=1}^{\ngen} [\mu_h^-,\, \mu_h^+],
      \qquad
      \mu_h^\pm := \frac{-\ell_h \pm \sqrt{\ell_h^2 + q_h\norm{P_L^\perp a_h}_2^2}}{q_h}.
    \end{equation}
  \item For each $h$, $0 \in \{\mu_h^-, \mu_h^+\} \iff P_L^\perp a_h = 0$, where
    $S_L(\mu)_{hh} = -\mu_h(2\ell_h + q_h\mu_h)$ vanishes on $\{0, -2\ell_h/q_h\}$.
    When $n = \dimin$, $\mathcal{M}(G_0,\dimin) = \{\mu : S_L(\mu) = 0\} \subset
    \prod_{h=1}^{\ngen}\{0, -2\ell_h/q_h\}$. Under single-wall orientation
    preservation at $B_0$ (Remark~\ref{rem:gram-PR}), $\mu = 0$ is the unique
    orientation-preserving point of $\mathcal{M}(G_0,\dimin)$, so
    $\mathcal{M}_{\mathrm{PR}}(G_0,\dimin) \subseteq \{0\}$ (with equality under
    full-cube positive path rank).
\end{enumerate}
\end{theorem}

\begin{proof}
If $\nu \in \operatorname{rec}(\mathcal{S}(G_0))$, then $\sum_{h=1}^{\ngen}\nu_h
E_h \succeq 0$ has zero upper-left $\dimin\times\dimin$ block
by~\eqref{eq:wall-matrix}, forcing $[\nu_1 w_1\ \cdots\ \nu_{\ngen} w_{\ngen}] =
0$ and $\nu = 0$, so $\mathcal{S}(G_0) \ni 0$ is compact.

\smallskip\noindent\textit{(1).}\
With $L\transp L \succ 0$ and $C_L(\mu) := (L\transp L)^{-1}(L\transp A +
W\transp\diag(\mu))$, block congruence by $\bigl(\begin{smallmatrix} I_{\dimin} &
-C_L(\mu) \\ 0 & I_{\ngen} \end{smallmatrix}\bigr)$ reduces $G(\mu)$ to
$\diag(L\transp L, S_L(\mu))$, giving $G(\mu) \succeq 0 \iff S_L(\mu) \succeq 0$
and $\rank G(\mu) = \dimin + \rank S_L(\mu)$. In gauge $\widetilde{L} =
[L;\,0_{(n'-n)\times\dimin}]$, $\widetilde{L}\transp\widetilde{a}_h = L\transp a_h
+ \mu_h w_h$ yields $P_{L,\RR^{n'}}\widetilde{a}_h$
and~\eqref{eq:schur-entries}--\eqref{eq:spectrahedral-box}.

\ppart{(2).}
Since $\mu_h^-\mu_h^+ = -\norm{P_L^\perp a_h}_2^2/q_h$, $0 \in \{\mu_h^-, \mu_h^+\}
\iff P_L^\perp a_h = 0$. At $n = \dimin$, $P_L^\perp = 0$ and $\rank G(\mu) \le
\dimin \iff S_L(\mu) = 0 \implies \mu_h \in \{0, -2\ell_h/q_h\}$. If $B_0$
satisfies single-wall orientation preservation, $\ell_h = 1 + w_h\transp L^{-1}a_h
> 0$ for all $h$, and in gauge $\widetilde{L} = L$ ($\widetilde{a}_h = a_h + \mu_h
L^{-\!\top}w_h$) we have $\det(L + \widetilde{a}_h w_h\transp)/\det L = \ell_h +
q_h\mu_h = -\ell_h < 0$ whenever $\mu_h = -2\ell_h/q_h$ (dropping rank at $s_h^*
= 1/(1+\ell_h) \in (0,1)$), ruling out every non-zero root.
\end{proof}

\begin{corollary}[Stress-matrix universal rigidity]
\label{cor:stress-rigidity}
In Theorem~\ref{thm:compactness-stress}, suppose $\rank B_0 = n < \dimin + \ngen$,
put $r := \dimin + \ngen - n \ge 1$, and let $Z \in \RR^{(\dimin+\ngen)\times r}$
span $\ker G_0$. If there exists $S_\Omega \in \Spec_{++}^r$ such that
$\Omega := Z S_\Omega Z\transp$ satisfies $\ip{\Omega}{E_h} = \tr(S_\Omega
Z\transp E_h Z) = 0$ for all $h$ and $\Lambda$ is injective, then
$\mathcal{S}(G_0) = \{0\}$ (so $\mathcal{M}(G_0,n') = \{0\}$ for all $n' \ge n$)
and $\ngen \le \binom{r+1}{2} - 1$ (hence $r \ge 2$). Thus $[B_0]$ is universally
rigid modulo $O(n')\ltimes\RR^{n'}$ across all $n' \ge n$ within~\eqref{eq:arr-map}
whenever $\mathcal{N}(\arr) = \spn\{E_1,\dots,E_{\ngen}\}$, and among formal
$2^{\ngen}$-chamber realisations for any $w_h \ne 0$. When $\arr$ is
codimension-$2$ generic, this rigidity extends to all continuous
piecewise-affine maps subordinate to $\arr$
(Theorem~\ref{thm:gram-spectrahedron}, Corollary~\ref{cor:spline-exhaustion}).
\end{corollary}

\begin{proof}
Given $\mu \in \mathcal{S}(G_0)$, we have $0 \le \tr(\Omega G(\mu)) = \tr(S_\Omega Z\transp
G_0 Z) + \sum_{h}\mu_h\ip{\Omega}{E_h} = 0$. Since $G(\mu) \succeq 0$ and
$S_\Omega \succ 0$, $\norm{G(\mu)^{1/2} Z S_\Omega^{1/2}}_F^2 = 0$, forcing
$G(\mu) Z = 0$, $\Lambda(\mu) = Z\transp G(\mu) Z = 0$, and $\mu = 0$. Moreover,
$S_\Omega \in (\im\Lambda)^\perp\setminus\{0\}$, so injectivity gives $\ngen =
\dim\im\Lambda \le \binom{r+1}{2} - 1$.
\end{proof}

\begin{remk}\label{rem:compactness-examples}
\begin{enumerate}
  \item Without single-wall orientation preservation in
    Theorem~\ref{thm:compactness-stress}(2), tangential fold-overs occur: for
    $\ngen = 1$, $L = I_{\dimin}$, $a_1 = 0$, $w_1 = e_1$ ($F(x) = x$, $\ell_1 =
    q_1 = 1$), $\mathcal{M}(G_0,\dimin) = \{0,-2\}$, where $\mu_1 = -2$ is the
    reflection fold-over $\widetilde{F}(x) = (-|x_1|, x_2, \dots, x_{\dimin})$
    (Remark~\ref{rem:PR-dichotomy}).
  \item In Corollary~\ref{cor:stress-rigidity}, the stress condition is
    incompatible with full-cube positive path rank when $r \le 2$: $\ngen \le
    \binom{r+1}{2} - 1$ rules out $r = 1$ and forces $\ngen = 2$, $n = \dimin$ at
    $r = 2$. With $\widetilde{C} := WL^{-1}A \in \RR^{2\times 2}$, $Z =
    [-L^{-1}A;\,I_2]$, and $S_\Omega = (s_{ij}) \succ 0$,~\eqref{eq:wall-congruence}
    gives $(1+\widetilde{C}_{11})s_{11} + \widetilde{C}_{12}s_{12} =
    \widetilde{C}_{21}s_{12} + (1+\widetilde{C}_{22})s_{22} = 0$, hence
    $(1+\widetilde{C}_{11})(1+\widetilde{C}_{22})s_{11}s_{22} =
    \widetilde{C}_{12}\widetilde{C}_{21}s_{12}^2$. Since $\det J(1,0)/\det L =
    1+\widetilde{C}_{11} > 0$, $\det J(0,1)/\det L = 1+\widetilde{C}_{22} > 0$,
    and $s_{11}s_{22} > s_{12}^2 > 0$, this forces $\det J(1,1)/\det L =
    (1+\widetilde{C}_{11})(1+\widetilde{C}_{22}) - \widetilde{C}_{12}\widetilde{C}_{21}
    < 0$.

    At $r = 3$, taking $(\dimin,\ngen,n) = (3,5,5)$, $A = I_5$, and $S_\Omega =
    I_3$ gives $Z = [I_3;\,-L]$ and $\tr(Z\transp E_h Z) = -2(L_h w_h +
    \norm{L_h}_2^2) = 0 \iff L_h w_h = -\norm{L_h}_2^2$ on rows $L_h, w_h\transp$.
    Set
    \[
      L = \begin{pmatrix} 0 & 1 & 1 \\ -1 & 0 & -1 \\ 1 & -1 & 0 \\ 1 & 1 & 1 \\ 0 & -1 & 0 \end{pmatrix},
      \qquad
      W = \begin{pmatrix} -2 & 0 & -2 \\ 3 & -2 & -1 \\ -2 & 0 & 1 \\ 0 & -3 & 0 \\ -2 & 1 & -1 \end{pmatrix},
    \]
    so $\det(L\transp L) = 8$, $L_h w_h = -\norm{L_h}_2^2$ ($1 \le h \le 5$), the
    five traceless matrices $Z\transp E_h Z$ have determinant $16 \ne 0$ in
    $(Y_{12}, Y_{13}, Y_{22}, Y_{23}, Y_{33})$ (injective $\Lambda$), and all ten
    $3$-minors of $W$ are non-zero. Full-cube positive path rank requires
    $\rank(L + \diag(s)W) = 3$ at every $s \in [0,1]^5$. It suffices that some
    $3\times 3$ minor is non-zero at each $s$, and three minors already cover
    the cube. Writing $J_{ijk}(s)$ for rows $\{i,j,k\}$ of
    $L + \diag(s)W$ on $[0,1]^5$,
    \begin{align*}
      \det J_{145}(s) &= -1 + 2s_5(1-s_1)(1-3s_4),\\
      \det J_{125}(s) &= 1 + s_2\bigl(4s_5(1-s_1) + 8s_1 - 3\bigr),\quad
      \det J_{235}(s)\big|_{s_2=1,\,s_5=1/2} = 3 + 3s_3.
    \end{align*}
    If $\det J_{145}(s) = 0$, then $1 - 3s_4 > 0$ and $s_5(1-s_1) =
    \frac{1}{2(1-3s_4)} \ge \frac{1}{2}$, forcing $\det J_{125}(s) \ge 1 - s_2 +
    8s_1 s_2 \ge 0$ with equality only at $(s_1,s_2,s_4,s_5) = (0,1,0,\frac{1}{2})$,
    where $\det J_{235}(s) = 3 + 3s_3 > 0$. Hence at every $s \in [0,1]^5$ one
    of the three displayed minors is non-zero, and full-cube positive path rank
    and $\mathcal{S}(G_0) = \{0\}$ hold simultaneously.
\end{enumerate}
\end{remk}

\begin{remk}\label{rem:bridge}
For $\rho \in \Pprof$, setting $s_h(x) := \rho'(\beta z_h(x)) \in [0,1]$ and
$\Upsilon(s) := [I_{\dimin};\,\diag(s)W]$ gives $g_\beta(x) =
\Upsilon(s(x))\transp G\,\Upsilon(s(x))$ and
\begin{equation}\label{eq:finite-beta-pencil}
  \Upsilon\bigl(s(x)\bigr)\transp\bigl(G(\mu) - G_0\bigr)\Upsilon\bigl(s(x)\bigr)
  = \sum_{h=1}^{\ngen} 2\mu_h\,s_h(x)\bigl(1-s_h(x)\bigr)\,w_h w_h\transp.
\end{equation}
For $\beta \ge \beta_\sharp$, evaluating~\eqref{eq:finite-beta-pencil} on each
wall-$h$ tube ($s_h \in (0,1)$, $s_{h'} \in \{0,1\}$ for $h' \ne h$) isolates
$\mu = 0$ uniquely from $g_\beta$. As $\beta \to \infty$, $s_h(1-s_h) \to 0$ on
open chambers, and $\mu$ is left free across $\mathcal{M}(G_0,n)$. Here
$h_{2k}(s)$ samples $(0,1)^{2k}$, where $\mu$ enters, whereas $K(\Phi_I)$ and
$\mathcal{M}(G_0,n)$ depend only on $\{0,1\}^{2k}$
(Proposition~\ref{prop:vertex-flex-rigid}(3)).
\end{remk}
\section{Discussion and Open Problems}\label{sec:conclusion}

We have shown that convex coordinate-wise smoothings of an arrangement immersion
concentrate their Gauss--Bonnet--Chern curvature onto the even-codimension flats
as a Hafnian density whose integral depends only on the incident chamber
metrics. In even dimension the global vertex sum is invariant under
translating the walls. On the realisation side, immersions sharing prescribed
chamber metrics form a compact rank-constrained spectrahedron whose generic
local dimension is governed by the corank $r = \max(0, \dimin+\ngen-n)$.
Theorem~\ref{thm:global-weak-star} relies on simplicity of $\arr(b)$ to
separate overlapping transition tubes at scale $\beta^{-1}$. Neither the
local orthant formula of Theorem~\ref{thm:banchoff-orthant} nor the global
vertex sum of Theorem~\ref{thm:offset-invariance}, however, requires
simplicity.

Several questions remain open. Does the weak-$*$ convergence of
Theorem~\ref{thm:global-weak-star}, with its local $O(\beta^{-1})$ rate, extend
to non-simple or central arrangements, where more than $\dimin$ walls meet
along a single flat? Outside the codimension-$2$-generic regime of
Remark~\ref{rem:q2-counterexamples}, can the chamber-nullspace dimension
$\dim\mathcal{N}(\arr)$ be expressed in terms of the affine oriented matroid of
$\arr$? In the same regime, when is the immersive moduli space $\mathcal{M}_{\mathrm{PR}}(G_0,n)$
connected? Certifying full-cube positive path rank~\eqref{eq:PR-global} is
co-NP-hard when $\dimin$ and $\ngen$ both grow, by the P-matrix
criterion~\eqref{eq:PR-pmatrix} and~\cite{coxson1994}. Is it tractable for
fixed $\dimin$, where $M_{\mathrm{P}} - I_{\ngen}$ has rank at most $\dimin$,
and what is its complexity in positive codimension $n > \dimin$? We close with the global counterpart of
Theorem~\ref{thm:generic-local-dim}.

\begin{openproblem}[Generic global uniqueness for $\ngen \le \binom{r+1}{2}$]
Theorem~\ref{thm:generic-local-dim} establishes local rigidity at $\mu = 0$ on
$\Gamma$ whenever $\ngen \le \binom{r+1}{2}$, and
Theorem~\ref{thm:compactness-stress}(2) gives global uniqueness in
$\mathcal{M}_{\mathrm{PR}}(G_0,\dimin)$ at $n = \dimin$. By contrast, the
unrestricted space $\mathcal{M}(G_0,\dimin)$ admits reflection fold-overs
already at $\ngen = 1$ (Remark~\ref{rem:compactness-examples}(1)).

At the
tight boundary $\ngen = \binom{r+1}{2}$, Proposition~\ref{prop:vertex-flex-rigid}(4)
($b_3 = -1$) exhibits a second realisation
$\mu^* \in \mathcal{M}_{\mathrm{face}}(G_0,3;b) \setminus
\mathcal{M}_{\mathrm{PR}}(G_0,3)$. It satisfies positive path rank on every
non-empty face of $\arr(b)$ but fails full-cube positive path
rank~\eqref{eq:PR-global} on the faces bordering the missing chamber. This
echoes Hendrickson's redundant-rigidity obstruction for bar
frameworks~\cite{hendrickson1992}. For generic
$(L_0,A_0,W) \in \Gamma$ satisfying full-cube positive path rank, does local
rigidity imply global uniqueness $\mathcal{M}_{\mathrm{PR}}(G_0,n) = \{0\}$
throughout the range $\ngen \le \binom{r+1}{2}$, in analogy with generic global
rigidity for bar frameworks~\cite{connelly2005, gortler2010}? In the
strictly overconstrained regime $\ngen \le \binom{r+1}{2} - 1$ with $r \ge 3$
(Remark~\ref{rem:compactness-examples}(2)), when does the stress certificate of
Corollary~\ref{cor:stress-rigidity} ensuring $\mathcal{S}(G_0) = \{0\}$ exist
generically?
\end{openproblem}

\backmatter

\section*{Declarations}
\bmhead{Funding}
No funds, grants, or other support were received for the preparation of
this manuscript.
\bmhead{Competing interests}
The author has no competing interests to declare that are relevant to the
content of this article.
\bmhead{Data availability}
This manuscript has no associated data. All numerical quantities appearing in
the examples are determined by explicit rational formulas or integer polynomial
equations derived in Appendix~\ref{app:moduli-proofs}.
\bmhead{AI use declaration}
Large language models (Google Gemini and Anthropic Claude) were used to assist
with drafting and editing manuscript prose, drafting and condensing proof
write-ups from the author's mathematical outlines, structural organisation, and
internal error-checking passes. Every mathematical definition, theorem
statement, construction, and proof was independently verified line by line by
the author, who takes full responsibility for the content of this article.

\PaperAppendices

\section{Auxiliary Bounds for Section~\ref{sec:banchoff}}
\label{app:banchoff-aux}

\begin{proposition}\label{prop:folding-bounds}
Assume local transversality and positive path rank at $I$, with $|I| = q \ge 1$.
Let $\varepsilon_a(\xi)$ be the indicators~\eqref{eq:ray-slopes} and
$\xi \sim \Norm(0,P_I)$.
\begin{enumerate}
  \item If $q$ is odd, $\prod_{a=1}^q\varepsilon_a(-\xi) =
    -\prod_{a=1}^q\varepsilon_a(\xi)$ and
    $\EE_\xi\bigl[\prod_{a=1}^q\varepsilon_a(\xi)\bigr] = 0$.
  \item Writing $\theta_a \in [0,\pi)$ for the angle between $r_a$ and
    $r_a + p_a$ in $\mathcal{T}_I^\perp$, $\EE_\xi|\varepsilon_a(\xi)| =
    \theta_a/\pi$, and for $q = 2k$, $|K(\Phi_I)| \le \min_{1\le a\le 2k}
    \theta_a/\pi < 1$ (so $K(\Phi_I) = 0$ if some $p_a = 0$).
  \item If $q = 2k$ and $\{1,\dots,2k\} = I_1 \sqcup I_2$ with $|I_1|, |I_2|$
    even and $\spn\{r_a, p_a : a \in I_1\} \perp \spn\{r_b, p_b : b \in I_2\}$
    in $\mathcal{T}_I^\perp$, then $K(\Phi_I) = K(\Phi_{I_1})\,K(\Phi_{I_2})$,
    where $K(\Phi_{I_j})$ is the Banchoff curvature of $\sum_{a\in I_j}(u_ar_a +
    (u_a)_+p_a)$.
\end{enumerate}
\end{proposition}

\begin{proof}
Item~1 follows from $\varepsilon_a(-\xi) = -\varepsilon_a(\xi)$ and symmetry of
$\xi$. For Item~2 we use Sheppard's formula~\cite{sheppard1899}, which gives
$\PP_\xi\{\ell_a^-(\xi)\ell_a^+(\xi) > 0\} = \theta_a/\pi$, with $\theta_a <
\pi$ since $X_a(s_a) \ne 0$ on $[0,1]$, together with
$|\prod_a\varepsilon_a| \le |\varepsilon_{a_0}|$. Item~3 follows from
independence of orthogonal Gaussian projections.
\end{proof}

\section{Proofs of Lemma~\ref{lem:tube-separation} and Theorem~\ref{thm:global-weak-star}}
\label{app:global-proofs}

Throughout this appendix we use $\lambda_* \le \lambda^*$ from
Section~\ref{subsec:weak-star}. Implicit constants in $O(\cdot)$ may depend on
$(\arr(b), L, A, \rho)$, and on $R_\varphi$ where we write $O_{R_\varphi}(\cdot)$,
but never on $\varphi$ or $\beta$.

\begin{proof}[Proof of Lemma~\ref{lem:tube-separation}]
If no such threshold exists, there are $\beta_n \to \infty$ and $x_n \in
\RR^{\dimin}$ admitting no admissible pair $(S^*, y)$ at scale $\beta_n$. Passing
to a subsequence, assume $T_{\beta_n}(x_n) = S$ and $(\sgn z_h(x_n))_{h=1}^{\ngen}$
are independent of $n$, and each $z_h(x_n)$ either converges to $0$ or satisfies
$\inf_n|z_h(x_n)| > 0$. Set $S' := \{h : z_h(x_n) \to 0\} \supseteq S$.

Suppose first that $S' \ne \emptyset$. Then $\Phi_{S'} \ne \emptyset$, since
otherwise $\inf_x\norm{W_{S'}x + b_{S'}}_\infty > 0$, and $W_{S'}$ has full row
rank by simplicity. The point $y_n := x_n - W_{S'}\pinv z_{S'}(x_n)$ lies in
$\Phi_{S'}$ and satisfies $\norm{x_n - y_n}_2 \le
\norm{W_{S'}\pinv}_{\mathrm{op}}\norm{z_{S'}(x_n)}_2 \to 0$. This holds even if
$\norm{x_n}_2 \to \infty$, because $x_n - y_n$ depends only on $z_{S'}(x_n)$.
If $S' = \emptyset$, put $y_n := x_n$. For $h \notin S'$, $\inf_n|z_h(x_n)| > 0$
and $|z_h(y_n) - z_h(x_n)| \to 0$ give $|z_h(x_n)| > \beta_n^{-1}$ and
$\sgn z_h(y_n) = \sgn z_h(x_n) \ne 0$ for large $n$. Thus $(S', y_n)$ is
admissible, a contradiction.

Now let $\beta \ge \beta_\sharp$ and let $(S^*, y)$ be any admissible pair for
$x$. Every $h \notin S^*$ has $|z_h(x)| > \beta^{-1}$ and $s_h = \rho'(\beta
z_h(x)) = \mathbf{1}\{z_h(x) > 0\} = \mathbf{1}\{z_h(y) > 0\}$, so $DF_\beta(x)
= J_{\mathsf{F}(y)}(s)$ with $I(\mathsf{F}(y)) = S^*$. Hence $|T_\beta(x)| \le
|S^*| \le \dimin$ with independent normals, and $\lambda_* I_{\dimin} \preceq
g_\beta(x) \preceq \lambda^* I_{\dimin}$.
\end{proof}

\begin{proof}[Proof of Theorem~\ref{thm:global-weak-star}]
For $\beta \ge \beta_\sharp$, Item~1 follows directly from
Lemma~\ref{lem:tube-separation}. Moreover, if $|z_h(x)| \le \beta^{-1}$
for every $h \in S$, then $S \subseteq T_\beta(x) \subseteq S^*$, so $\Phi_S
\supseteq \Phi_{S^*} \ne \emptyset$. Since the strict cores
$\{T_J^\circ(\beta)\}_{J \subseteq \{1,\dots,\ngen\}}$ partition $\RR^{\dimin}$
and $\dd\mu_\beta^{(2k)} \equiv 0$ on $T_\emptyset^\circ(\beta)$, we analyse the
non-empty strict cores according to the cardinality $p = |J|$.

\pstep{1 (Strict cores with $p < 2k$ and $p > 2k$)}
Consider first a connected subset $U_J \subset T_J^\circ(\beta) \ne \emptyset$
on which the off-$J$ sign pattern $\sigma := (\mathbf{1}\{z_h > 0\})_{h\notin J}$
is constant. Lemma~\ref{lem:tube-separation} supplies $S^* \supseteq J$ with
$\Phi_J \ne \emptyset$, $\rank W_J = p$, and $y \in \Phi_{S^*}$ such that
$J_0(\sigma) + \sum_{h\in J} s_h a_h w_h\transp$ has rank $\dimin$ throughout
$[0,1]^J$; here $J_0(\sigma) := L + \sum_{h\notin J}\sigma_h a_h w_h\transp$.
Every $x \in U_J \subset T_J^\circ(\beta)$ satisfies $(2\sigma_h - 1)\beta
z_h(x) > 1$, and hence $\rho_\beta(z_h(x)) = \sigma_h(z_h(x) -
\beta^{-1}c_\rho)$, for all $h \notin J$. Fix any basepoint $x_0 \in \Phi_J$,
irrespective of the values $z_h(x_0)$. Relative to $x_0$, the restriction
$F_\beta|_{U_J}$ matches the local model~\eqref{eq:local-model} with $\sigma$
in place of $\sigma^0$ and $c_\beta(\sigma) := Lx_0 + c_0 + \sum_{h\notin
J}\sigma_h(z_h(x_0) - \beta^{-1}c_\rho)a_h$. When $p < 2k$,
Lemmas~\ref{lem:shear-splitting} and~\ref{lem:translation-gauss} split $g_\beta$
into $g_{\Psi,\beta}(u) \oplus g_\perp$ with $g_\perp$ constant and $\dim u = p
< 2k$. Every $2k$-fold antisymmetric contraction
in~\eqref{eq:smoothed-gbc-measure} then contains at least one flat tangential
index, so $\dd\mu_\beta^{(2k)} \equiv 0$ on $T_J^\circ(\beta)$.

When $p > 2k$, every point $x$ in the truncated tube $T_J^\circ(\beta)\cap
[-R_\varphi,R_\varphi]^{\dimin}$ satisfies $|T_\beta(x)| \le \dimin$, $\lambda_*
I_{\dimin} \preceq g_\beta(x) \preceq \lambda^* I_{\dimin}$, $0 \le \rho''_\beta
\le \beta\norm{\rho''}_\infty$, and
\[
  \partial_a\partial_e F_\beta(x)
  = \sum_{h\in T_\beta(x)} \rho''_\beta(z_h(x))\, a_h w_{h,a}w_{h,e}.
\]
Since $g_\beta$ is the pullback metric of the immersion $F_\beta$, the Gauss
equation $R_{\beta,abef} = \ip{\mathrm{II}_{ae}}{\mathrm{II}_{bf}} -
\ip{\mathrm{II}_{af}}{\mathrm{II}_{be}}$ expresses its curvature through the
normal components of $\partial_a\partial_e F_\beta = O(\beta)$. Hence
$|R_{\beta,abef}(x)| = O(\beta^2)$. Raising indices with $g_\beta^{-1} \preceq
\lambda_*^{-1}I_{\dimin}$ gives $|R^{cd}_{\beta,ab}(x)| = O(\beta^2)$ and
$|\dd\mu_\beta^{(2k)}/\dd x| = O(\beta^{2k})$. On the other hand, Fubini's
theorem along $\ker W_J \oplus (\ker W_J)^\perp$ gives
$\vol_{\dimin}\bigl(T_J(\beta)\cap \allowbreak
[-R_\varphi,R_\varphi]^{\dimin}\bigr) = \allowbreak O_{R_\varphi}(\beta^{-p})$.
As $p \ge 2k+1$ and $\beta \ge 1$, we have $\beta^{2k-p} \le \beta^{-1}$, so
summing over all $J$ with $|J| > 2k$ contributes
$O_{R_\varphi}(\norm{\varphi}_\infty\beta^{-1})$.

\pstep{2 (Top degree $\dimin = 2k$)}
When $\dimin = 2k$ and $\beta \ge \beta_\sharp$, no strict core has $|J| >
\dimin$, and those with $|J| = \dimin$ are the vertex boxes $T_I(\beta) = v_I +
W_I^{-1}[-\beta^{-1},\beta^{-1}]^{\dimin}$. These boxes are pairwise disjoint,
because any $x \in T_I(\beta)\cap T_{I'}(\beta)$ with $I \ne I'$ would satisfy
$|T_\beta(x)| \ge |I\cup I'| > \dimin$. Under the affine chart $u = W_I(x - v_I)$,
diffeomorphism invariance of the unoriented density $\dd\mu_\beta^{(\dimin)}$
cancels $|\det W_I|$ against $\dd x = |\det W_I|^{-1}\dd u$, pulling
$\dd\mu_\beta^{(\dimin)}$ back to $e_{\dimin}(g_{\Psi,\beta})(u)$ on
$\Rtube_\beta$. Its integral over $\Rtube_\beta$ equals $K(v_I)$ by
Proposition~\ref{prop:hafnian-euler}(3) and Theorem~\ref{thm:banchoff-orthant}.
This proves~\eqref{eq:vertex-box-mass}. Item~2 follows by combining
$|\varphi(x) - \varphi(v_I)| = O(\operatorname{Lip}(\varphi)\beta^{-1})$ on
$T_I(\beta)$ with $\int_{T_I(\beta)}|\dd\mu_\beta^{(\dimin)}| =
\int_{[0,1]^{\dimin}} |h_{\dimin}(s)|\,\dd s < \infty$.

\pstep{3 (Shaved shear tubes for $2k < \dimin$)}
Assume now that $2k < \dimin$. Fix $I$ with $|I| = 2k$ and an open
$(\dimin-2k)$-cell $\Phi_{I,\alpha} \subset \Phi_I^{\mathrm{reg}}$ with off-$I$
sign pattern $\sigma^0_{I,\alpha}$, basepoint $x_0 \in \Phi_{I,\alpha}$, and
orthonormal basis $Q_{\perp,I}$ of $T\Phi_I = \ker W_I$. Composing
$\Theta_{\beta,I}^{-1}(u,y) = (u, y - \sum_{a=1}^{2k}\varphi_{a,\beta}(u_a))$
from Lemma~\ref{lem:shear-splitting} (built from $J_{0,I,\alpha} :=
J_0(\sigma^0_{I,\alpha})$) with the adapted affine chart $\iota_I$ of
Definition~\ref{def:local-model}(1) defines the ambient diffeomorphism
\[
  x_{\beta,I}(u,y)
  := (\iota_I\circ\Theta_{\beta,I}^{-1})(u,y)
   = x_0 + V_I u + Q_{\perp,I}\Bigl(y - \sum_{a=1}^{2k}\varphi_{a,\beta}(u_a)\Bigr)
  \in \RR^{\dimin},
\]
with global inverse $x_{\beta,I}^{-1}(x) = \bigl(u,\, Q_{\perp,I}\transp(x-x_0)
+ \sum_{a=1}^{2k}\varphi_{a,\beta}(u_a)\bigr)$ for $u = W_I(x-x_0) = z_I(x)$
(since $Q_{\perp,I}\transp V_I = 0$).

On any $\beta$-admissible region, Lemma~\ref{lem:shear-splitting}(2) gives
$(x_{\beta,I})^*g_\beta(u,y) = g_{\Psi,\beta}(u) \oplus g_{\perp,I,\alpha}$,
and~\eqref{eq:smoothed-gbc-measure} is invariant under diffeomorphisms, so we
may evaluate it in the coordinates $(u,y)$. Because $g_{\perp,I,\alpha}$ is
constant, the Riemann tensor of the product metric is supported on the
$u$-block, so in~\eqref{eq:smoothed-gbc-measure} the generalised Kronecker
delta restricts both index strings to permutations $\pi, \pi'$ of
$\{1,\dots,2k\}$ with $\delta^{\pi(1)\cdots\pi(2k)}_{\pi'(1)\cdots\pi'(2k)} =
\sgn(\pi)\sgn(\pi')$. We raise the indices of each factor,
$R^{j_{2\ell-1}j_{2\ell}}_{i_{2\ell-1}i_{2\ell}} = \sum_{c,d}
g_{\Psi,\beta}^{j_{2\ell-1}c}g_{\Psi,\beta}^{j_{2\ell}d}R_{i_{2\ell-1}i_{2\ell}cd}$,
and sum over $\pi'$. This produces $\det(g_{\Psi,\beta}^{-1})\,\mathsf{P}_{2k}(R)$,
which combines with $\sqrt{\det(g_{\Psi,\beta}\oplus g_{\perp,I,\alpha})} =
\sqrt{\det g_{\Psi,\beta}}\sqrt{\det g_{\perp,I,\alpha}}$ to give the
$(\det g_{\Psi,\beta})^{-1/2}\mathsf{P}_{2k}(R)$ of~\eqref{eq:coord-pfaffian}.
Hence the measure pulls back as $(x_{\beta,I})^*\dd\mu_\beta^{(2k)}(u,y) =
e_{2k}(g_{\Psi,\beta})(u)\sqrt{\det g_{\perp,I,\alpha}}\,\dd y$ (with
$e_{2k}(g_{\Psi,\beta})$ denoting the unoriented density).

We next choose the region on which this chart is used. Since $|u_a| \le
\beta^{-1}$ and $0 \le \rho_\beta(u_a) \le 2\beta^{-1}$ on
$\Rtube_\beta$ (as $0 \le \rho \le 2$ on $[-1,1]$), the constant
\[
  C_1 := \max_{(I,\alpha)}\Bigl(\sqrt{2k}\,\norm{V_I}_{\mathrm{op}}
  + \norm{(J_{0,I,\alpha}Q_{\perp,I})\pinv}_{\mathrm{op}}
    \sum_{a=1}^{2k}\bigl(\norm{J_{0,I,\alpha}v_a}_2 + 2\norm{a_{h_a}}_2\bigr)\Bigr)
\]
satisfies, for every $(I,\alpha)$, every $x \in T_I(\beta)$ with $u = z_I(x)
\in \Rtube_\beta$ and projection $x_I := x - V_I u \in \Phi_I$, and every shear
basepoint $p(y) := x_0 + Q_{\perp,I}y \in \Phi_I$,
\[
  \norm{x - x_I}_2 \le C_1\beta^{-1},
  \qquad
  \norm{x_{\beta,I}(u,y) - p(y)}_2
  = \Bigl\lVert V_I u - Q_{\perp,I}\sum_{a=1}^{2k}\varphi_{a,\beta}(u_a)\Bigr\rVert_2
  \le C_1\beta^{-1}.
\]
To keep the sheared tube inside the strict core $T_I^\circ(\beta)$ we shave a
collar off each cell. Write $\partial\Phi_{I,\alpha}$ for the relative boundary
of $\Phi_{I,\alpha}$ in $\Phi_I$, with the convention
$\dist(\cdot,\emptyset) := +\infty$. For $h \notin I$, the
restriction of the affine function $z_h$ to the affine subspace $\Phi_I$ has
gradient $P_{T\Phi_I}w_h$. If $\Phi_{I\cup\{h\}} = \Phi_I\cap\wall{h} \ne
\emptyset$, simplicity gives $c_{I,h} := \norm{P_{T\Phi_I}w_h}_2 > 0$ and
$|z_h(q)| = c_{I,h}\dist(q,\Phi_{I\cup\{h\}})$ for $q \in \Phi_I$. Since the
hyperplane $\Phi_{I\cup\{h\}}$ of $\Phi_I$ is disjoint from the cell
$\Phi_{I,\alpha}$, every $p \in \Phi_{I,\alpha}$ has $\dist(p,\Phi_{I\cup\{h\}})
\ge \dist(p,\partial\Phi_{I,\alpha})$. If instead $\Phi_{I\cup\{h\}} =
\emptyset$, then $z_h$ is a non-zero constant $z_h(\Phi_I)$ on $\Phi_I$. Set
\[
  C_2 := 1 + \max\Bigl\{\max_{\Phi_{I\cup\{h\}}\ne\emptyset}
  \frac{1 + \norm{w_h}_2 C_1}{c_{I,h}},\
  \max_{\Phi_{I\cup\{h\}}=\emptyset}\frac{1 + \norm{w_h}_2 C_1}{|z_h(\Phi_I)|}\Bigr\},
  \quad
  \beta_0 := \max(\beta_\sharp, C_1 + C_2),
\]
the maxima being taken over all $(I,\alpha)$ and $h \notin I$, with
$\max\emptyset := 0$. Let $\beta \ge \beta_0$, and let $p(y)$ lie in the
shaved cell $\Phi_{I,\alpha}^{\mathrm{sh}}(\beta) := \{p \in \Phi_{I,\alpha} :
\dist(p,\partial\Phi_{I,\alpha}) \ge C_2\beta^{-1}\}$. Then $|z_h(p(y))| >
(1 + \norm{w_h}_2C_1)\beta^{-1}$ for every $h \notin I$; in the second case
this uses $\beta > C_2$. Along the segment $[p(y), x_{\beta,I}(u,y)]$ with
$u \in \Rtube_\beta$, $z_h$ varies by at most $\norm{w_h}_2C_1\beta^{-1}$, so
this segment avoids $T_{\{h\}}(\beta)$ for every $h \notin I$. The endpoint
$x_{\beta,I}(u,y)$ therefore carries the off-$I$ sign pattern of
$\Phi_{I,\alpha}$, with $|z_h| > \beta^{-1}$ for $h \notin I$ and $|z_{h_a}| =
|u_a| \le \beta^{-1}$ for $a \le 2k$. So $T_\beta(x_{\beta,I}(u,y)) = I$, and
the region is $\beta$-admissible at $I$.

Consequently the shaved shear tubes
\[
  \mathcal{T}_{I,\alpha}^{\mathrm{sh}}(\beta)
  := x_{\beta,I}\Bigl(\Rtube_\beta \times \bigl\{y \in \RR^{\dimin-2k} :
     p(y) \in \Phi_{I,\alpha}^{\mathrm{sh}}(\beta)\bigr\}\Bigr)
\]
lie in $T_I^\circ(\beta)$ with off-$I$ sign pattern $\sigma^0_{I,\alpha}$ and
are pairwise disjoint across $(I,\alpha)$. Indeed, the strict cores
$T_I^\circ(\beta)$ are disjoint for distinct $I$, and distinct cells $\alpha$
have distinct off-$I$ sign patterns. Since $C_1\beta^{-1} < 1$, we also have
$\varphi(x_{\beta,I}(u,y)) = \varphi(p(y)) = 0$ whenever $p(y) \notin
B_{R_\varphi} := B(0,\sqrt{\dimin} R_\varphi + 1)$.

On $\Phi_{I,\alpha}^{\mathrm{sh}}(\beta) \cap B_{R_\varphi}$,
Theorem~\ref{thm:banchoff-orthant} gives $\int_{\Rtube_\beta}
e_{2k}(g_{\Psi,\beta})(u)\,\dd u = K(\Phi_{I,\alpha})$. Using
$\int_{\Rtube_\beta}|e_{2k}(g_{\Psi,\beta})|\,\dd u = \int_{[0,1]^{2k}}
|h_{2k}(s)|\,\dd s < \infty$, $|\varphi(x_{\beta,I}(u,y)) - \varphi(p(y))| \le
C_1\operatorname{Lip}(\varphi)\beta^{-1}$, and $\Haus{\dimin-2k}(\Phi_I\cap
B_{R_\varphi}) < \infty$, we obtain
\[
  \left|\int_{\mathcal{T}_{I,\alpha}^{\mathrm{sh}}(\beta)} \varphi\,\dd\mu_\beta^{(2k)}
        - \int_{\Phi_{I,\alpha}^{\mathrm{sh}}(\beta)}
          \varphi\,K(\Phi_{I,\alpha})\sqrt{\det g_{\perp,I,\alpha}}\,\dd\Haus{\dimin-2k}\right|
  = O_{R_\varphi}\bigl(\operatorname{Lip}(\varphi)\beta^{-1}\bigr).
\]

\pstep{4 (Collar remainder around lower-dimensional flats)}
Take any $\beta \ge \beta_0$ and any $x \in (T_I^\circ(\beta)\cap
[-R_\varphi,R_\varphi]^{\dimin}) \setminus \bigcup_\alpha
\mathcal{T}_{I,\alpha}^{\mathrm{sh}}(\beta)$. There are two cases. If the
segment $[x,x_I]$ meets some $\wall{h}$ with $h \notin I$, then
$\Phi_{I\cup\{h\}} \ne \emptyset$, and $|z_h(x_I)| \le \norm{w_h}_2
C_1\beta^{-1}$ implies $\dist(x_I,\Phi_{I\cup\{h\}}) = O(\beta^{-1})$.
Otherwise $x_I$ shares the off-$I$ sign pattern of $x$ and lies in a unique
cell $\Phi_{I,\alpha}$. Its shear basepoint $p(y)$ is not in
$\Phi_{I,\alpha}^{\mathrm{sh}}(\beta)$ and satisfies $\norm{p(y) - x_I}_2 \le
2C_1\beta^{-1}$, which gives $\dist(x_I,\partial\Phi_{I,\alpha}) < (C_2 +
2C_1)\beta^{-1}$. In either case, both $x_I$ and $(\Phi_{I,\alpha}\setminus
\Phi_{I,\alpha}^{\mathrm{sh}}(\beta))\cap [-R_\varphi,R_\varphi]^{\dimin}$ lie
within distance $O(\beta^{-1})$ of the $(\dimin-2k-1)$-skeleton
$\bigcup_{h\notin I}(\Phi_I\cap\wall{h})$ inside $\Phi_I\cap B_{R_\varphi}$.

This skeleton is a union of at most $\ngen - 2k$ affine
$(\dimin-2k-1)$-subspaces of $\Phi_I$, so its $O(\beta^{-1})$ collar
$K_I(\beta)$ in $\Phi_I\cap B_{R_\varphi}$ has $\Haus{\dimin-2k}$-measure
$O_{R_\varphi}(\beta^{-1})$. We now apply Fubini's theorem in the adapted chart
$\iota_I(u,x_\perp) = x_0 + V_I u + Q_{\perp,I}x_\perp$, where $u \in
\Rtube_\beta$ and $|\det D\iota_I| = (\det(W_I W_I\transp))^{-1/2}$. The
smoothed collar set then has $\dimin$-dimensional volume at most $|\det
D\iota_I|(2\beta^{-1})^{2k} \Haus{\dimin-2k}(K_I(\beta)) =
O_{R_\varphi}(\beta^{-2k-1})$. On the smoothed side we combine this with
$|\dd\mu_\beta^{(2k)}/\dd x| = O(\beta^{2k})$; on the unsmoothed side we use
$|K(\Phi_{I,\alpha})| < 1$ from Proposition~\ref{prop:folding-bounds}(2). Both
collar remainders are therefore $O_{R_\varphi}(\norm{\varphi}_\infty\beta^{-1})$
for $\beta \ge \beta_0$. Uniform boundedness of $|\dd\mu_\beta^{(2k)}/\dd x|$ on
$[\beta_\sharp, \beta_0]$ absorbs $\beta \in [\beta_\sharp, \beta_0]$ into
$C(R_\varphi)$, proving~\eqref{eq:weak-star-rate} for all $\beta \ge
\beta_\sharp$. Tracking the dependence on $R_\varphi$ through
$\Haus{\dimin-2k}(\Phi_I\cap B_{R_\varphi})$, the collar measure, and the
volume bound of Step~1 shows that $C(R_\varphi) = O(R_\varphi^{\dimin-2k})$.
\end{proof}

\section{Moduli Genericity Proofs and Degree-Four Flex Verification}
\label{app:moduli-proofs}

\begin{lemma}\label{lem:obstruction-genericity}
Fix $\dimin \ge 1$, $\ngen \ge 1$, $n \ge \dimin$, and $r := \max(0,\,
\dimin + \ngen - n)$. For $B = [L\ A] \in \RR^{n\times(\dimin+\ngen)}$ of rank
$\min(n,\dimin+\ngen)$ and $W \in \RR^{\ngen\times\dimin}$ with kernel basis
$Z(B) \in \RR^{(\dimin+\ngen)\times r}$, $\rank\Lambda_{(B,W)} \le
\min\bigl(\ngen,\, \binom{r+1}{2}\bigr)$, and there exists a non-empty
Zariski-open set $\Gamma \subset \RR^{n\times\dimin} \times \RR^{n\times\ngen}
\times \RR^{\ngen\times\dimin}$ on which $\rank[L\ A] = \min(n,\dimin+\ngen)$,
$\rank L = \dimin$, $w_h \ne 0$ for all $h$, and
\begin{equation}\label{eq:obstruction-max-rank}
  \rank\Lambda_{(B,W)}
  = \dim\spn\bigl\{Z(B)\transp E_h(W) Z(B)\bigr\}_{h=1}^{\ngen}
  = \min\!\left(\ngen,\, \binom{r+1}{2}\right).
\end{equation}
\end{lemma}

\begin{proof}
Put $d_r := \min(\ngen,\binom{r+1}{2}) \ge \rank\Lambda_{(B,W)}$. When $r = 0$,
$\Gamma := \{\det(L\transp L)\det(B\transp B)\prod_{h}\norm{w_h}^2 \ne 0\}$
suffices. For $r \ge 1$ and each $n$-subset $J \subset \{1,\dots,\dimin+\ngen\}$,
set $\delta_J(B) := \det B_{\cdot,J}$ and $(Z_J)_{J^c,\cdot} := \delta_J(B)I_r$,
$(Z_J)_{J,\cdot} := -\operatorname{adj}(B_{\cdot,J})\,B_{\cdot,J^c}$ (spanning
$\ker B$ on $\{\delta_J(B) \ne 0\}$). Writing $c_J(L,A,W)$ for the sum of squares
of the $d_r\times d_r$ minors of $[Z_J\transp E_h(W)Z_J]_{h=1}^{\ngen}$,
Lemma~\ref{lem:wall-congruence} gives~\eqref{eq:obstruction-max-rank} on $\Gamma
:= \{P_{\mathrm{gen}} \ne 0\}$ for $P_{\mathrm{gen}}(L,A,W) := \det(L\transp L)
(\prod_{h}\norm{w_h}^2)\sum_{|J|=n}(\delta_J(B)c_J(L,A,W))^2$.

It remains to exhibit one point of $\Gamma$. Taking $Z_L^* := 0$
gives $Z\transp E_h Z\big|_{Z_L=0} = -2\eta_h\eta_h\transp$. Since $n \ge
\dimin$ we have $r \le \ngen$, so we may take the rows of $Z_A^* \in
\RR^{\ngen\times r}$ to be $e_i$ ($1 \le i \le r$), then
$\min(\ngen - r, \binom{r}{2})$ vectors $e_i + e_j$ ($i < j$), and $0$ for the
remaining $\ngen - d_r$ rows. This choice spans $d_r$ dimensions in $\Spec^r$
and has $Z_A^{*\transp}Z_A^* \succeq I_r$. Now take $C_A \in
\RR^{(n-\dimin)\times\ngen}$ with $\ker C_A = \im Z_A^*$, $L^* :=
[I_{\dimin};\,0]$, $A^* := [0;\,C_A]$, and $w_h^* \ne 0$. When $n = \dimin$
this means $r = \ngen$, $Z_A^* = I_{\ngen}$ and $A^* = 0$. Then $\ker[L^*\ A^*]
= \im[0;\,Z_A^*]$ and $P_{\mathrm{gen}}(L^*,A^*,W^*) > 0$.
\end{proof}

\begin{lemma}\label{lem:local-moduli-manifold}
Let $B_0 = [L_0\ A_0] \in \RR^{n\times(\dimin+\ngen)}$ have $\rank B_0 =
\min(n, \dimin+\ngen)$, $G_0 := B_0\transp B_0$, $r := \max(0, \dimin+\ngen-n)
= \dim\ker G_0$, $\im Z = \ker G_0$, $\Lambda(\mu) := \sum_{h=1}^{\ngen} \mu_h
Z\transp E_h Z$, and $r_\Lambda := \rank\Lambda$. On $\mathcal{M}_{\mathcal{V}} :=
\mathcal{M}(G_0,n) \cap \mathcal{V}$ for an open neighbourhood $\mathcal{V}
\subset \RR^{\ngen}$ of $0$:
\begin{enumerate}
  \item If $r = 0$, $\mathcal{M}_{\mathcal{V}} = \mathcal{V}$ has dimension $\ngen$.
  \item If $r \ge 1$ and $r_\Lambda = \binom{r+1}{2}$, $\mathcal{M}_{\mathcal{V}}$ is
    a real-analytic submanifold of dimension $\ngen - \binom{r+1}{2}$ with
    $T_0\mathcal{M}_{\mathcal{V}} = \ker\Lambda$.
  \item If $r \ge 1$ and $r_\Lambda = \ngen$, $\mathcal{M}_{\mathcal{V}} = \{0\}$.
  \item In all cases, $\dim\ker\Lambda = \ngen - r_\Lambda$ and
    $\mathcal{M}_{\mathcal{V}}$ lies in a real-analytic submanifold of dimension
    $\ngen - r_\Lambda$ with tangent space $\ker\Lambda$ at $0$.
\end{enumerate}
\end{lemma}

\begin{proof}
When $r = 0$, $G_0 \succ 0$ is an open condition, so shrinking $\mathcal{V}$ so
that $G(\mu) \succ 0$ on $\mathcal{V}$ gives $\mathcal{M}_{\mathcal{V}} =
\mathcal{V}$ in Item~1. Suppose $r \ge 1$. Then $\rank G_0 = n$, and by
continuity of $\mu \mapsto G(\mu)$ we may shrink $\mathcal{V} \ni 0$ so that
$G(\mu)$ lies in the open domain $\mathcal{U}$ of
Lemma~\ref{lem:psd-schur-chart} for all $\mu \in \mathcal{V}$. On
$\mathcal{V}$, Lemma~\ref{lem:psd-schur-chart} equates $\mu \in
\mathcal{M}(G_0, n)$ ($G(\mu) \succeq 0$ and $\rank G(\mu) \le n$) with
$f_{\ker}(\mu) := \Sigma_{\ker}(G(\mu)) = 0$. Here $f_{\ker} : \mathcal{V} \to
\Spec^r$ is real-analytic, $f_{\ker}(0) = 0$, and $Df_{\ker}(0) = \Lambda$
by~\eqref{eq:schur-chart-diff}. Choose a coordinate projection $\pi : \Spec^r
\to \RR^{r_\Lambda}$ that restricts to an isomorphism on $\im\Lambda \subset
\Spec^r$. Then $\rank D(\pi\circ f_{\ker})(0) = \rank(\pi\circ\Lambda) =
r_\Lambda$, and $\ker(\pi\circ\Lambda) = \ker\Lambda$ has dimension $\ngen -
r_\Lambda$. By the real-analytic implicit function theorem, after shrinking
$\mathcal{V}$ further, $(\pi\circ f_{\ker})^{-1}(0)$ is a real-analytic
submanifold of $\mathcal{V}$ of dimension $\ngen - r_\Lambda$ with tangent space
$\ker\Lambda$ at $0$, containing $\mathcal{M}_{\mathcal{V}} = f_{\ker}^{-1}(0)$
and proving Item~4. In Item~2, $r_\Lambda = \binom{r+1}{2} = \dim\Spec^r$, so
$\pi$ is an isomorphism and $\mathcal{M}_{\mathcal{V}} = f_{\ker}^{-1}(0) =
(\pi\circ f_{\ker})^{-1}(0)$; in Item~3, $r_\Lambda = \ngen$, so $(\pi\circ
f_{\ker})^{-1}(0) = \{0\}$ and $\mathcal{M}_{\mathcal{V}} = \{0\}$ since $0 \in
\mathcal{M}_{\mathcal{V}}$.
\end{proof}

\begin{lemma}\label{lem:joint-genericity}
Fix $\dimin \ge 2$, $\ngen \ge 1$, and $n \ge \dimin$. Then:
\begin{enumerate}
  \item Matrices $W \in \RR^{\ngen\times\dimin}$ with non-zero, pairwise
    non-parallel rows form a non-empty Zariski-open set.
  \item For any $W$, non-simple offsets $b \in \RR^{\ngen}$ form the finite
    measure-zero union $\bigcup \{\pi_S^{-1}(\im W_S) : \emptyset \ne S \subseteq
    \{1,\dots,\ngen\},\, \rank W_S < |S|\}$ ($\pi_S(b) := b_S$).
  \item Hence $(L, A, W) \in \Gamma$, simplicity of $\arr(b)$, and codimension-$2$
    genericity hold jointly for Lebesgue-a.e.\ $(L, A, W, b)$.
  \item If $b = 0$ and $\ngen > \dimin$, $\arr(0)$ is not simple since $\rank W \le
    \dimin < \ngen$ at $0 \in \bigcap_h\wall{h}$.
\end{enumerate}
\end{lemma}

\begin{proof}
Item~1 holds via $w_h = e_1 + h e_2$. For Item~2, non-simplicity at $x$ means
$\rank W_S < |S|$ for $S = I(\mathsf{F}(x))$ and $b_S = -W_S x \in \im W_S
\subsetneq \RR^S$. Item~3 follows from Theorem~\ref{thm:gram-spectrahedron}(3),
Lemma~\ref{lem:obstruction-genericity}, and Fubini's theorem, and Item~4 is
checked at $0$.
\end{proof}

\begin{proof}[Proof of Proposition~\ref{prop:vertex-flex-rigid}]
\textit{(1).}\
Clearing denominators by $(1+t^2)(9t^2+1)$ gives $4\sin\phi\sin\psi =
3(\cos\phi+\cos\psi)$, $3\sin\phi + 4\cos\phi\sin\psi = 5\sin\psi$, and
$3\sin\psi + 4\sin\phi\cos\psi = 5\sin\phi$, proving $c_1(t)\transp c_2(t) = 0$;
the chamber metrics and the formulas for $\mu_1(t), \mu_2(t)$ follow by
expanding $G_t = B_t\transp B_t$.

\ppart{(2).}
For $J(s;t) = [e_1+s_1 a_1(t)\ \ e_2+s_2 a_2(t)]$, the $(1,3)$- and
$(2,3)$-minors $m_{13} = \frac{4}{5}s_2(1+s_1)\sin\psi(t)$ and $m_{23} =
-\frac{8}{5}s_1\sin\phi(t)$ vanish together on $[0,1]^2$ only at $s = (0,0)$
($\rank L = 2$). Moreover, $\det[e_1,e_2,a_1(t)] = \frac{8}{5}\sin\phi(t) > 0$,
$z_t = (-\sin\psi,0,\sin\psi,-2\sin\phi)\transp$ spans $\ker B_t$, and
$z_t\transp E_1 z_t = -4\sin^2\psi < 0$. In the coordinates $x_1 := \frac{3}{5}
+ \frac{5}{8}\mu_1$ and $x_2 := \frac{5}{13} + \frac{5}{4}\mu_2$ (so that $x_1
= \cos\phi$ and $x_2 = \cos\psi$ along the flex),
Theorem~\ref{thm:compactness-stress}(1) gives
\[
  S_L(\mu) = \frac{8}{25}\begin{pmatrix} 8(1-x_1^2) & 3(x_1+x_2) \\
  3(x_1+x_2) & 2(1-x_2^2) \end{pmatrix},
\]
which cuts out the nodal quartic $16(1-x_1^2)(1-x_2^2) = 9(x_1+x_2)^2$ for
$x_1,x_2 \in [-1,1]$. Its ordinary double points are at $(1,-1)$ and $(-1,1)$,
with branch slopes $-9$ and $-1/9$. Setting $x_1 = \cos\phi$, $x_2 = \cos\psi$
modulo $(\phi,\psi)\mapsto(-\phi,-\psi)$ splits $\phi \in (0,\pi)$ into $3tu =
1$ and $tu = -3$. On the fold-over branch $tu = -3$, $m_{13} =
\frac{4}{5}s_2(1-3s_1)\sin\psi$ and
$m_{23} = -\frac{8}{5}s_1(1-2s_2)\sin\phi$ vanish at $s =
(\frac{1}{3},\frac{1}{2})$, where $\rank J = 1$. As $t \to 0^+$ or $t \to
\infty$, instead, row~$3$ of $B_t$ vanishes and $\rank J$ drops at
$(0,\frac{5}{9})$ and $(\frac{5}{13},0)$. Since $\ip{J_{10}e_1}{J_{00}e_1} =
\frac{8}{5}\cos\phi(t)$ is strictly decreasing, $B_t \not\sim_{O(3)} B_{t'}$
($t \ne t'$).

By contrast, the piecewise-linear homeomorphism $\Psi_{\mathrm{swap}}(x_1,x_2)
:= \bigl(x_2 - \frac{1}{2}(x_2)_+,\, x_1 + (x_1)_+\bigr)$ swaps the two
coordinate walls. Its differentials on $\cell{00},\cell{10},\cell{01},\cell{11}$
are
$\bigl(\begin{smallmatrix}0&1\\1&0\end{smallmatrix}\bigr)$,
$\bigl(\begin{smallmatrix}0&1\\2&0\end{smallmatrix}\bigr)$,
$\bigl(\begin{smallmatrix}0&1/2\\1&0\end{smallmatrix}\bigr)$,
$\bigl(\begin{smallmatrix}0&1/2\\2&0\end{smallmatrix}\bigr)$, and it pulls
back the chamber metrics of $F_{1/(3t)}$ to those of $F_t$. It also matches
$(e_1,e_2,c_1(t),c_2(t))$ with
$(e_2,e_1,2c_2(\frac{1}{3t}),\frac{1}{2}c_1(\frac{1}{3t}))$ under the
coordinate reflection $Q =
\diag\bigl(\bigl(\begin{smallmatrix}0&1\\1&0\end{smallmatrix}\bigr), 1\bigr)
\in O(3)\setminus SO(3)$ ($\det Q = -1$, since $\phi(\frac{1}{3t}) =
\psi(t)$). Hence $F_t(\RR^2) \cong_{O(3)} F_{1/(3t)}(\RR^2)$.

\ppart{(3).}
Since the chamber metrics are constant along the flex, so are the $M_\tau$,
and~\eqref{eq:codim2-arcsin} gives $K(0) = \frac{1}{\pi}\arcsin(3/5)$. In
codimension one ($m = 1$, $I = \{1,2\}$, $P_I = I_3$), $\Pi^\perp(s)$ is the
orthogonal projector onto the line spanned by $X_1(s_1)\times X_2(s_2)$ and
$\det M(s) = \norm{X_1(s_1)\times X_2(s_2)}_2^2$, so $N_{12}(s) =
\det[X_1,X_2,a_1]\det[X_1,X_2,a_2]/\det M(s)$ and
\[
  h_2(s)
  = \frac{\det[X_1,X_2,a_1]\,\det[X_1,X_2,a_2]}{2\pi\,(\det M(s))^{3/2}}.
\]
At $s^* = (\frac{1}{2},\frac{1}{2})$ both $t = \frac{1}{3}$ and $t = 1$ have
$\det[X_1,X_2,a_1]\det[X_1,X_2,a_2] = \frac{144}{125}$, whereas
$\det M(s^*;\frac{1}{3}) = \frac{297}{400} < \frac{329}{400} =
\det M(s^*;1)$, whence $h_{2,1/3}(s^*) > h_{2,1}(s^*)$.

\ppart{(4).}
Fix $t_0 = \frac{1}{2}$, $w_3 = (1,1)\transp$, and $a_3 = e_3$, so that $a_1 =
(-\frac{1}{25}, -\frac{6}{5}, \frac{32}{25})\transp$ and $a_2 = (-\frac{3}{5},
-\frac{9}{13}, \frac{48}{65})\transp$. The columns of $J(s)$ are $c^{(1)}(s) =
(1-\frac{1}{25}s_1, -\frac{6}{5}s_1, \frac{32}{25}s_1+s_3)\transp$ and
$c^{(2)}(s) = (-\frac{3}{5}s_2, 1-\frac{9}{13}s_2,
\frac{48}{65}s_2+s_3)\transp$. If $c^{(2)}(s) = \lambda c^{(1)}(s)$ on
$[0,1]^3$, rows~$1,2$ force $s_2 > 0$, $\lambda < 0$, $s_1 > 0$, which
contradicts row~$3$. So full-cube positive path rank holds. In the integer
basis $z_1 = (15,0,-15,26,0)\transp$, $z_2 = (-1,-30,-25,0,32)\transp$ of
$\ker B$, \eqref{eq:wall-congruence} gives the $(Y_{11},Y_{12},Y_{22})$-rows
$30(-30,-37,-40)$, $52(-26,-15,0)$, and $96(0,5,-42)$. Their reduced $3\times 3$
integer matrix $\begin{pmatrix} -30 & -37 & -40 \\ -26 & -15 & 0 \\
0 & 5 & -42 \end{pmatrix}$ has determinant $26{,}704 \ne 0$, so $\Lambda$ is
injective and $\dim_0\mathcal{M}(G_0,3) = 0$.

We next locate the points of $\mathcal{M}(G_0,3)$ that can be immersive. Write
$\mu' := (\mu_1,\mu_2)$. By Theorem~\ref{thm:compactness-stress}(1) (here
$L\transp L = I_2$), $\mu \in \mathcal{M}(G_0,3)$ if and only if $S_L(\mu)
\succeq 0$ and $\rank S_L(\mu) \le 1$. As $E_3$ only involves the row and
column of $a_3$, the leading $2\times 2$ block of $S_L(\mu)$ is the Schur
complement of Item~2 evaluated at $\mu'$. On the face $s_3 = 0$, which
corresponds to the vertex $v_{12} = 0$ of $\arr(b)$ where $z_3 = -1 < 0$, the
realisation of $G(\mu)$ restricts to the $(2,2,3)$ configuration with
parameters $\mu'$.
Positive path rank at $v_{12}$ therefore forces $\mu'$ onto the admissible arc
of Item~2: $\mu' = (\mu_1(t),\mu_2(t))$ with $3tu = 1$ and $t > 0$. Both
$\mathcal{M}_{\mathrm{PR}}(G_0,3)$ and $\mathcal{M}_{\mathrm{face}}(G_0,3;b)$
are thus contained in the set of $\mu = (\mu_1(t),\mu_2(t),\mu_3)$ with $t >
0$, $\mu_3 \in \RR$, and $\rank S_L(\mu) \le 1$, which we now determine. Along
the arc, Item~2 gives $S_{11} = \sigma_1^2$, $S_{12} = \sigma_1\sigma_2$, and
$S_{22} = \sigma_2^2$ with
\[
  \sigma_1 := \tfrac{8}{5}\sin\phi(t) = \frac{16t}{5(1+t^2)},
  \qquad
  \sigma_2 := \tfrac{4}{5}\sin\psi(t) = \frac{24t}{5(9t^2+1)}
\]
(the heights of the crease rays $c_1(t)$, $c_2(t)$ above the plane of $L$).
For the remaining entries,~\eqref{eq:schur-entries} with $L\transp a_3 = 0$,
$w_3\transp L\transp a_1 = -\frac{31}{25}$, and $w_3\transp L\transp a_2 =
-\frac{84}{65}$ gives
\[
  S_{13} = \tfrac{32}{25} + \mu_3\bigl(\tfrac{31}{25} - \mu_1\bigr),
  \qquad
  S_{23} = \tfrac{48}{65} + \mu_3\bigl(\tfrac{84}{65} - \mu_2\bigr),
  \qquad
  S_{33} = 1 - 2\mu_3 - 2\mu_3^2.
\]
Since $S_{11} > 0$, $S_L(\mu)$ has rank at most $1$ if and only if its second
and third columns are multiples of the first, that is, $S_{11}S_{22} =
S_{12}^2$ (automatic on the arc), $S_{11}S_{23} = S_{12}S_{13}$, and
$S_{11}S_{33} = S_{13}^2$. A rank-one symmetric matrix with $S_{11} > 0$ is
positive semidefinite, so these conditions characterise $\mathcal{M}(G_0,3)$
over the arc. Dividing the second condition by $\sigma_1$ and inserting
$\mu_1(t)$, $\mu_2(t)$,
\[
  \sigma_1 S_{23} - \sigma_2 S_{13}
  = \frac{24t}{5(1+t^2)(9t^2+1)}
    \Bigl[(1+t^2)\mu_3 - \frac{256(1-4t^2)}{325}\Bigr],
\]
which is linear in $\mu_3$ with coefficient $\sigma_2 > 0$; hence
\[
  \mu_3 = \mu_3(t) := \frac{256(1-4t^2)}{325(1+t^2)} = \tfrac{16}{13}\mu_1(t),
  \qquad\text{and then}\qquad
  S_{13} = \frac{32(89 + 186t^2 - 543t^4)}{1625(1+t^2)^2}.
\]
Substituting into the third condition and clearing denominators,
\[
  S_{11}S_{33} - S_{13}^2
  = -\frac{256\,(4t^2-1)\,p(t^2)}{1625^2\,(1+t^2)^4},
  \qquad
  p(x) := 294{,}849\,x^3 + 203{,}198\,x^2 - 451{,}015\,x - 31{,}684.
\]
The factor $4t^2 - 1$ recovers $\mu = 0$ at $t = t_0$. By Descartes' rule of
signs, $p$ has exactly one positive root $x_* = t_*^2$, and $p(\frac{1}{9}) <
0$, $p(0.9812) < 0 < p(0.9813)$, and $p(1) = 15{,}348 > 0$ place $t_*^2 \in
(0.9812, 0.9813) \subset (\frac{1}{9}, 1)$. Consequently $\mu^* :=
(\mu_1(t_*),\mu_2(t_*),\mu_3(t_*)) \approx (-0.9449, 0.3296, -1.1629)$ is the
only point of $\mathcal{M}(G_0,3)$ over the admissible arc other than $\mu =
0$, and $G(\mu^*)$ has rank $3$. (In the full space $\mathcal{M}(G_0,3)$,
neither flat-folded endpoint of Item~2 lifts, since $S_{13} = S_{23} = 0$ is
inconsistent at both. The fold-over arc $tu = -3$, on the other hand, carries
two non-immersive lifts, which fail positive path rank at $v_{12}$. They are
determined by $\mu_3 = -\frac{192(9t^2+29)}{325(13t^2+45)}$ together with a
quartic in $t^2$ with sign pattern $(-,+,+,+,-)$ that is positive at $t = 1$.
Thus $|\mathcal{M}(G_0,3)| = 4$ and
$\mathcal{M}_{\mathrm{face}}(G_0,3;b) = \{0,\mu^*\}$.)

It remains to compare $\mu^*$ with the faces of $\arr(0,0,-1)$. This
arrangement has seven chambers (every sign pattern except $\sigma = 001$), nine
edges, and the three vertices $v_{12} = (0,0)$, $v_{13} = (0,1)$, and $v_{23} =
(1,0)$. Every face is incident to a vertex. Positive path rank at every
non-empty face therefore follows from positive path rank on the three vertex
squares of $[0,1]^3$: $\{s_3 = 0\}$ (at $v_{12}$, where $\sigma_3 = 0$),
$\{s_2 = 1\}$ (at $v_{13}$, where $\sigma_2 = 1$), and $\{s_1 = 1\}$ (at
$v_{23}$, where $\sigma_1 = 1$). Let $\widetilde{B} = [L\ \widetilde{A}]$ realise $G(\mu^*)$ in
$\RR^3$, put $\widetilde{J}(s) := L + \widetilde{A}\diag(s)W$, and write
$\widetilde{J}_\sigma := \widetilde{J}(\sigma)$ for the chamber Jacobians. On
$\{s_3 = 0\}$, $\widetilde{J}(s)$ is the Jacobian of the $(2,2,3)$
configuration at $t = t_*$, which has positive path rank by Item~2. For the
other two squares we use the bounded chamber $\triangle := \cell{110}$, which
is incident to all three vertices, and the function
\[
  d(s) := \det\bigl(\Upsilon_{110}\transp G(\mu^*)\Upsilon(s)\bigr)
  = \det\bigl(\widetilde{J}_\triangle\transp\widetilde{J}(s)\bigr),
  \qquad \Upsilon(s) := [I_2;\,\diag(s)W].
\]
If $\rank\widetilde{J}(s) \le 1$ then $d(s) = 0$, so it suffices to show that
$d > 0$ on both squares. Along the admissible arc with $\mu_3 = \mu_3(t)$,
expanding $G(\mu(t)) = G_0 + \sum_{h=1}^3 \mu_h(t) E_h$ gives
\[
\begin{aligned}
  \Upsilon_{110}\transp G(\mu(t))\Upsilon(s)
  &= \begin{pmatrix} \frac{8}{5}\cos\phi(t) & -\frac{6}{5} \\[2pt]
      -\frac{3}{5} & \frac{4}{5}\cos\psi(t) \end{pmatrix}
    + s_1\begin{pmatrix} 4 - \frac{8}{5}\cos\phi(t) \\[2pt] \frac{3}{5} \end{pmatrix} e_1\transp \\
  &\quad + s_2\begin{pmatrix} \frac{6}{5} \\[2pt] 1 - \frac{4}{5}\cos\psi(t) \end{pmatrix} e_2\transp
    + s_3\begin{pmatrix} \frac{32}{25} + \mu_3(t) \\[2pt] \frac{48}{65} + \mu_3(t) \end{pmatrix}(1\ \ 1).
\end{aligned}
\]
On any vertex square with free coordinates $(s_a,s_b)$, this has the form $K_0
+ s_a k_a w_a\transp + s_b k_b w_b\transp$, and the coefficients of $s_a^2$ and
$s_b^2$ in its determinant are the determinants of the rank-one matrices $k_a
w_a\transp$ and $k_b w_b\transp$, which vanish. Hence $d$ is bilinear in
$(s_a,s_b)$, equals the bilinear interpolation of its corner values $d_\sigma :=
\det(\widetilde{J}_\triangle\transp\widetilde{J}_\sigma)$, and is positive on
the square as soon as these four values are. We evaluate the $2\times 2$
determinants directly from the formula above; no square roots are involved.
This gives $d_{110} = \det\diag(4,1) = 4$ and
\begin{gather*}
  d_{010} = \frac{8(1-t^2)}{5(1+t^2)},\qquad
  d_{100} = \frac{16(9t^2-1)}{5(9t^2+1)},\qquad
  d_{111} = \frac{4(989 - 611t^2)}{325(1+t^2)},\\
  d_{011} = \frac{8(1493 - 1216t^2 - 149t^4)}{1625(1+t^2)^2},\qquad
  d_{101} = \frac{16(313 + 10{,}244t^2 - 9{,}909t^4)}{1625(1+t^2)(9t^2+1)},
\end{gather*}
all of which are positive for $\frac{1}{3} < t < 1$: the two quadratics in
$t^2$ are positive at $t^2 = 0$ and $t^2 = 1$, the first being decreasing and
the second concave. Since $t_* \in (\frac{1}{3},1)$, the corner chambers
$\{010, 011, 110, 111\}$ of $\{s_2 = 1\}$ and $\{100, 101, 110, 111\}$ of
$\{s_1 = 1\}$ all have $d_\sigma > 0$. This proves $\mu^* \in
\mathcal{M}_{\mathrm{face}}(G_0,3;b)$. The only small margin is $d_{010} =
\frac{8}{5}\cos\phi(t_*) \approx 0.015$, which vanishes exactly when the crease
ray $c_1$ becomes orthogonal to $e_1$ at $t = 1$.

Finally, we check that $\mu^* \notin \mathcal{M}_{\mathrm{PR}}(G_0,3)$. Realise
$G(\mu^*)$ in gauge $L = [e_1\ e_2]$ by $\tilde{a}_1 = c_1(t_*) - e_1$,
$\tilde{a}_2 = c_2(t_*) - e_2$, and $\tilde{a}_3 = (\mu_3(t_*),\mu_3(t_*),
\tilde{a}_{3,3})\transp$, where $\tilde{a}_{3,3} := S_{13}(t_*)/\sigma_1(t_*)$
and $x_* := t_*^2 \in (0.9812,0.9813)$. With $q(x) := 543x^2 - 186x - 89$ we
get $q(x_*) \in (251,252) > 0$ and $t_*\tilde{a}_{3,3} =
-\frac{2q(x_*)}{325(1+x_*)} < 0$, together with the identity $\tilde{a}_1 = e_1
+ \kappa\tilde{a}_2$ for $\kappa := \frac{2(9x_*+1)}{3(1+x_*)} > 0$.
Collinearity $c^{(2)}(s) = \lambda c^{(1)}(s)$ of the columns of
$\widetilde{J}(s)$ reads
$-\lambda(1+s_1)e_1 + e_2 + (s_2 - \lambda\kappa s_1)\tilde{a}_2 +
(1-\lambda)s_3\tilde{a}_3 = 0$. Rows~$2$ and~$3$ fix $\eta := s_2 -
\lambda\kappa s_1$ and $\zeta := (1-\lambda)s_3$ over $\mathbb{Q}(x_*)$ as
\[
  \eta = \frac{5(9x_*+1)\,q(x_*)}{D(x_*)},
  \qquad
  \zeta = \frac{3900\,x_*(1+x_*)}{D(x_*)},
\]
where $D(x_*) := 3072\,x_*(4x_*-1) + 9(x_*+1)q(x_*) > 8790 + 4470 = 13{,}260$.
As $5(9x_*+1)q(x_*) < 50\cdot 252 = 12{,}600$ and $3900x_*(1+x_*) < 7800$, we
have $\eta,\zeta \in (0,1)$. Row~$1$ then sets $\lambda(1+s_1) =
-\frac{3}{5}\eta + \zeta\mu_3(t_*) =: \lambda_0$. Because $3(9x_*+1) >
9(x_*+1)$, this gives
$-\lambda_0 = (3072\,x_*(4x_*-1) + 3(9x_*+1)q(x_*))/D(x_*) > 1$. Moreover
$\lambda_1 := \lambda_0 + \eta/\kappa =
-(3072\,x_*(4x_*-1) + \frac{3}{2}(13x_*-3)q(x_*))/D(x_*) \in (\lambda_0,0)$,
and $\lambda_0/\lambda_1 - 1 = (\eta/\kappa)/(-\lambda_1) \in (0,1)$, as
$\eta/\kappa = \frac{15}{2}(1+x_*)q(x_*)/D(x_*) < 3800/D(x_*) < -\lambda_1$.
Thus $\rank\widetilde{J}(s(\lambda)) = 1$ for $\lambda \in
(\lambda_0,\lambda_1)$ along the curve
$s(\lambda) := (\lambda_0/\lambda - 1,\, \eta + \kappa(\lambda_0-\lambda),\,
\zeta/(1-\lambda)) \in (0,1)^3$. This curve joins the unrealised faces
$\{s_1 = 0\}$ and $\{s_2 = 0\}$ around $\cell{001}$, running from
$(0,\eta,\zeta/(1-\lambda_0)) \approx (0, 0.9290, 0.2568)$ to
$(\lambda_0/\lambda_1-1,0,\zeta/(1-\lambda_1)) \approx (0.2989, 0, 0.2939)$.
Together with the reduction to the admissible arc above, this gives
$\mathcal{M}_{\mathrm{PR}}(G_0,3) = \{0\}$. At $b_3 = 0$,
Remark~\ref{rem:q2-counterexamples} gives $\mathcal{N}(\arr) =
\spn\{E_1,E_2,E_3\}$, and local rigidity modulo $O(3)$ follows from the
injectivity of $\Lambda$ as before.
\end{proof}
\PaperAppendicesEnd

% Bibliography for Paper A v3.

\end{document}